# A tour d'horizon of de Casteljau's work

Andreas Müller

*Hochschule Kempten, Bahnhofstraße 61, 87435 Kempten, Germany*

ABSTRACT

Whilst Paul de Casteljau is now famous for his fundamental algorithm of curve and surface approximation, little is known about his other findings. This article offers an insight into his results in geometry, algebra and number theory.

Related to geometry, his classical algorithm is reviewed as an index reduction of a polar form. This idea is used to show de Casteljau's algebraic way of smoothing, which long went unnoticed. We will also see an analytic polar form and its use in finding the intersection of two curves. The article summarises unpublished material on metric geometry. It includes theoretical advances, e.g., the 14-point strophoid or a way to link Apollonian circles with confocal conics, and also practical applications such as a recurrence for conjugate mirrors in geometric optics. A view on regular polygons leads to an approximation of their diagonals by golden matrices, a generalisation of the golden ratio.

Relevant algebraic findings include matrix quaternions (and anti-quaternions) and their link with Lorentz' equations. De Casteljau generalised the Euclidean algorithm and developed an automated method for approximating the roots of a class of polynomial equations. His contributions to number theory not only include aspects on the sum of four squares as in quaternions, but also a view on a particular sum of three cubes. After a review of a complete quadrilateral in a heptagon and its angles, the paper concludes with a summary of de Casteljau's key achievements.

The article contains a comprehensive bibliography of de Casteljau's works, including previously unpublished material.

**Contents**

*E-mail address:* andreas.mueller@hs-kempten.de.

## Introduction

Paul de Faget de Casteljau (Fig. 1) is widely recognised today for his groundbreaking algorithm that enables the generation of polynomial curves and surfaces. However, his additional contributions to geometry, algebra or number theory remain less explored. Despite the limited accessibility of many of his writings, this article offers a comprehensive overview of his work, drawing closely from his articles, books, and personal communications.

This paper is divided into 19 sections and includes 30 figures. It is organised into four main parts: Introduction (1), Geometry (2-13), Algebra and Number Theory (14-18), and Conclusion (19). As we delve into the various sections of this paper, we will explore the central themes and methods present in de Casteljau's work. We will recall his concept of 'simple poles' and 'generalised poles' also known as blossoms or polar forms. We will learn about the close link between interpolation and approximation in blossoming, which de Casteljau represented by the characteristics $(n, c, r)$ of degree $n$, continuity $c$, and restitution $r$. Furthermore, we will see his approaches to tolerances, smoothing, and interpolation. In homage to his brother, de Casteljau worked on the structure of foundational theorems in projective geometry, a topic we will revisit in this article. His later works shifted focus towards metric geometry. We will uncover a previously unpublished introduction to the subject, see applications in geometric optics, gain insights into strophoids using bifocal coordinates, and explore de Casteljau's theorem on golden ratios derived from regular polygons. De Casteljau's engagement with algebra and number theory encompasses matrix quaternions, a generalised Euclidean algorithm in $n$ variables and an idea on the sum of three cubes.

Throughout the sections, addenda serve to spotlight divergent ideas and alternative pathways stemming from de Casteljau's contributions. These serve partly as commentary and partly as avenues for opening new perspectives or interpretations. Similar to the core ideas, they are rooted in the diverse writings of de Casteljau.

The overarching goal of this article is to give room to his various findings, and to pay homage to a mathematician whose innovative ideas often lacked widespread recognition. Before delving into his original ideas, let us briefly trace his life and career.

## 1. Mindset and career

> *No one in 1958 could have foreseen that one day mathematicians would be looking for the equation of an artistic or aesthetic form. Until then, they limited themselves to drawing the curves representing the equations, while being fascinated, above all, by the singularities. Now the singularity has become an obstacle. Without the insistence of Monsieur de la Boixière, and my reluctance to leave Citroën without having been able to offer them a tangible solution, what would have been the history of CAD?* — de Casteljau (Müller, 2024, p.42)

Paul de Faget de Casteljau's life took a fateful turn when German troops invaded France in 1940.[1] Born on 19 November 1930, he was just nine years old when his mother decided to leave their hometown Besançon for the small château near Agen on the banks of the Garonne river, where her parents had lived during her childhood. Her father had been in the French Navy all his life: he had studied at the prestigious *École Navale*, was involved in the Crimean war as a young cadet, and left the fleet with the rank of Rear-Admiral. In his library, the young Paul de Casteljau found a broad collection of scientific books, on one side the chemistry books that the grandfather of his grandfather wrote as an early founder of higher education institutions during Robespierre's government; on the other, the literature that his grandfather had used to prepare for the entry exams at the *École Navale*. He enjoyed reading the higher mathematics titles (partly without really understanding, as he said) when he was sent out on Sundays or holidays to look after the geese or sheep in the meadow. Without a pencil, he had to visualise everything in his head and was content to think hard. He would later say, that the language that he learned intuitively, was Mathematics, not French nor German.

French mathematics in the early 19th century, and in particular geometry, was characterised by the works of "the three L's", Lagrange (1736-1813), Laplace (1749-1827), and Legendre (1752-1833), as well as by Monge (1746-1818), Fourier (1768-1830),

[1] In personal letters, de Casteljau mockingly reflects that his fate was based on *"the 'mistake' of Adolf who had decided to do 'Panzer' sightseeing as our roads were too bad"* (personal communication to W. Boehm, Christmas 2010, and to the author, 20 March 2002).

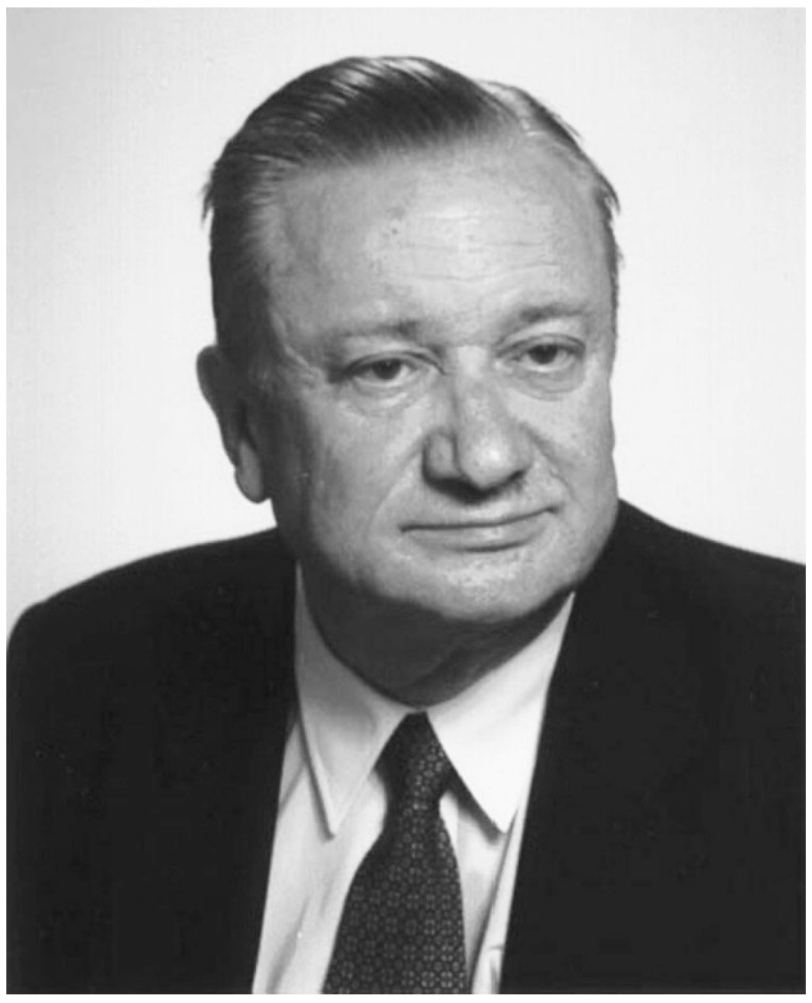

**Fig. 1.** Paul Bertrand René Lucien de Faget de Casteljau (1930-2022).

Cauchy (1789-1857) and Poncelet (1790-1867).[2] It was influenced by the Swiss Euler (1707-1783), the Germans Leibniz (1646-1716) and Gauss (1777-1855), by the Norwegian Abel (1802-1829) and the Irish Hamilton (1805-1865). Geometry was mainly Euclidean or projective (called descriptive) and it was the blooming period of the synthetic geometry with the Swiss Steiner (1796-1863), the German von Staudt (1798-1867), and the Italian Mascheroni (1750-1800). Later, that 19th century would see an immense growth in geometry, stemming in particular from Riemann (1826-1866), Klein (1849-1925) and Poincaré (1854-1912), including the development of non-Euclidean geometry.[3] Assuming that his grandfather's library was not kept up-to-date after his entry to the Navy, the geometry that left a mark on Paul de Casteljau as a ten year-old, was based on the literature of the early 19th century, thus less "arithmetised" (Klein, 1895) or algebraic than we know it today. Continued fractions were in fashion, similarly logarithms and Delambre's formulas for spherical triangles.

When the war ended and his family moved back to Besançon, Paul finished his baccalaureate and entered the preparation classes in his Catholic college. These classes, common in France, prepare students for the challenging entrance exams of the top universities. Both, the texts from his meadow readings as well as the maths preparation classes equipped him well to pass the national entry exams to two of the best French universities, *École Polytechnique (X)* and *École Normale Supérieure (ENS)*, of which he chose the latter. Amongst the 25 students in his scientific class was the later Nobel Prize winner Pierre-Gilles de Gennes together with many others who went on to become Physics professors (e.g., André Authiers, Philippe Nozières, Roland Omnès).[4] The literature track of his year hosted students such as philosopher Michel Serres whom he knew from his school in Agen, and sociologist Michel Bourdieu. The conversations on campus in *Rue d'Ulm* or in the *ENS* dormitories must have offered plenty of inspiration.

The preparation classes already included much of 'modern' mathematics, as his teachers Bhenmerah and Pham Tinh Quat were well acquainted with the approaches of *Nicolas Bourbaki*, the group of young French mathematicians, who published structuralistic textbooks under this invented single name. Yet at *ENS*, he was confronted even with its founding father, Henri Cartan. Many students of his class did not become friends with this abstract approach and decided to focus more on physics than on (modern) maths (Plévert, 2011), as did de Casteljau.

In his letters, he points out some of his preferred literature, such as the *Traité de Géométrie Analytique* by Salmon (1884), *Géométrie* by Deltheil and Caire (1950), or the unfinished *Courbes Géométriques remarquables* by Brocard and Lemoyne (1919). He underlines in particular *Les Géométries* by the Belgian Godeaux (1947) due to his approach to the geometric locus and the choice to inspire the reader instead of enumerate theorems.[5]

Three years of military service in France and Algeria interrupted his stay as a research assistant in the lab of the solid-state physicist André Guinier, where he also published his first article (de Casteljau and Friedel, 1956). Yet after his return, the position was occupied and he decided to apply for a job in the private sector. One job ad was anonymously looking for 'young people presenting well' and only by way of the reply did he learn that he had applied to Citroën.

On 1 April 1958, de Casteljau started at Citroën in the Tooling department, where he pursued the challenge of turning the handmade design models into a mathematical form that would be able to produce a tool milled by a numerically controlled machine.

[2] A more elaborated view can be found in the recommended lectures of Felix Klein (1926–1927) on the development of mathematics in the 19th century.

[3] The first approaches on non-Euclidean geometry were published by Lobachevsky in 1829 and Bolyai in 1831, yet it was largely discussed after Riemann presented his geometry in 1854, cf. (Klein, 1928).

[4] De Gennes' biography describes well the atmosphere with his classmates: "The school was their whole life: they ate together in the *pot* (the refectory) morning, noon and evening, worked together in the *turnes* (study rooms for four or five students) and slept in the dormitory (they were only entitled to their own room after the second year)." (Plévert, 2011).

[5] de Casteljau, personal communication, 20 April 1995.

His initial approach, documented at the French patent office INPI only months after his entry (de Casteljau, 1959), focused on the trajectory of the mill as a 3D-polyline. His famous algorithm allowed from early on to recursively calculate polynomial curves or triangular surfaces within the boundaries of their control points, or *pôles* as they were called at Citroën. This property is connected to what we today call the convex hull. His approach was later integrated into *SADUSCA*,[6] the internal computer system to design and machine Citroën's bodywork (Krautter and Parizot, 1971). Less known from these days is his approach to trigonometric functions and a recursive scheme for their calculation.

What would later be called a polynomial curve in Bézier form,[7] was at Citroën termed *courbe à pôles* (pole curve), a name invented by his department head Jean de la Boixière, in reference to *courbe à points* (point curve), which described the digitised (or at least sampled) result of the freeform curves and surfaces, made by the hands of respected designers.[8] Soon after the first work report, de Casteljau extended his approach to tensor product surfaces, and the intellectual property for the curves and surfaces *à pôles* was registered with the bailiff in 1963, coinciding with the year of his wedding. The practice of designing curves and surfaces with the help of control vertices was also taught in Citroën's design school since the early 1960s.[9]

At about the same time and not far away, the engineer Pierre Bézier (1910-1999) at Renault was also working on the question of how to mathematically describe the curves given by designers, such that the tools could be formed in an exact way, based on NC milling. Bézier had joined Renault in 1933, and after a first management success in production method development, he decided to work in technical development and particularly on numerical control. As we know today, he had knowledge of de Casteljau's *courbes à pôles* approach at large by some defectors from Citroën to Renault, yet he was not in possession of the precise theoretical background (Gardan, 2023). So Bézier came up with his own way of constructing those curves. His approach by intersecting quarter cylinders differs from de Casteljau's 'barycentric method' and it is even more different in the case of surfaces: to construct triangular surfaces, Bézier coalesces two vertices (which poses some challenges regarding continuity or surface degree), while de Casteljau built triangular surfaces based on barycentric coordinates from the beginning.[10] However, due to different corporate policies on communicating intellectual property, Bézier already gained internationally prominence in the late 1960's (Bézier, 1966, 1967, 1968); while de Casteljau was first mentioned in France in the 1970's (Krautter and Parizot, 1971; Parizot, 1971–1975) (cf. Fig. 3), and was only recognised internationally from 1977 after being discovered by Boehm in 1975[11] (Böhm and Gose, 1977; Boehm, 1977, 1981; Boehm et al., 1984). Citroën's secrecy also had an influence on knowledge networks: while *Unisurf*, the CAD/CAM system of Renault and Peugeot was in part developed in collaboration with French universities, Citroën's *SADUSCA* depended on the know-how of a small group, not larger than ten persons including de Casteljau, Bouffard-Vercelli, Deschamps, Ménard, Loschutz, Heinrich, Sureau, Giaume. Today, we honour both CAD pioneers with the notions of the **Bézier form** of a polynomial curve or surface and the **de Casteljau algorithm**.

Citroën was renowned for its innovative prowess, with the founder André Citroën already having invested more in technology than he could afford from the revenues of selling cars, which led to its first takeover, by Michelin in 1934. The chevrons maintained their innovative image with cars such as the *Traction Avant* (1934-1957), *DS* (1955-1975), *SM* (1970-1975), *GS* (1970-1984), *CX* (1974-1991) and the commercial *Type H* (1948-1981), which all became icons of French car design and technology. Yet the high investments were again challenged by the markets when the Oil Crisis led customers to demand more affordable and mass-produced cars. When Citroën got into financial trouble, Peugeot first had to acquire 38% of its shares in 1974, and ultimately assumed full control of the company in 1976, the whole deal being enforced by de Gaulle's government. The unwanted merger imposed various frictions for employees as well as for management, including de Casteljau.

The early years at Citroën had invited competition on the brightest solutions, where a gifted mathematician like de Casteljau could be seen as a competitive advantage, in particular in combination with a talented implementation team. Yet these new ideas sometimes posed a challenge for the existing processes. The price for a successful idea often lay in derision or exclusion. His department waited for opportunities along the lines of a strike by wood modellers or seemingly unsolvable challenges such as a more efficient cam design in order to prove their value to the organisation. They could combine their mathematical strength with a growing understanding of the designer challenges and could eventually build a practical example, which seemingly was the only way to convince in the automotive industry at that time (de Faget de Casteljau, 1999a; Müller, 2024).[12]

It was during this pioneering period that de Casteljau not only laid all the groundwork for the 'mathematisation' of curves and surfaces – the early days of digital transformation (Carpo, 2017, p.55ff.), (Smith, 2021, p.270ff.), but he also improved the design of cam shafts (de Faget de Casteljau, 1995a) and made his mark in tooling processes (de Faget de Casteljau, 1997), while promoting numerically stable algorithms. His colleague Bouffard-Vercelli put together the conversational software SPAC/CAR (*Système de Programmation Automatique Citroën / CARrosserie*) with its versions CAR-C (*Catalogue*), CAR-G (*Géométrie*) and CAR-U (*Usinage*

[6] *Système d'Aide à la Définition et à l'Usinage des Surfaces de Carrosserie*: Assistance System for the Definition and Machining of Body Surfaces.

[7] Or dCB (de Casteljau Bézier) curve, as Wildberger (2021) salomonically calls the curves, as they are so much connected to both, Paul de Casteljau as well as Pierre Bézier.

[8] The term *pôles* also referred to its occurrence in POLynomials, POLar form, interPOLation, cf. (de Faget de Casteljau, 1992a).

[9] A former Citroën designer recalls (translated from French): "In 1961, André Estaque [head of bodywork design] enrolled me in a course on the mathematisation of bodywork: *les courbes à pôles* [the pole curves]. These curves are supported by external points that originate from tangents. In a way, it's the image of a marquee whose interior volume is supported by exterior support points. By acting on these points, you can vary the volume ... " (Nicolas, 2008).

[10] Further similarities and differences of the two CAD pioneers are described well by Piegl (1993); Boehm and Müller (1999); Laurent and Sablonnière (2001); Lagrange (2008).

[11] Communication from W. Boehm to Citroën, dated 4 September 4 1975, response by J. Krautter and S. Parizot on 15 September 1975, followed by a reply from W. Boehm on 29 September 1975.

[12] Pierre Bézier experienced similar resistance in Renault, cf. (Bézier, 1982a,b; Rabut, 2002)

- machining) (Krautter and Parizot, 1971). The software was far from being interactive, and the calculations were performed on mainframe systems with only a fraction of today's computing power. These early systems were able to create wireframe models of selected parts (cf. Fig. 2) and were capable of driving NC milling machines to generate the tools that later stamped metal sheets into car bodies.

After a heyday of in-house CAD/CAM software development in the early 1970s, a growing software industry led all carmakers towards the acquisition of ready-made software, written and licensed by the large global CAD/CAM providers such as Dassault or Computervision (later acquired by PTC) (Urban-Galindo and Ripalles, 2019). In-house mathematicians or programmers were not as much needed by carmakers as before, and de Casteljau even received a formal note not to work on CAD topics anymore.

And yet, on Pentecost 1979,[13] he had an idea how to combine his curves with the smoothness conditions of spline curves. Different from the *pôles simples* so far in use, he called this approach the *pôles généralisées*, which Ramshaw (1987, 1989) later introduced as *blossoms*. But due to his CAD ban, de Casteljau did not talk about it. Only upon the occasion of a possible concentration of all French CAD knowledge into one French CAD company did he mention his idea as a potential business, and would subsequently present it to the Peugeot-Citroën management. From then on, the CAD ban was lifted, and some 25 years after his famous algorithmic idea, de Casteljau was invited to publish his ideas with a tailwind from Citroën's CEO, Xavier Karcher, who wrote the foreword of his first book.

The idea to unite existing know-how and experience under the roof of an all-French CAD house did not materialise. Instead Peugeot, which still had a cooperation with Renault and was involved in the development of Unisurf 3, voted to acquire CAD licences from the external vendor Computervision (CV),[14] although their main product CADDS3 lacked the capabilities to model free-form curves and surfaces in 3D. Peugeot SA (PSA) decided to support the company with its proprietary SPAC knowledge, with some PSA engineers joining Computervision and some PSA employees enhancing the software with their existing approaches to form the new system "CV+". A re-writing of the software turned the so-far conversational cathode-ray tube system into an interactive CAD/CAM system on a graphics workstation, which entered the market as CADDS4 in the early 1980s.

De Casteljau did not play a role in the Computervision activities, but instead acted as a scientific advisor, warmly addressed as *'Maître'* by his colleagues. As a by-product, his interests turned into non-CAD topics such as quaternions or other problems from Euclidean metric geometry. To his surprise, his work on quaternions proved to be more CAD-related than initially thought (applied to programming the kinetics of robots), and his ideas were in demand once again at Citroën.

When de Casteljau left Citroën for his retirement, he had published three books (de Casteljau, 1985, 1987; de Faget de Casteljau, 1990) and had become visible to the research community. His first invitation to a CAGD conference was to Biri (Norway) in 1991 (de Faget de Casteljau, 1992a), after which several other conferences and articles followed (de Casteljau, 1993; de Faget de Casteljau, 1994, 1995a, 1997; de Casteljau, 1998; de Faget de Casteljau, 1999b, 2000, 2001a,b). During latter decades, de Casteljau corresponded with CAGD researchers through handwritten letters or greeting cards, which often included tiny drawings of his current ideas. He has bequeathed two autobiographies, one short and rather coded version from 1992 when he left Citroën (de Faget de Casteljau, 1992b; Müller, 1995; de Faget de Casteljau, 1999a), and one longer version written in 1997 (Müller, 2024).

Paul de Faget de Casteljau passed away on 24 March 2022, in Versailles (cf. Table 1). He was a member of the *Académie des Sciences, Lettres et Arts de l'Ardèche*, and he came to honour beyond the naming of his fundamental algorithm, receiving the 1987 Seymour Cray Prize from the French National Centre for Scientific Research, the John Gregory Memorial Award (awarded in Dagstuhl 1996), being honoured as a pioneer in Computer Graphics by SIGGRAPH in Orlando in 1998, and receiving the 2012 Bézier Award from the Solid Modeling Association. In 1997, he received a *Doctor honoris causa* from the University of Berne, Switzerland, which honoured his *'fundamental scientific discoveries and inventions in the field of geometric modelling as well as his ideas with which he paved the way for computer supported construction and production."* (Bieri and Prautzsch, 1999)

De Casteljau was proud of his ability to reduce his texts to a minimum, such that the results were like a "summarised memory aid", as one of his former teachers called this early on. To de Casteljau, this ability valued the essence of a topic, *"like a cut diamond, once gotten rid of its surrounding gangue. In geometry, this consists of extracting from a figure only the fundamental features, such that one cannot remove one nor add another, without deflating the figure"* (de Faget de Casteljau, 1999b, p.110). Regrettably, most of his publications were not translated into English, nor were they digitised or are easily accessible. Additionally, several of his ideas were never formally published. This makes it more challenging to understand his works or to build on his ideas. Hence, this article seeks to provide an overview of his main works, inviting the reader to unearth more 'diamonds' in the quoted original sources.

## On geometry

This second part will cover the geometrical aspects: the initial algorithm (including focal splines) and its extension to polar forms, algebraic smoothing as a blend of interpolation and approximation, and an integrated view on core theorems in projective geometry. Topics in metric geometry include a link between Apollonian circles and confocal conics, some ideas on geometric optics, his 14-point strophoid (related to Feuerbach's nine-point circle), and a remarkable recursive scheme for regular polygons.

[13] de Casteljau, personal communication to W. Boehm, 20 September 1992.

[14] At that time, Computervision was still managed by its founders, the Paris-born Philippe Villers, educated at Harvard and MIT, and the entrepreneur Martin Allen, a graduate from the University of California. Computervision dominated the CAD industry at the time with a good third of the market share. The company hired talented personalities like Ken Versprille, the father of NURBS, Sam Geisberg, who later started PTC, or his brother Vladimir Geisberg, later CAD VP at Prime Computer. Prime acquired Computervision in 1988, and was itself acquired by PTC in 1998. Prime/Computervision had invested in the product data firm Windchill Technology, PTC bought the remaining shares and integrated its personnel. Windchill founder Jim Heppelmann today is PTC's CEO.

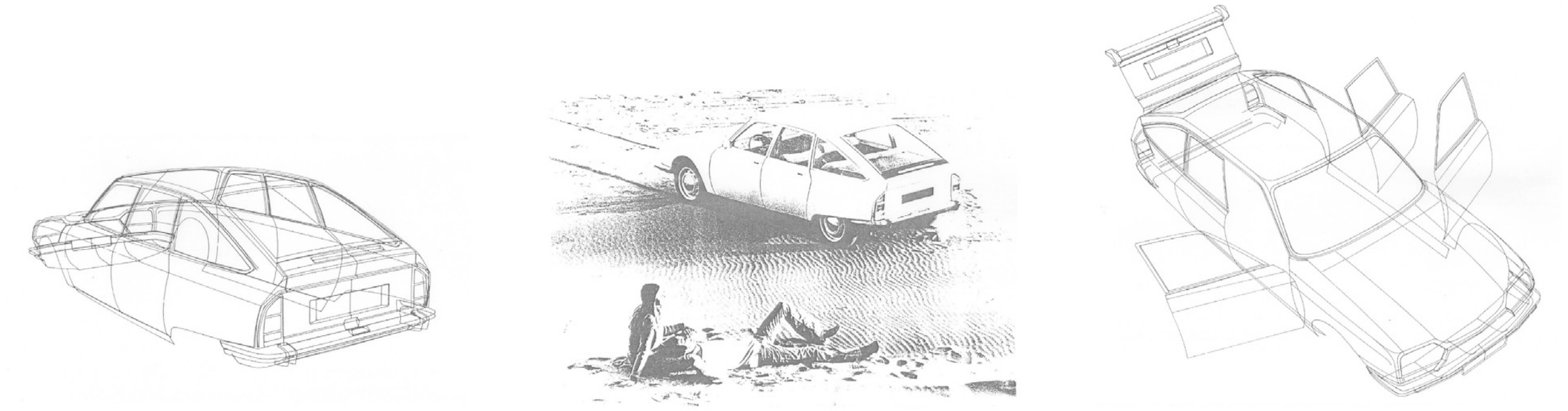

**Fig. 2.** Wireframe models and image of the GS, the first Citroën model with a complete mathematical description of its body (de Casteljau, personal communication, 29 September 1998).

**Table 1**
A timeline featuring significant events from Paul de Casteljau's life.

| Year | Event |
|---|---|
| 1930 | Born on 19 November in Besançon, Franche-Comté, *École maternelle, Collège St Jean* |
| 1940 | Flees with his mother to the place of his grandfather in Agen, *Collège St Caprais* |
| 1945 | Family returns to Besançon, *Collège St Jean, Lycée Victor Hugo* |
| 1951-1955 | Student at the *École Normale Supérieure, Rue d'Ulm*, graduating as *Licencié ès Sciences* |
| 1956 | His very first publication, on Noble Metals (de Casteljau and Friedel, 1956) |
| 1955-1958 | Military service in France and Algeria |
| 1958 | Starts at Citroën on 1 April |
| 1959 | Report *Enveloppe Soleau 40.040* is registered at INPI Paris by Citroën |
| 1958-1960 | First machined part based on de Casteljau's algorithm (Poitou, 1989) |
| 1963 | Notes on *Courbes et Surfaces à Pôles* are registered with the bailiff by Citroën |
| 1971 | Krautter and Parizot (1971) publish the Citroën approach at a SIA conference |
| 1974 | First mention in French in *Enyclopédie des Métiers* (Parizot, 1971–1975) |
| 1975 | Boehm discovers de Casteljau |
| 1977 | First mention in German by Boehm (1977) |
| 1981 | First mention in English by Boehm (1981) |
| 1985 | First book *Formes à Pôles* |
| 1987 | Book *Les Quaternions* |
| 1990 | Book *Le Lissage* |
| 1991 | First conference appearance in Biri, Norway |
| 1992-2001 | Publishes ten articles, partly in hardly accessible conference proceedings, primarily in French |
| 1997 | *Doctor Philosophiae honoris causa* from University of Berne, Switzerland |
| 2022 | Passes away on 24 March in Versailles |

## 2. The fundamental calculus

*Instead of drawing the curve that represented an equation, we had to find the equation that defined a form!* — de Casteljau[15]

The internal report, which de Casteljau wrote in the winter of 1958/59, was registered as an *Enveloppe Soleau* in March 1959 at the *Institut Nationale de la Propriété Industrielle (INPI)*, thus providing evidence of an early instance of his idea. The report contains some remarkable approaches beyond of course the original algorithm:

- **barycentric coordinates:** Differently from the Bernstein notation today, de Casteljau worked from the beginning with barycentric coordinates $p+q=1$. He even headlines his approach *Méthode Barycentrique*. His natural generalisation on surfaces in 1959 needs one more coordinate and thus describes triangular surfaces.
- **recurrence relation:** For the sake of effective storage management, the approach was meant to be recursive. As the initial paper works with the parameter $p=q=1/2$, we could also speak of a subdivision scheme (Boehm and Müller, 1999, p.590).
- **difference scheme:** The evaluation of a point is based on the well-known triangular scheme of the control points. Yet, he not only re-iterates the affine combinations, but also works with additions of finite differences to gain $C^k$ continuity as Stärk (1976) later did independently.
- **poles:** The first report already introduces the notion of *pôles* for the control points – a name given by his superior de la Boixière, which de Casteljau will use throughout all further documents (e.g., *Courbes et Surfaces à Pôles* (de Casteljau, 1963)).

[15] de Casteljau, personal communication to C. Rabut (translated from French), December 1999.

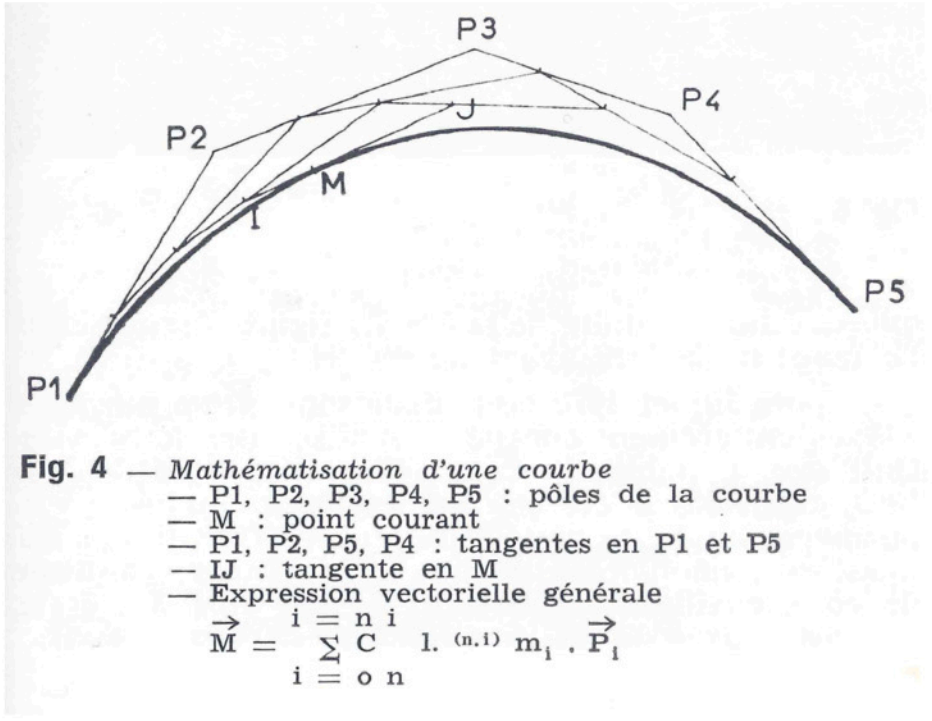

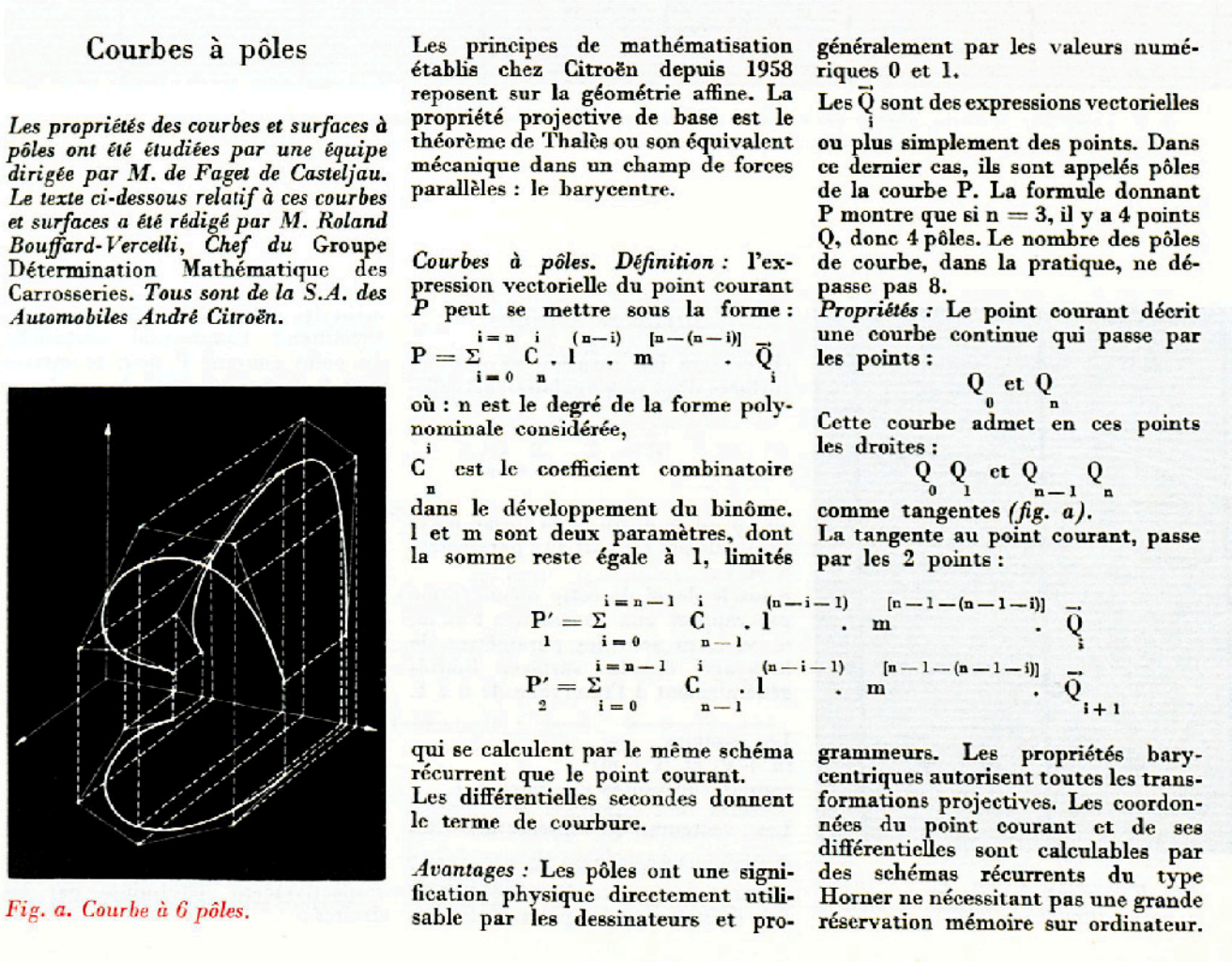

Courbes à pôles

*Les propriétés des courbes et surfaces à pôles ont été étudiées par une équipe dirigée par M. de Faget de Casteljau. Le texte ci-dessous relatif à ces courbes et surfaces a été rédigé par M. Roland Bouffard-Vercelli, Chef du Groupe* Détermination Mathématique des Carrosseries. *Tous sont de la S.A. des Automobiles André Citroën.*

*Fig. a. Courbe à 6 pôles.*

Les principes de mathématisation établis chez Citroën depuis 1958 reposent sur la géométrie affine. La propriété projective de base est le théorème de Thalès ou son équivalent mécanique dans un champ de forces parallèles : le barycentre.

*Courbes à pôles. Définition :* l'expression vectorielle du point courant $P$ peut se mettre sous la forme :

$$P = \sum_{i=0}^{i=n} C^{i}_{n} \,.\, 1^{(n-i)} \,.\, m^{[n-(n-i)]} \,.\, \vec{Q}_i$$

où : n est le degré de la forme polynominale considérée,

$C^{i}_{n}$ est le coefficient combinatoire dans le développement du binôme. l et m sont deux paramètres, dont la somme reste égale à 1, limités généralement par les valeurs numériques 0 et 1.

Les $\vec{Q}_i$ sont des expressions vectorielles ou plus simplement des points. Dans ce dernier cas, ils sont appelés pôles de la courbe P. La formule donnant P montre que si n = 3, il y a 4 points Q, donc 4 pôles. Le nombre des pôles de courbe, dans la pratique, ne dépasse pas 8.

*Propriétés :* Le point courant décrit une courbe continue qui passe par les points :

$$Q_0 \text{ et } Q_n$$

Cette courbe admet en ces points les droites :

$$Q_0\,Q_1 \text{ et } Q_{n-1}\,Q_n$$

comme tangentes *(fig. a)*.

La tangente au point courant, passe par les 2 points :

$$P'_1 = \sum_{i=0}^{i=n-1} C^{i}_{n-1} \,.\, 1^{(n-i-1)} \,.\, m^{[n-1-(n-1-i)]} \,.\, \vec{Q}_i$$

$$P'_2 = \sum_{i=0}^{i=n-1} C^{i}_{n-1} \,.\, 1^{(n-i-1)} \,.\, m^{[n-1-(n-1-i)]} \,.\, \vec{Q}_{i+1}$$

qui se calculent par le même schéma récurrent que le point courant.

Les différentielles secondes donnent le terme de courbure.

*Avantages :* Les pôles ont une signification physique directement utilisable par les dessinateurs et programmeurs. Les propriétés barycentriques autorisent toutes les transformations projectives. Les coordonnées du point courant et de ses différentielles sont calculables par des schémas récurrents du type Horner ne nécessitant pas une grande réservation mémoire sur ordinateur.

**Fig. 3.** Left: First public presentation of the algorithm by Citroën at the *Salon de l'Automobile* 1971 in Paris (without reference to de Casteljau) (Krautter and Parizot, 1971, p.583). Right: First public reference to de Casteljau in the *Encyclopédie des Métiers*, probably in 1974 (Bouffard-Vercelli, 1971–1975, p.257).

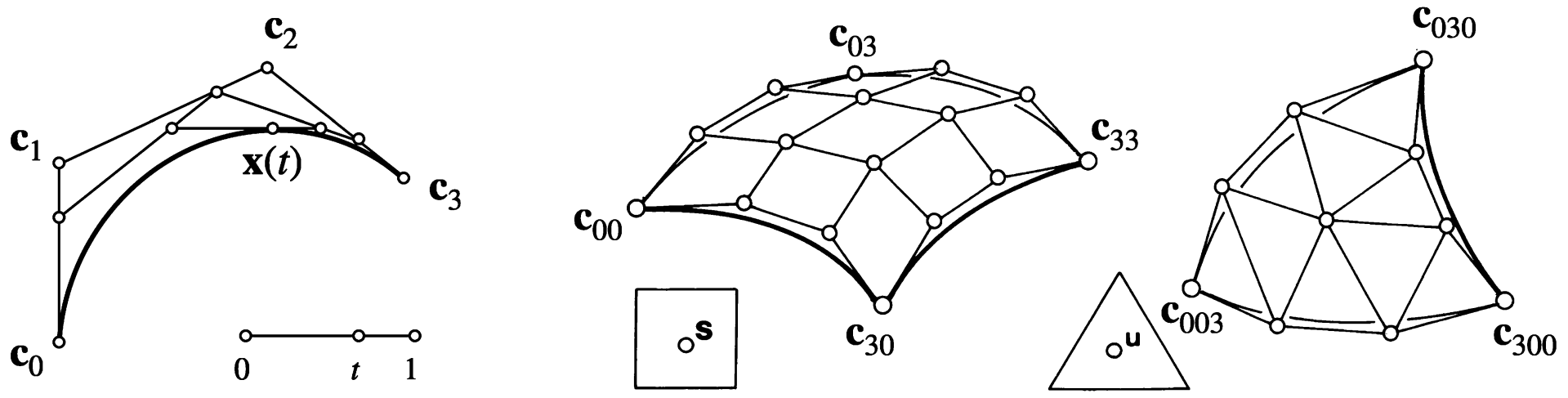

**Fig. 4.** The de Casteljau's algorithm, $n = 3$: barycentric combinations of two control points define the next family of control points, the final point lies on the curve $\mathbf{x}(t)$; the same applies to triangular surfaces with three control points each, and to quadrangular surfaces with four control points each, images adapted from (Boehm and Müller, 1999).

- **periodic functions:** The application of central differences to get a mainly additive calculation scheme is not limited to polynomials but also extends to periodic functions: Executing the angle addition theorem on $F_{n\pm 1} = \sin(n\varphi \pm \varphi)$, we get $F_{n+1} + F_{n-1} = 2F_n \cos\varphi$, and from here

$$\Delta_2 F_{n-1} = F_{n+1} - 2F_n + F_{n-1} = 2F_n(\cos\varphi - 1) = -4\sin^2\frac{\varphi}{2} F_n = -kF_n \tag{1}$$

  Knowing $F_0 = 0$ and $\Delta_1 = F_1 - F_0 = \sin\varphi$, we can generate a trigonometric table with constant step size and given precision.

More observations and details can be taken from Rabut (2024) in this special issue.

## 3. "The" de Casteljau algorithm

*This algorithm is not only interesting in the area of CAD, but also as an example of recurrence. From my point of view, it is of great elegance with a very telling geometric interpretation. In general, I started telling students that this algorithm was as beautiful as a poem by Rimbaud.* — Gardan (2023)

Within just one year of starting at Citroën, de Casteljau documented in 1959 the affine de Casteljau algorithm (Fig. 4) with reference to barycentric coordinates $p, q$, in a triangular scheme[16]:

[16] While de Casteljau originally denoted the control points as $A, B, C, D$, which in the first step lead to $G, H, K$, then to $I, J$ and finally to $M$ on the cubic curve, we will use the notation $\mathbf{c}_0, \mathbf{c}_1, \mathbf{c}_2, \mathbf{c}_3$, which in step 1 leads to the points $\mathbf{c}^1_0, \mathbf{c}^1_1, \mathbf{c}^1_2$ etc.

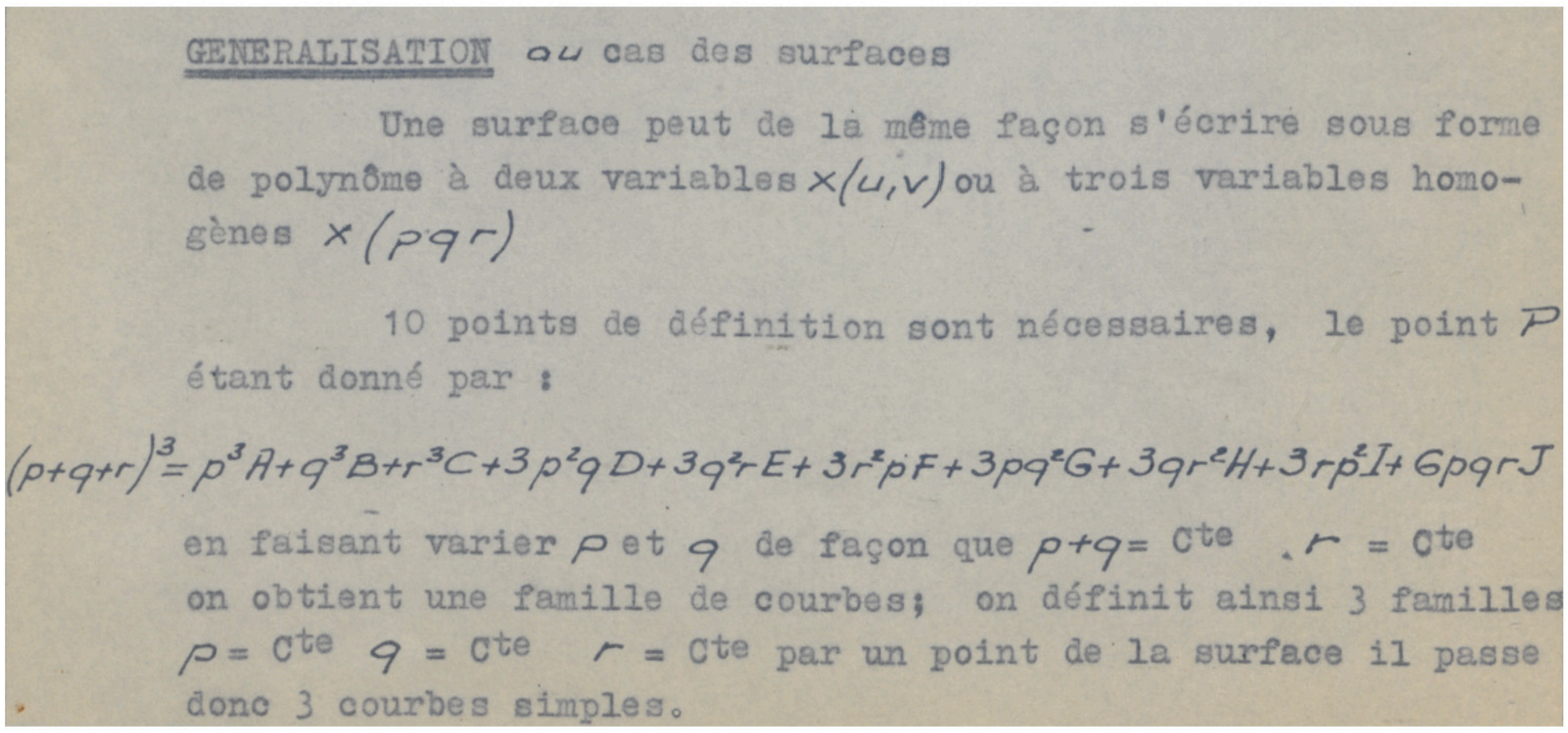
GENERALISATION ou cas des surfaces

Une surface peut de la même façon s'écrire sous forme de polynôme à deux variables $x(u,v)$ ou à trois variables homogènes $x(pqr)$

10 points de définition sont nécessaires, le point $P$ étant donné par :

$$(p+q+r)^3 = p^3A+q^3B+r^3C+3p^2qD+3q^2rE+3r^2pF+3pq^2G+3qr^2H+3rp^2I+6pqrJ$$

en faisant varier $p$ et $q$ de façon que $p+q=$ Cte , $r$ = Cte on obtient une famille de courbes; on définit ainsi 3 familles $p$ = Cte $q$ = Cte $r$ = Cte par un point de la surface il passe donc 3 courbes simples.

**Fig. 5.** de Casteljau's initial approach to (triangular) surfaces: ten points $A, B, \dots, I, J$ together with barycentric coordinates $p, q, r$ define a point $P$ on the surface; fixing one coordinate lets the two other trace a line on the surface, hence three such lines are passing by $P$. (de Casteljau, 1959) (replication with kind permission of Stellantis).

$$\begin{array}{llll}
\mathbf{c}_0 & & & \\
\mathbf{c}_1 & \mathbf{c}_0^1 = [p\mathbf{c}_0 + q\mathbf{c}_1] \cdot \frac{1}{p+q} & & \\
\vdots & \vdots & \ddots & \\
\mathbf{c}_n & \mathbf{c}_{n-1}^1 = [p\mathbf{c}_{n-1} + q\mathbf{c}_n] \cdot \frac{1}{p+q} & \cdots & \mathbf{c}_0^n = [p\mathbf{c}_0^{n-1} + q\mathbf{c}_1^{n-1}] \cdot \frac{1}{p+q}
\end{array}$$

For $n = 3$, the curve $\mathbf{x}(p,q) = \mathbf{c}_0^3(p,q)$ could thus be calculated as

$$\mathbf{x}(p,q) = \frac{p^3\mathbf{c}_0 + 3p^2q\mathbf{c}_1 + 3pq^2\mathbf{c}_2 + q^3\mathbf{c}_3}{(p+q)^3}.$$

By substituting $t = q/(p+q)$, he got the now well-known form of a polynomial curve

$$\mathbf{x}(t) = \sum_{i=0}^{n} B_i^n(t)\,\mathbf{c}_i$$

with polynomials $B_i^n(t) = \binom{n}{i}(1-t)^{n-i}t^i$, which only later were identified as Bernstein polynomials (cf. Farouki (2012)).

### 3.1. The case of surfaces

The surface case was already included in de Casteljau's initial document from 1959; by using barycentric coordinates, he first formulated triangular surfaces (cf. Fig. 5). Tensor product surfaces appear then in the Citroën manual by de Casteljau (1963). Triangular and quadrangular surfaces can be denoted in affine as well as in projective form.

The fundamental importance of this algorithm for CAGD and its applications can barely be overestimated, as the reader can appreciate from the contributions of de Casteljau (1993), Farin (1993), Boehm and Müller (1999), Prautzsch (2023).

### 3.2. The rational case

In *Le Lissage*, de Faget de Casteljau (1990, p.24) explains:

*The fundamental reason for choosing the pole forms is that any intermediate construction point, given by barycentric definition, starting with the previous poles, is on the sheet of paper, if the previous ones are there.*

Although the theory easily extends to the projective space and rational curves and surfaces, he saw the risk of possible singularities between start and end point, and the ensuing problems, e.g., when steering a milling cutter. Therefore de Casteljau preferred to always come back to the simple control points and to apply his original algorithm.

Instead of rational curves, de Casteljau often worked with finite continued fractions, which represent rational numbers. Avoiding infinite continued fractions also allows avoidance of irrational or transcendent numbers.

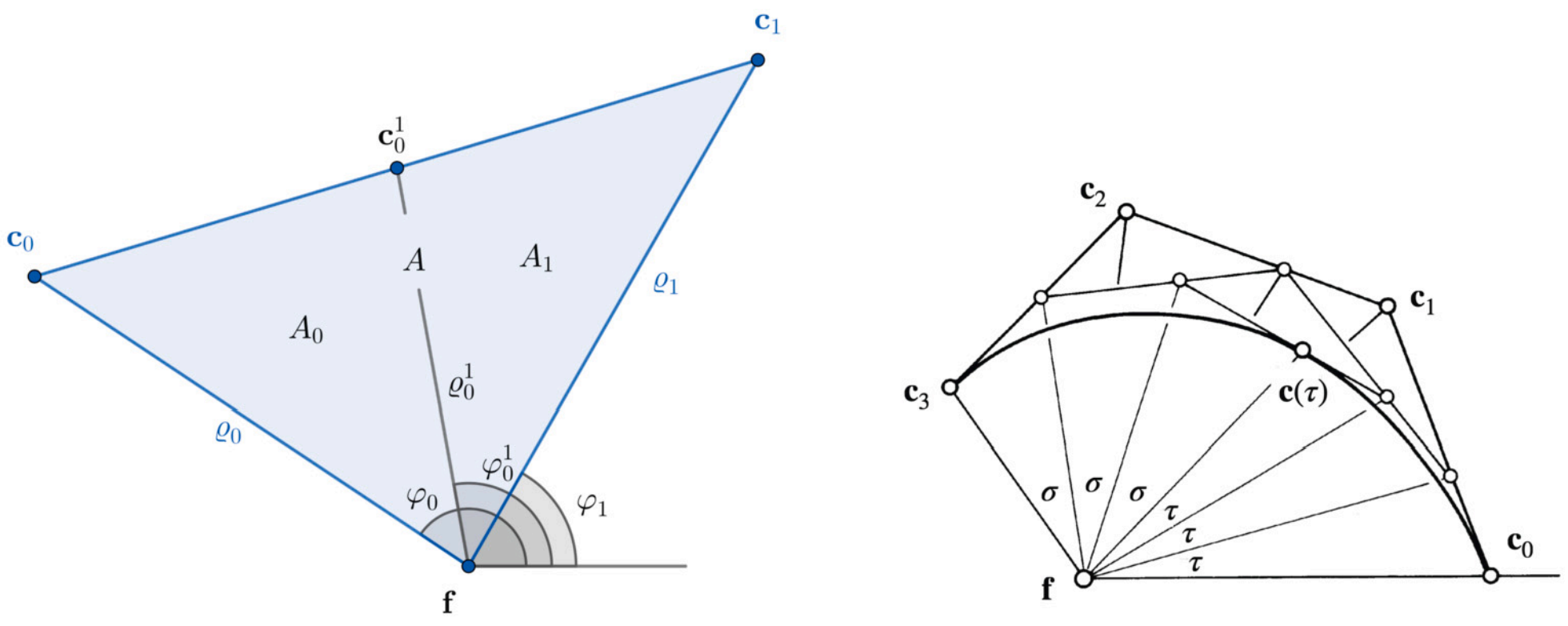

Fig. 6. Left: The area $A$ is divided into two areas $A_0$ and $A_1$. Right: de Casteljau's algorithm with angular parameter $\tau$, $n = 3$.

### 3.3. The angular case

A very different utilisation of his algorithm emerged with his presentation of 'focal splines' (de Faget de Casteljau, 1994). Seen from a 'focus' $\mathbf{f}$, the area $A$ of the triangle $\mathbf{c}_0\mathbf{f}\mathbf{c}_1$ can be split into two smaller areas $A_0, A_1$ by a separating line $\mathbf{f}\mathbf{c}_0^1$ (see Fig. 6). By using some auxiliary angle $\varphi_1$, the areas are calculated as

$$2A = \varrho_0\varrho_1 \sin(\varphi_0 - \varphi_1)$$
$$2A_0 = \varrho_0\varrho_0^1 \sin(\varphi_0 - \varphi_0^1)$$
$$2A_1 = \varrho_0^1\varrho_1 \sin(\varphi_0^1 - \varphi_1)$$

and with $A = A_0 + A_1$, it follows

$$\frac{\sin(\varphi_0 - \varphi_1)}{\varrho_0^1} = \frac{\sin(\varphi_0 - \varphi_0^1)}{\varrho_1} + \frac{\sin(\varphi_0^1 - \varphi_1)}{\varrho_0}$$

or in its additive recursive notation

$$\frac{1}{\varrho_0^1} = \frac{1}{\varrho_0}\frac{\sin\sigma}{\sin(\sigma+\tau)} + \frac{1}{\varrho_1}\frac{\sin\tau}{\sin(\sigma+\tau)} \qquad (2)$$

A recursive angular subdivision delivers a point on a focal curve; for $n = 2$ this is a conic with one focus in $\mathbf{f}$, see Fig. 6, which easily extends to spline curves. Focal splines have further been studied with an emphasis on their polar coordinates as $p$-Bézier curves or $p$-splines by Sanchez-Reyes (1990, 2002) or Casciola and Morigi (1996).

**Addendum**

The choice of a specific parameter in de Casteljau's algorithm leads him also to the algorithm of Vernier, which generates a *Vernier scale*, as described, for example, by Randow (2010). De Casteljau provides a simple access to Randow's generalisation in *Le Lissage*: Take a parameter range from $p/q$ to $r/s$; the evaluation at $(p+q)/(r+s)$ delivers the poles of the Vernier algorithm. (de Faget de Casteljau, 1990, p. 28)

## 4. Blossoms are polar forms are osculants

*Pole theory is entirely algebraic. The calculation of B-spline coefficients is based on analysis, and the trouble with analysis is that it often leads to the conditions imposed at the outset being forgotten at the end of the calculation. The use of algebra does, however, allow a complete inventory of all possible cases to be made.* — de Faget de Casteljau (1986, p.34)

De Casteljau's geometric approach of designing a curve (or surface) by use of a control polygon (or net) was kept an industry secret within Citroën; its intellectual property was solely documented internally by de Casteljau (1959, 1963). Yet, Pierre Bézier from Renault, who was aware about the method at large, arrived at a similar geometric approach on curves, and disseminated the idea in the late 1960s, particularly at conferences in the United States and the U.K. (Bézier, 1966, 1968, 1970, 1971). This inspired Robin Forrest (1972)[17] to see the relationship with the Bernstein basis, as well as Rich Riesenfeld (1973) and Ken Versprille (1975, 2013) to generalise it to integral B-Spline or Non-Uniform Rational B-Spline curves and surfaces. After Boeing integrated the then abbreviated NURBS, the question arose in France at Citroën in the late 1970s, how to integrate the B-Spline's advantage of smoothness into their SPAC system.

[17] Forrest remembers (personal communication to Alvy Ray Smith, 1 January 2024). *"At General Motors I was given a demonstration of 'the crazy way Renault (pronounced Wren Alt!) design curves'. [...] I re-wrote Bézier's mathematics in terms of vertices and discovered Bernstein polynomials were the basis. I knew about Bernstein*

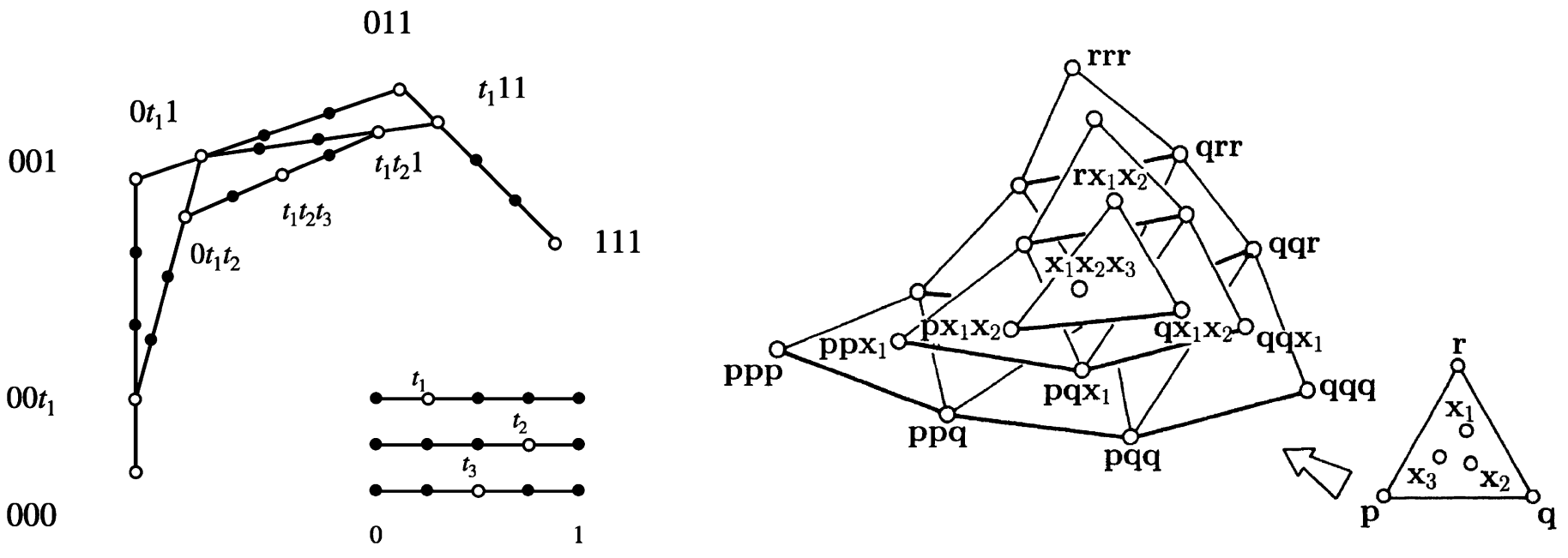

**Fig. 7.** Polar form for a cubic curve (left) and for a triangular surface (right), $n = 3$; image from (Prautzsch et al., 2002).

De Casteljau, who at that time did not know about Riesenfeld's findings, revisited his original approach from the 1960s and found a novel way to marry both, continuity and local control, first documented by de Faget de Casteljau (1983) in a manuscript, which was presumably intended to be submitted to the Academie of Sciences. His intention was to not lose the information on interpolation for the price of smooth $C^k$ continuity, and vice versa. The result were his *pôles généralisées*, inextricably linked with a quasi-interpolation scheme.

### *4.1. Pôles généralisées*

Let us introduce the polar form in the algebraic way that de Casteljau proposed in *Formes à Pôles*. The elementary symmetric functions are defined as

$$\sigma_1^n = x_1 + x_2 + \cdots + x_n$$

$$\sigma_2^n = x_1x_2 + x_1x_3 + \cdots + x_2x_3 + \cdots + x_{n-1}x_n$$

$$\vdots$$

$$\sigma_n^n = x_1x_2 \cdots x_n,$$

with the abbreviation $\sigma_i^n := \sigma_i(x_1, \ldots, x_n)$. Their recursive property

$$\sigma_i^n = \sigma_i^{n-1} + x_n \cdot \sigma_{i-1}^{n-1}$$

can be used to generate symmetric functions of higher degrees.

The fundamental theorem on elementary symmetric functions by Waring (1762) reads *"Each rational symmetric function in $x_1, \ldots, x_n$ can be written as a polynomial in $\sigma_1^n, \ldots, \sigma_n^n$"*, which results in the fundamental equation for a polar form

$$\mathbf{p}_{1\ldots n} := \mathbf{p}(x_1, \ldots, x_n) = \mathbf{p}_0 + \sum_{i=1}^{n} \mathbf{p}_i \sigma_i^n. \tag{3}$$

We will now assume an increasing order for the $n$ parameters, $t_1 \le t_2 \le \cdots \le t_n$, and we will sometimes identify the index sequence with its polar form $\mathbf{p}_{123} = \mathbf{p}(x_1, x_2, x_3) = [123]$.

Three algebraic properties of elementary symmetric functions translate into geometric properties of polar forms (Fig. 7):

- **diagonality:**

  The elementary symmetric functions reproduce binomial coefficients if $x_j = x$.

  $$\sigma_i^n = \sigma_i(x, x, \ldots, x) = \binom{n}{i} x^i$$

  The polar form reproduces the curve if $t_k = t$.

  $$\mathbf{p}_{t,t,\cdots,t} = [\mathtt{tt}\ldots\mathtt{t}] = \mathbf{p}(t) = \mathbf{p}_0 + \sum_{i=1}^{n} \mathbf{p}_i \binom{n}{i} t^n$$

*polynomials from approximation theory where they were largely used in a proof of the Weierstrass Approximation Theorem because they converged slowly but replicated all derivatives, the clue to their variation diminishing property. This was [already] documented in CAD Group Document 24 in mid-1969 [...]"*

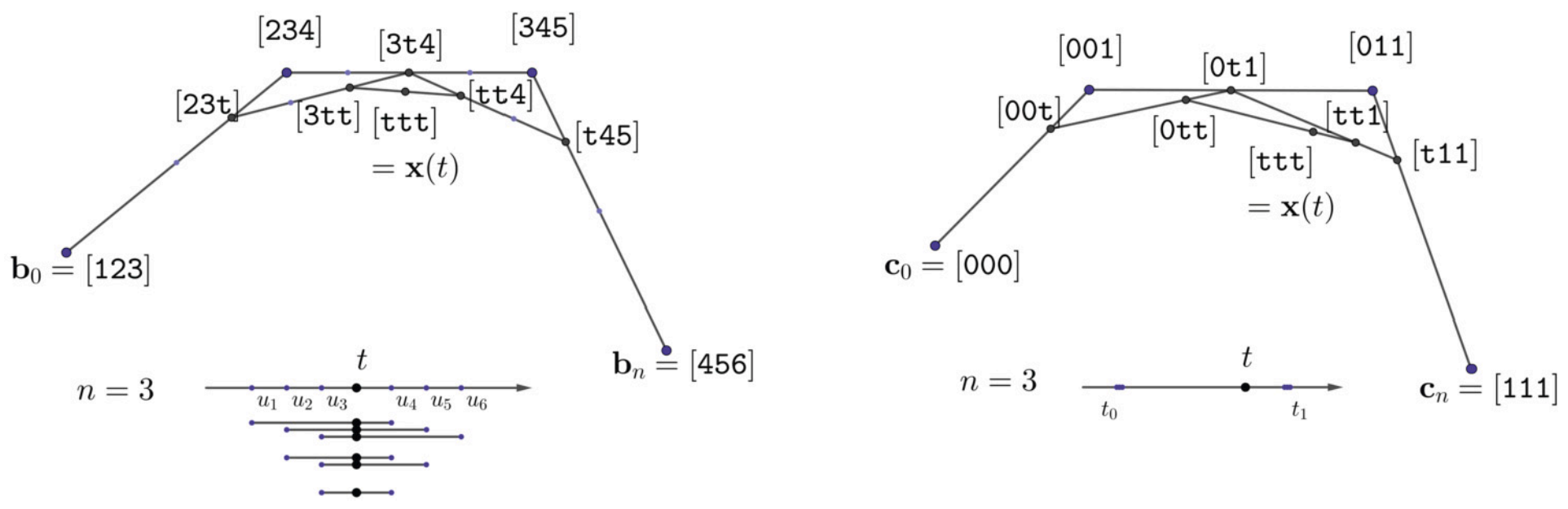

**Fig. 8.** Two famous curve evaluation algorithms and their poles: index reduction with progressive poles $\mathbf{b}_i$ (de Boor, left) and simple poles $\mathbf{c}_i$ (de Casteljau, right).

- **multi-affine:**

A linear combination of symmetric functions equals the symmetric function of its linear combination.

Fixing one parameter $t$ in the polar form of degree $n$ results in an affine combination of curves of degree $n-1$, the *osculants* at $\mathbf{p}(t)$.

$$\begin{aligned}\mathbf{p}(t_1,\cdots,t_{n-1})(t) &= \mathbf{p}(t_1,\cdots,t_{n-1})(t_0)\cdot(1-\alpha)+\mathbf{p}(t_1,\cdots,t_{n-1})(t_n)\cdot\alpha \\ &= \mathbf{p}(t_0,t_1,\cdots,t_{n-1})\cdot(1-\alpha)+\mathbf{p}(t_1,\cdots,t_{n-1},t_n)\cdot\alpha \qquad \text{with} \quad \alpha=\frac{t-t_n}{t_0-t_n}\end{aligned}$$

$\mathbf{p}(t_1,\cdots,t_{n-1})$ depends on $(n-1)$ variables, different to $\mathbf{p}(t_0,\cdots,t_{n-1})$ and $\mathbf{p}(t_1,\cdots,t_n)$ with its $n$ variables, which lets de Casteljau call the operation *index reduction*.

- **symmetry:**

A symmetric function is invariant under permutations of their arguments.

The evaluation of a polar form is invariant of the order of arguments (A-frame property).

### 4.2. Some remarkable relations

Let us now look at the influence of de Casteljau's *index reduction* in three remarkable settings. Note that a change of parameters is contragredient to a change of points.

*The algorithm of de Boor*

Let us choose a series of $2n$ indices in an increasing order and a parameter $t\in[u_n,u_{n+1}]$. We evaluate the polar form of consecutive indices

$$\mathbf{b}_i=\mathbf{p}(u_{i+1},\ldots,u_{i+n}),\quad i=0,\ldots,n.$$

The index reduction here corresponds to the algorithm of de Boor (1972) and Cox (1972)[18] (cf. Fig. 8), the points $\mathbf{b}_i$ are control points of the spline curve

$$\mathbf{x}(t)=\sum_{i=0}^{n}\mathbf{b}_i N_i^n(t)$$

The B-Splines $N_i^n(t)$ form a polynomial basis, thus their polar form can be formed from elementary symmetric functions $\sigma_i^n$.

*The algorithm of de Casteljau*

Let us now coalesce the first and last $n$ indices into two: $u_1=\cdots=u_n=:t_0,\quad t_1:=u_{n+1}=\cdots=u_{2n}$ and choose a parameter $t\in[t_0,t_1]$. We again evaluate the polar form of consecutive indices, keeping in mind their multiplicity. Index reduction here corresponds to the algorithm of de Casteljau (1959)[19] (cf. Fig. 8), the points $\mathbf{c}_i$ are the control points of the curve

$$\mathbf{x}(t)=\sum_{i=0}^{n}\mathbf{c}_i B_i^n(t), \tag{4}$$

they are results of the polar form $\mathbf{p}(\underbrace{t_0,\ldots,t_0}_{n-i},\underbrace{t_1,\ldots,t_1}_{i})$.

[18] The recurrence relation is due to Cox, de Boor and Mansfield, which due to de Boor was already known in the 1930's (de Boor and Pinkus, 2003).

[19] The recurrence relation was found by de Casteljau, which due to Boehm (1993) was already discovered by Haase (1870).

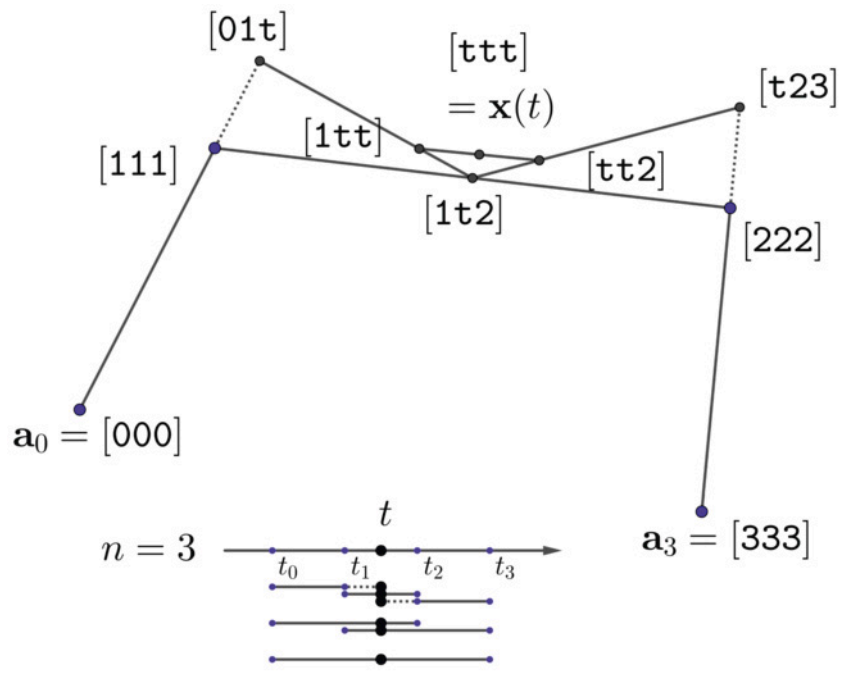

**Fig. 9.** Aitken's curve interpolation algorithm and its poles $\mathbf{a}_i$.

The polynomials $B_i^n(t)$ form a polynomial basis, thus their polar form can also be formed from elementary symmetric functions $\sigma_i^n$.

*The algorithm of Aitken-Neville*

Let us now associate $(n+1)$ consecutive indices with points $\mathbf{a}_0, \ldots, \mathbf{a}_n$, such that

$$\mathbf{x}(t) = \sum_{i=0}^{n} \mathbf{a}_i L_i^n(t) \tag{5}$$

with $t \in [t_0, t_n]$. Index reduction corresponds to the algorithm of Aitken (1932)[20] (cf. Fig. 9) to generate an interpolation polynomial of degree $n$. The Lagrange polynomials $L_i^n(t)$ form a polynomial basis and can thus be formed from the $\sigma_i^n$. We also find the recursion

$$L^i_{01\ldots(n-1)n}(t) = \alpha_0^n L^i_{01\ldots(n-1)}(t) + (1-\alpha_0^n) L^i_{1\ldots(n-1)n}(t)$$

$$\text{where} \quad \alpha_k^l = \frac{t - t_l}{t_k - t_l}$$

$$\text{and} \quad L_i^n(t) = L^i_{0\ldots n}(t) = \frac{(t-t_0)\cdots(t-t_{i-1})(t-t_{i+1})\cdots(t-t_n)}{(t_i-t_0)\cdots(t_i-t_{i-1})(t_i-t_{i+1})\cdots(t_i-t_n)}$$

which turns to a construction of intermediate points

$$\mathbf{a}_i^r = \mathbf{a}_{i+1}^{r-1} \cdot \alpha_{i+r}^i + \mathbf{a}_i^{r-1} \cdot (1 - \alpha_{i+r}^i)$$

### *4.3. Primitive, simple and progressive poles*

Consider given parameters $t_i$ and a sequence of their indices $i = \ldots, 1, 1, 2, 3, 3, 4, 4, 5, 6, 6, 6, 7, \ldots$ A given degree $n$ leads to a sequence of poles with consecutive indices, which de Casteljau calls **primitive poles**, in our case:

$$(n=3) \qquad \mathbf{p}_{112}, \mathbf{p}_{123}, \mathbf{p}_{233}, \mathbf{p}_{334}, \mathbf{p}_{344}, \mathbf{p}_{445}, \mathbf{p}_{456}, \mathbf{p}_{566}, \mathbf{p}_{666}, \mathbf{p}_{667}, \cdots$$

The poles with $n$ identical indices are points on the curve $\mathbf{x}(t)$, in our example $\mathbf{p}_{666}$. More than $n$ consecutive indices lead to identical poles.

Maximal insertion of further indices to the sequence leads to an $n$-fold repetition of a parameter, which corresponds to the **simple poles** of the curve in the respective segment. With $t$ between the parameters $t_3$ and $t_4$, we have simple poles $\mathbf{p}_{333}, \mathbf{p}_{334}, \mathbf{p}_{344}, \mathbf{p}_{444}$, and de Casteljau's algorithm can be visualised in triangular form:

$\mathbf{c}_0 = [333]$
$\mathbf{c}_1 = [334] \quad [33\mathtt{t}]$
$\mathbf{c}_2 = [344] \quad [3\mathtt{t}4] \quad [3\mathtt{tt}]$
$\mathbf{c}_3 = [444] \quad [\mathtt{t}44] \quad [\mathtt{tt}4] \quad [\mathtt{ttt}] = \mathbf{x}(t)$

Poles that correspond to a sequence of indices without repetition are called **progressive poles**. Given $i = \ldots, 1, 2, 3, 4, 5, 6, 7, \ldots$ we have progressive poles $\mathbf{p}_{123}, \mathbf{p}_{234}, \mathbf{p}_{345}, \mathbf{p}_{456}$ and can thus construct a point on the curve by de Boor's algorithm, again visualised in triangular form:

[20] The recurrence relation was found by Aitken and has been improved to its consecutive knot version by Neville (1934), as Miller (1968) reports.

$\mathbf{b}_0 = [123]$
$\mathbf{b}_1 = [234] \quad [23t]$
$\mathbf{b}_2 = [345] \quad [3t4] \quad [3tt]$
$\mathbf{b}_3 = [456] \quad [t45] \quad [tt4] \quad [ttt] = \mathbf{x}(t)$

**Addenda**

- The polynomials $B_i^n(t) = \binom{n}{i}(1-t)^{n-i}t^i$ are widely known as Bernstein polynomials, which often lets us forget about Bernstein's origin in probability theory seemingly unrelated to geometry. Instead, de Faget de Casteljau (1986, p.39) suggests a relation to Euler's beta and gamma functions, *"Representation with poles is to standard representation [in monomial form], what the Euler Beta function is to the Gamma function"*:

$$B(p,q) = \int_0^1 x^{p-1}(1-x)^{q-1}dx, \qquad \Gamma(p) = \int_0^\infty x^{p-1}e^{-x}dx = (p-1)!$$

$$B(p,q) = \frac{\Gamma(p)\Gamma(q)}{\Gamma(p+q)} = \frac{(p-1)!(q-1)!}{(p+q-1)!} = \frac{1}{(p+q-1)\binom{p+q-2}{p-1}}, \qquad B(n-k+1,k+1) = \frac{1}{(n+1)\binom{n}{k}}$$

- Despite the algebraic proximity of triangular patches (due to barycentric coordinates), de Casteljau considered the bi-parametric quadrangular surfaces to be the natural extension to curves, with the understanding that parameters should be adjusted when algebraic smoothing (next section) is applied. At that stage, knot removal or data reduction also comes into play (Goldman and Lyche, 1993). He had in mind to extend the polar forms to triangular surfaces, but decided that this would be too complicated for draftsmen.[21] The approach on triangular blossoms was instead brought forward by Ramshaw (1987) and Dahmen et al. (1992).

## 5. Algebraic smoothness, an application of polar forms

*It's more interpolation than smoothing, to be honest. The difference is that the degree of interpolation is deliberately limited. This is the role of the notion of degree of restitution, which is less novel than you might think. Present in older works on numerical interpolation, it is neglected by mathematicians who much prefer to substitute the term* et cætera *to indicate that one has just put one's finger in an infernal spiral that will never stop, as is the case for series, sequences or unlimited recurrences.* — de Faget de Casteljau (1990, p.138)

The following idea of de Casteljau is based on the way car makers worked in the 1950s and 1960s: designers would develop a clay model of the exterior, which was sampled, 'mathematised' and only then could be used for the subsequent steps such as tooling. The sampling result was called *'courbe à points'*, which turned into a *'courbe à pôles'* after mathematicians such as de Casteljau had worked on them.[22] We will now benefit from the fact that the same curve $\mathbf{x}(t)$ can be computed in at least three different ways, partly interpolative partly approximative, cf. section 4. In fact, de Faget de Casteljau (1983) already considered this proximity in his first manuscript on polar forms.

De Casteljau introduces the *degré de restitution*, $r$: assuming that $q$ points were sampled of an algebraic curve that was defined by polynomials of degree $r$ or above. The subsequent operations are to preserve this algebraic nature, and we should always be able to rebuild (or restitute) the original curve in a rigorous way, or to get as close to it as possible in terms of differentiation and continuity.

For simplicity, we assume that the curve $\mathbf{x}(t)$ is sampled at integer steps, leading to $q$ points $\mathbf{a}_1 = \mathbf{a}(t_1), \ldots, \mathbf{a}_q = \mathbf{a}(t_q)$ with $\Delta t_i = t_{i+1} - t_i = 1$, and we will identify $t_i = i$, w.l.o.g.[23] We want to construct a spline curve with $s$ segments of degree $n$ and an order of continuity $c = n - p$, which interpolates to the degree of $r$.

The general idea is as follows: we first interpolate $r+1$ points $\mathbf{a}_i$ (step 1); we then take into account the desired continuity and raise the indices to multiplicity $p$, such that we arrive at spline control points $\mathbf{b}_j$ (step 2); finally we translate those points towards control points $\mathbf{c}_k$, which we could use to perform the (classic) de Casteljau algorithm (step 3).

### 5.1. *The example* $(5,3,3)$

We will now illustrate de Casteljau's procedure for the case $(n=5, c=3, r=3)$ and $q=6$ in more detail:

#### *Step 1: Lagrange interpolation*

Let us label the six sampled points as $\mathbf{a}_{-2}, \mathbf{a}_{-1}, \mathbf{a}_0, \mathbf{a}_1, \mathbf{a}_2, \mathbf{a}_3$, with corresponding indices $-2, -1, 0, 1, 2, 3$ – de Casteljau chooses the coordinates without loss of generality. As the degree of restitution is $r=3$, we will evaluate four consecutive points each for a Lagrange interpolation, $L_{-2-101}, L_{-1012}, L_{0123}$.

---

[21] J.-L. Loschutz, personal communication, 2023.

[22] In his autobiography, de Casteljau mentions that opposing Peugeot colleagues instead referred to the homophonous *'formes à poiles'* (hairy shapes), whereas supporting Citroën colleagues called them *'formes à Paul'* (Paul's shapes) (Müller, 2024).

[23] The approach also works for non uniform knot vectors, as shown in *Formes à Pôles* on pages 96/97.

For the sake of simplicity, we pick $L_{-1012}$ with its node polynomial $(x+1)x(x-1)(x-2)$. The barycentric form of the Lagrange polynomial is $L_{-1012} = \sum_{i=-1}^{i=2} l_i \mathbf{a}_i$, and its weights $l_i$ can be calculated by dividing the node polynomial by one of its constituting factors (Berrut and Trefethen, 2004). De Casteljau's tabular arrangement by polynomial coefficients will help to maintain an overview, especially when dealing with higher degrees (de Faget de Casteljau, 1990, p.68):

$$\begin{array}{rrrrrll}
L_{-1012} : (x+1)x(x-1)(x-2) & x^3 & x^2 & x^1 & x^0 & & \\
\div(x+1): & 1 & -3 & \mathbf{2} & 0 & [\div -\mathbf{6}] & \to l_{-1} = \frac{(x-0)(x-1)(x-2)}{(-1-0)(-1-1)(-1-2)} = \frac{1x^3-3x^2+\mathbf{2}x^1+0x^0}{\mathbf{-6}} \\
\div x: & 1 & -2 & -1 & 2 & [\div +2] & \\
\div(x-1): & 1 & -1 & -2 & 0 & [\div -2] & \\
\div(x-2): & 1 & 0 & -1 & 0 & [\div +6] & \to l_2 = \frac{1x^3+0x^2-1x^1+0x^0}{+6}
\end{array}$$

Thus we get the Lagrange interpolant

$$L_{-1012}(x) = \frac{x^3-3x^2+2x}{-6}\mathbf{a}_{-1} + \frac{x^3-2x^2-x+2}{2}\mathbf{a}_0 + \frac{x^3-x^2-2x}{-2}\mathbf{a}_1 + \frac{x^3-x}{6}\mathbf{a}_2$$

*Step 2: Polar forms*

Let us raise the multiplicity of the indices to $p = 2$ such that the final spline curve gets the desired continuity $c = n - p = 3$. Evaluating the Lagrange polar at $n = 5$ consecutive knots will result in six spline control points $\mathbf{b}_j$, with each of the three Lagrange polynomials defining two such poles:

| | | | | | | | |
|---|---|---|---|---|---|---|---|
| | $-2$ | $-1-1$ | 00 | | | | $\mathbf{b}[-2,-1,-1,0,0]$ |
| $L_{-2-101}$ | | $-1-1$ | 00 | 1 | | | $\mathbf{b}[-1,-1,0,0,1]$ |
| | | $-1$ | 00 | 11 | | | $\mathbf{b}[-1,0,0,1,1]$ |
| $L_{-1012}$ | | | 00 | 11 | 2 | | $\mathbf{b}[0,0,1,1,2]$ |
| | | | 0 | 11 | 22 | | $\mathbf{b}[0,1,1,2,2]$ |
| $L_{0123}$ | | | | 11 | 22 | 3 | $\mathbf{b}[1,1,2,2,3]$ |

We will now substitute the monomials by their elementary symmetric functions, e.g., $x$ by $\sigma_1^5 = x_1 + x_2 + x_3 + x_4 + x_5$ and divide by their quantity. Finally, we evaluate the polar form of the Lagrange interpolant, in our case $L_{-1012}$, at the above-assigned consecutive indices $(x_1, x_2, x_3, x_4, x_5) = (-1, 0, 0, 1, 1)$ and $(0, 0, 1, 1, 2)$:

| | $[\div 10]$ | $[\div 10]$ | $[\div \mathbf{5}]$ | $[\div 1]$ |
|---|---|---|---|---|
| | $\sigma_3^5$ | $\sigma_2^5$ | $\sigma_1^5$ | 1 |
| $L_{-1012}[-1,0,0,1,1]$ : | $-1$ | $-1$ | $\mathbf{1}$ | 1 |
| $L_{-1012}[0,0,1,1,2]$ : | 2 | 5 | 4 | 1 |

From here, we can evaluate the polar forms of the above Lagrange polynomials, making use of de Casteljau's tabular arrangement (the components of the second summand are highlighted in the above two tables to better follow the calculation):

$$\begin{aligned}
L_{-1012}[-1,0,0,1,1] &= \left(\frac{1}{1}\cdot\frac{0}{-6} + \frac{\mathbf{1}}{\mathbf{5}}\cdot\frac{\mathbf{2}}{\mathbf{-6}} + \frac{-1}{10}\cdot\frac{-3}{-6} + \frac{-1}{10}\cdot\frac{1}{-6}\right)\mathbf{a}_{-1} + \left(\frac{2}{2} + \frac{-1}{10} + \frac{2}{20} + \frac{-1}{20}\right)\mathbf{a}_0 \\
&\quad + \left(\frac{0}{-2} + \frac{-2}{-10} + \frac{1}{-20} + \frac{-1}{-20}\right)\mathbf{a}_1 + \left(\frac{0}{6} + \frac{-1}{30} + \frac{0}{60} + \frac{-1}{60}\right)\mathbf{a}_2 \\
&= \frac{1}{60}\begin{bmatrix}-6 & 57 & 12 & -3\end{bmatrix}\begin{bmatrix}\mathbf{a}_{-1} & \mathbf{a}_0 & \mathbf{a}_1 & \mathbf{a}_2\end{bmatrix}^t \\
L_{-1012}[0,0,1,1,2] &= \frac{1}{60}\begin{bmatrix}-3 & 12 & 57 & -6\end{bmatrix}\begin{bmatrix}\mathbf{a}_{-1} & \mathbf{a}_0 & \mathbf{a}_1 & \mathbf{a}_2\end{bmatrix}^t
\end{aligned}$$

We have the same results for $L_{-2-101}$ and $L_{0123}$ in relation to their shifted control points, which allows us to write the spline control points $\mathbf{b}_j$ in matrix form $B = HA$:

$$\begin{bmatrix}\mathbf{b}_{-2}\\ \mathbf{b}_{-1}\\ \mathbf{b}_0\\ \mathbf{b}_1\\ \mathbf{b}_2\\ \mathbf{b}_3\end{bmatrix} = \begin{bmatrix} L_{-2-101}[-2,-1,-1,0,0] \\ L_{-2-101}[-1,-1,0,0,1] \\ L_{-1012}[-1,0,0,1,1] \\ L_{-1012}[0,0,1,1,2] \\ L_{0123}[0,1,1,2,2] \\ L_{0123}[1,1,2,2,3] \end{bmatrix} = \frac{1}{60}\begin{bmatrix} -6 & 57 & 12 & -3 & & \\ -3 & 12 & 57 & -6 & & \\ & -6 & 57 & 12 & -3 & \\ & -3 & 12 & 57 & -6 & \\ & & -6 & 57 & 12 & -3 \\ & & -3 & 12 & 57 & -6 \end{bmatrix}\begin{bmatrix}\mathbf{a}_{-2}\\ \mathbf{a}_{-1}\\ \mathbf{a}_0\\ \mathbf{a}_1\\ \mathbf{a}_2\\ \mathbf{a}_3\end{bmatrix}$$

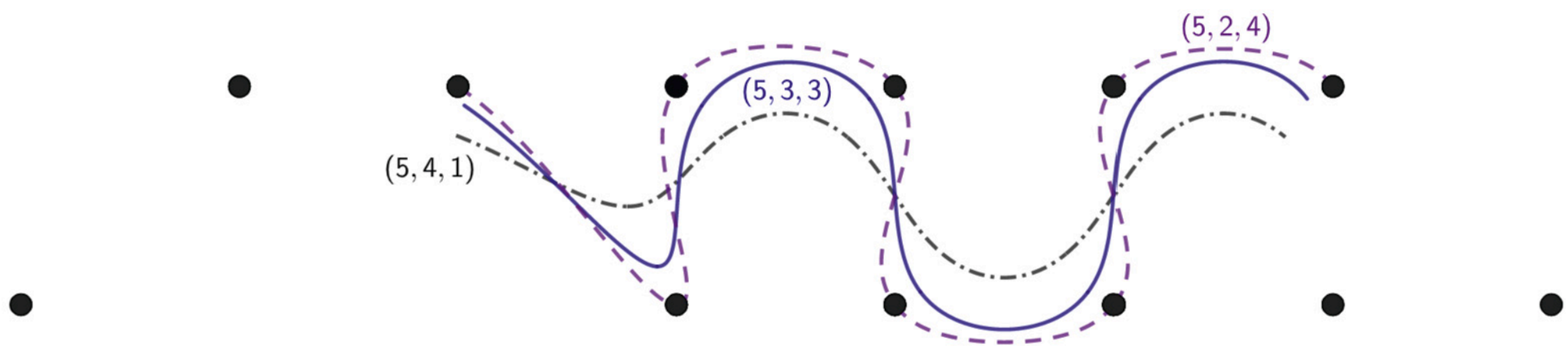

**Fig. 10.** Algebraic smoothing with characteristics $(n, c, r)$: all spline segments are curves of degree $n = 5$, (5,2,4) interpolates the given points with highest restitution $r = 4$, (5,4,1) has highest continuity $c = 4$, (5,3,3) mediates between continuity and interpolation.

### *Step 3: Index repetition*[24]

The final step is to construct simple control points $\mathbf{c}_k$. Applying the de Casteljau algorithm to these points will lead to points on the curve.

In our case, we will repeat indices until only two different values remain, which leads us to the respective simple poles (we simplify again and use the parameters instead of the polar form).

$$\begin{bmatrix} \mathbf{-2,-1,-1,0,0} \\ \mathbf{-1,-1,0,0,1} \\ -1,0,0,1,1 \\ 0,0,1,1,2 \\ 0,1,1,2,2 \\ 1,1,2,2,3 \end{bmatrix} \xrightarrow{K_1} \begin{bmatrix} \mathbf{-1,-1,0,0,0} \\ -1,0,0,0,1 \\ 0,0,0,1,1 \\ 0,0,1,1,1 \\ 0,1,1,1,2 \\ 1,1,1,2,2 \end{bmatrix} \xrightarrow{K_2} \begin{bmatrix} -1,0,0,0,0 \\ 0,0,0,0,1 \\ 0,0,0,1,1 \\ 0,0,1,1,1 \\ 0,1,1,1,1 \\ 1,1,1,1,2 \end{bmatrix} \xrightarrow{K_3} \begin{bmatrix} 0,0,0,0,0 \\ 0,0,0,0,1 \\ 0,0,0,1,1 \\ 0,0,1,1,1 \\ 0,1,1,1,1 \\ 1,1,1,1,1 \end{bmatrix}$$

The operation builds on de Casteljau's index reduction (cf. section 4), as we can see in the highlighted relation

$$[-1,-1,0,0,0] = \frac{(1-0)[-2,-1,-1,0,0] + (0-(-2))[-1,-1,0,0,1]}{1-(-2)} = \frac{1}{3}[-2,-1,-1,0,0] + \frac{2}{3}[-1,-1,0,0,1]$$

The matrices $K_1, K_2, K_3$ are multiplied and lead to the overall $(n = 5, c = 3, r = 3)$ matrix $C = K_3 K_2 K_1 B = K(HA) = (KH)A$, which relates the control points of the algorithms of de Casteljau $\mathbf{c}_k$, de Boor $\mathbf{b}_j$, and Aitken-Neville $\mathbf{a}_i$ (Fig. 11):

$$\begin{bmatrix} \mathbf{c}_0 \\ \mathbf{c}_1 \\ \mathbf{c}_2 \\ \mathbf{c}_3 \\ \mathbf{c}_4 \\ \mathbf{c}_5 \end{bmatrix} = \frac{1}{2} \begin{bmatrix} 1 & 1 & & & & \\ & 2 & & & & \\ & & 2 & & & \\ & & & 2 & & \\ & & & & 2 & \\ & & & & 1 & 1 \end{bmatrix} \cdot \frac{1}{2} \begin{bmatrix} 1 & 1 & & & & \\ & 1 & 1 & & & \\ & & 2 & & & \\ & & & 2 & & \\ & & & 1 & 1 & \\ & & & & 1 & 1 \end{bmatrix} \cdot \frac{1}{6} \begin{bmatrix} 2 & 4 & & & & \\ & 3 & 3 & & & \\ & & 4 & 2 & & \\ & & 2 & 4 & -6 & \\ & & & 3 & 3 & \\ & & & & 4 & 2 \end{bmatrix} \begin{bmatrix} \mathbf{b}_{-2} \\ \mathbf{b}_{-1} \\ \mathbf{b}_0 \\ \mathbf{b}_1 \\ \mathbf{b}_2 \\ \mathbf{b}_3 \end{bmatrix}$$

$$= \frac{1}{240} \begin{bmatrix} -7 & 28 & 198 & 28 & -7 & \\ -3 & -4 & 198 & 60 & -11 & \\ & -20 & 168 & 108 & -16 & \\ & -16 & 108 & 168 & -20 & \\ & -11 & 60 & 198 & -4 & -3 \\ & -7 & 28 & 198 & 28 & -7 \end{bmatrix} \begin{bmatrix} \mathbf{a}_{-2} \\ \mathbf{a}_{-1} \\ \mathbf{a}_0 \\ \mathbf{a}_1 \\ \mathbf{a}_2 \\ \mathbf{a}_3 \end{bmatrix}$$

## 5.2. *Optimal smoothing*

The parameters are related, they can be deduced from the number $q$ of given points and the number $s$ of desired segments, which leads to the configuration Table 2:

- $r = q - s \rightarrow$ degree of restitution (by Lagrange interpolants): $r$
- $q = r + \frac{m}{p} \rightarrow$ points to define a segment: $m$, multiplicity of the indices: $p$
- $c = n - p \rightarrow$ degree of each segment: $n = m - 1$, continuity of the spline curve: $c$

At the borders of the table, we can identify two well-known techniques: the high continuity $c = n - 1$ leads to B-spline interpolation, while the high restitution degree $r = n - 1$ leads to Riabenki splines, which in the case of a uniform knot vector are Catmull-Rom splines (Riabenki, 1974; Marchouk, 1980; Ramshaw, 1991).

[24] Also known as knot insertion.

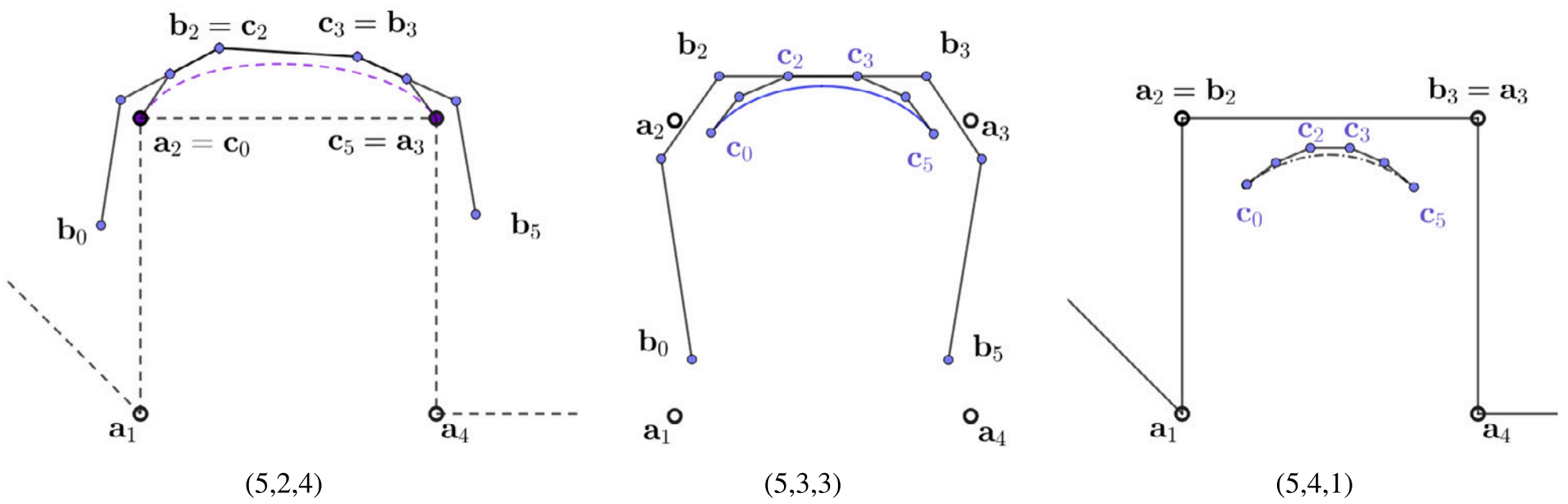

**Fig. 11.** Algebraic smoothing with characteristics $(n, c, r)$ – sampled points $\mathbf{a}_i$ on dashed lines, spline control points $\mathbf{b}_j$ and segment control points $\mathbf{c}_k$ on solid lines, curves as in Fig. 10. Left: spline control points $\mathbf{b}_j$ for (5,2,4) partly coincide with the segment control points $\mathbf{c}_k$ Centre: control points for (5,3,3) are combinations of both, $\mathbf{a}_i$ and $\mathbf{b}_j$ Right: spline control points $\mathbf{b}_j$ for (5,4,1) coincide with sampled points $\mathbf{a}_i$.

**Table 2**
Configuration table for smoothing, by de Faget de Casteljau (1990, p.73): Given $q$ points for $s$ segments, the optimal smoothing is performed with characteristics $(n, c, r)$: curve degree $n = m - 1$, continuity $c$, degree of restitution $r$ and knot multiplicity $p$. The table shows only $(n, c, r)$ as $p = n - c$ can be deduced; the emphasised configurations are mentioned in more detail by de Casteljau (1985); de Faget de Casteljau (1990). The case $(3, 1, 1)$ is a Catmull-Rom spline.

| $(n, c, r)$ | $s = 2$ | 3 | 4 | 5 | 6 |
|---|---|---|---|---|---|
| $q = 4$ | **(3,1,2)** | (3,1,1) CR spline | **(3,2,1)** B-spline | | |
| 6 | **(5,2,4)** | **(5,3,3)** | **(3,2,3)** B-spline | **(4,3,1)** B-spline | **(5,4,1)** B-spline |
| 8 | (7,3,6) | **(5,3,5)** | (7,5,4) | **(4,3,3)** B-spline | **(5,4,3)** B-spline |
| 10 | (9,4,8) | **(8,5,7)** | **(7,5,6)** | (9,7,5) | **(5,4,5)** B-spline |
| 12 | Riabenki: $n = 2p - 1$ | **(11,7,9)** | (11,8,8) | (9,7,7) | (11,9,6) |
| 14 | $c = p - 1$ | (11,7,11) | (11,8,10) | (9,7,9) | (11,9,8) |
| 16 | $r = 2p - 2$ | (14,9,13) | (15,11,12) | (14,11,11) | (11,9,10) |
| 18 | $p = q/2$ | (17,11,15) | (15,11,14) | (14,11,13) | (17,14,12) |
| 20 | | (17,11,17) | (19,14,16) | (19,15,15) | (17,14,14) |

The entries in this table are optimal in the sense of computational efficiency, as de Casteljau (1985, ch. 9) shows in *Formes à Pôles*. Taking the example of a target degree $n = 11$, he sums up the algebraic conditions for continuity and restitution and compares them with the number of coefficients in the transformation matrix. The spline case $(11, 10, 11)$ demands 264 coefficients for 397 conditions and interpolates $q = 22$ points, the Riabenki case $(11, 5, 10)$ demands 144 coefficients for 210 conditions and interpolates $q = 12$ points, while the optimum is the configuration $(11, 8, 10)$ with 168 coefficients for 267 conditions and $q = 14$ interpolation points. Thus Table 2 does not show all possible configurations, but only the most efficient ones.

### 5.3. *Applications*

Of the many applications of polar forms, let us quote a paragraph from *Formes à Pôles* (de Casteljau, 1985, p.110/111, translated from French):

*"The notion of poles can be used to define curves or surfaces, solutions to differential equations. This would bring us closer to the finite element methods, with the equations being transformed in such a way as to obtain relations on the poles. Certain trials allow to believe in the future of this research.*

*In practice, mastery of the $(n, c, r)$ characteristic allows the definition of all sorts of mathematical or physical tables, which was also the aim of splines. In addition to the continuity of first derivatives, the notion of degree of restitution guarantees some of the local series expansions. This avoids the inconveniences that arise from approximations using Tchebycheff polynomials. We can name ballistics and celestial mechanics as possible applications.*

*The use of several parameters allows for covering continuum physics, aerodynamics perhaps; the physicist will then be able to choose the right compromise between local study and a study covering a wider range, all thanks to the judicious choice of the $(n, c, r)$ characteristic."*

De Casteljau's general theory from section 4 is the subject of various publications under the name of *blossoming* by Ramshaw (1987, 1989), Seidel (1988, 1991), Dahmen et al. (1992), Stefanus and Goldman (1992), Gormaz (1993), Farin (1993), Mazure (1999, 2005). Yet, de Casteljau's algebraic smoothing approach has so far been quoted only by Lyle Ramshaw, in dense form in a final section in (1989) and in his SIGGRAPH course notes in (1991).

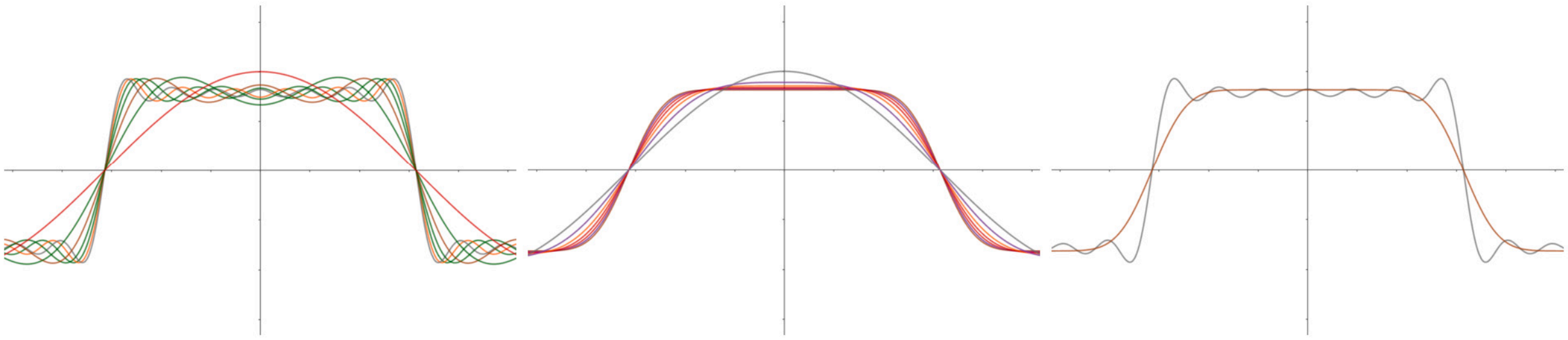

**Fig. 12.** Trigonometric Smoothing with the example of a square wave. Left: Fourier series, $p = 0..6$ Centre: Newton series, $p = 1..7$ Right: Direct comparison with 6 terms.

### Addendum

As already mentioned in the remarks of section 4, de Casteljau preferred the bi-parametric surfaces for polar forms, and consequently for this smoothing approach. Future research might address the connection of this algebraic smoothness to knot removal as well as to triangular patches.

## 6. Trigonometric smoothing

*A satisfactory smoothing requires the simultaneous regularisation of the function and its derivative.* — de Casteljau (Müller, 2024, p.42)

In *Le Lissage*, de Casteljau dedicates a chapter to boundary crossing and approaches the approximation of a square wave. Typically, it is approximated by orthogonal functions as

$$f(\varphi) = \cos\varphi - \frac{\cos 3\varphi}{3} + \frac{\cos 5\varphi}{5} + (-1)^p \frac{\cos(2p+1)\varphi}{2p+1} + \cdots \text{ etc.}$$

This Fourier series provokes the Gibbs phenomenon of oscillations that first overshoot the value of $\pm 1$ before swinging to a relatively stable value.

Instead, in order to support smoother derivatives, de Casteljau proposes an expansion of $\dfrac{\cos\varphi}{\sqrt{1-\sin^2\varphi}}$ around $\varphi = 0$, which leads to the following form

$$Kf(\varphi) = \cos\varphi - \frac{p-1}{p+1}\left[\frac{\cos 3\varphi}{3} - \frac{p-2}{p+2}\left[\frac{\cos 5\varphi}{5} - \frac{p-i}{p+i}\left[\frac{\cos(2i+1)\varphi}{2i+1} - \cdots - \frac{1}{2p-1}\frac{\cos(2p-1)\varphi}{2p-1}\right]\right]\right]$$

The result is an approximation without Gibbs phenomenon, as visualised in Fig. 12.

## 7. Tolerances

*Between 1958 and 1962, there were a large number of problems that needed a solution, even if it could be demonstrated that none existed: tolerance, speed and acceleration of stepper motors, pulse rate, cutting and scanning of curves or tiles, connections, approach to the workpiece, even worse avoid standstill at the end of the course (the cutter makes a hole), position of the convexity relative to the normal, asymptotics, etc.* — de Casteljau[25]

In 1966, Paul de Casteljau wrote the internal 17-pages report *Calcul de la Tolérance*, parts of which he published in a later article (de Casteljau, 1966; de Faget de Casteljau, 1997). The main question deals with the practical application of surfaces and the milling of their stamping tools.

Given a planar cubic arc $\varphi(x)$ between parameters $x_0$ and $x_1$, he calculates its tangents and the deviation from the line $l$ through $\varphi(x_0)$ and $\varphi(x_1)$. The tangents intersect their mutual parameters at $x_0$ in a distance of $d_1$ and at $x_1$ in a distance of $d_0$. Those distances are referred to as *tendencies*, as they predict the direction of the machine. If a tendency exceeds a certain threshold, the length of the arc should be shortened.

The distances to line $l$ can be seen as a function $f(x)$, its derivative leads to the maximum deviation between $x_0$ and $x_1$, which turns out to be, $(n = 3)$:

$$E_{max} = \frac{1}{9}\left[|d_0 + d_1| + \frac{d_0^2 - d_0 d_1 + d_1^2}{|d_0 - d_1| + 2\sqrt{\Delta}}\right] \qquad \text{with } \Delta = d_0^2 - d_0 d_1 + d_1^2$$

[25] de Casteljau, personal communication to W. Boehm (translated from French), 9 November 1995.

An extrapolation of the arc towards $x_0 + \rho(x_1 - x_0)$ leads to tendencies $D_0, D_1$, which help to predict future tendencies and to decide on the length of the arc:

$$D_0 = \rho^2[d_0 + (1-\rho)(d_0 - d_1)] \qquad D_1 = \rho^2[d_1 + 2(1-\rho)(d_0 - d_1)]$$

Moving to the case of surfaces, de Casteljau designates $M$ as the point to be machined, $P$ as the point on the cutter surface, and $B$ as the mathematically defined point on the targeted surface. With $N$ denoting the unit normal at $B$, $M$ or $P$, we can express their relations as follows:

$$P = M + R\overrightarrow{N} \qquad M = B + p\overrightarrow{N}$$

where $R$ is the radius of the cutter, and $p$ the distance between the defined surface and the machined one. A curve segment on the targeted surface shall start at $M_0$ and end at $M_1$.

De Casteljau then introduces three tolerances and approaches their deviations:

1. Normal tolerance $T_1$: A spherical cutter centred at $P$ cuts a parallel surface at $P + T_1$ without intersecting the targeted surface; $T_1$ is referred to as the normal tolerance.
   The tendencies result in
   $$d_0 = \overrightarrow{N_0}.\overrightarrow{P_0P_1} \qquad d_1 = -\overrightarrow{N_1}.\overrightarrow{P_0P_1} = \overrightarrow{N_1}.\overrightarrow{P_1P_0}$$
   and thus the corresponding normal deviation turns
   $$E_N = \frac{1}{9}\left[|d_0 + d_1| + \frac{d_0^2 - d_0d_1 + d_1^2}{|d_0 + d_1| + \sqrt{\Delta}}\right] \qquad \text{with } \Delta = d_0^2 - d_0d_1 + d_1^2$$
2. Geodesic tolerance $T_2$: A cutter with radius $R + T_2$ cuts the targeted surface. The geodesic tolerance $T_2$ represents the minimum distance that the cutter surface should maintain from the targeted surface.
   Let us define the geodesic normal $\vec{n}$ as the unit vector of
   $$\left(\frac{\partial M}{\partial u}\right)^2 \overrightarrow{\frac{\partial M}{\partial v}} - \left(\frac{\partial M}{\partial u}\frac{\partial M}{\partial v}\right)\overrightarrow{\frac{\partial M}{\partial u}}$$
   Then the tendencies are
   $$\delta_0 = \overrightarrow{n_0}.\overrightarrow{M_0M_1} \qquad \delta_1 = \overrightarrow{n_1}.\overrightarrow{M_1M_0}$$
   and the geodesic deviation is
   $$E_G \approx \frac{h_m^2}{2R} \quad \text{with} \quad h_m = \frac{1}{9}\left[|\delta_0 + \delta_1| + \frac{\delta_0^2 - \delta_0\delta_1 + \delta_1^2}{|\delta_0 + \delta_1| + \sqrt{\delta_0^2 - \delta_0\delta_1 + \delta_1^2}}\right]$$
3. Tolerance between grooves $T_3$: After the cutter has completed its work, the surface that is parallel to the machined one at a distance of $T_3$ completely disappeared. $T_3$ is referred to as the tolerance between grooves.
   For a given groove deviation $E_S$, de Casteljau finds the groove width as
   $$\Delta v = 2\sqrt{\frac{2eE_S}{eg - f^2}}$$
   where $e, f, g$ are calculated by using the first fundamental form from differential geometry and its derivative.

While de Casteljau further investigated the interdependencies of the deviations, the most significant among them is a constraint on the sum

$$E_N + E_G \leq T$$

which he ultimately employed to advance the cutter step-by-step. We will not delve into the detailed algorithm here.

## 8. Intersections

*The notion of digital machining is linked to that of tolerance, which conditions the validity of the interpolation (linear or not) between the calculation points. We are quickly led to consider the calculation of the difference between linear interpolation and parabolic interpolation. But experience quickly demonstrates the weakness of this solution (inflections, asymptotes).* — de Casteljau (1966)

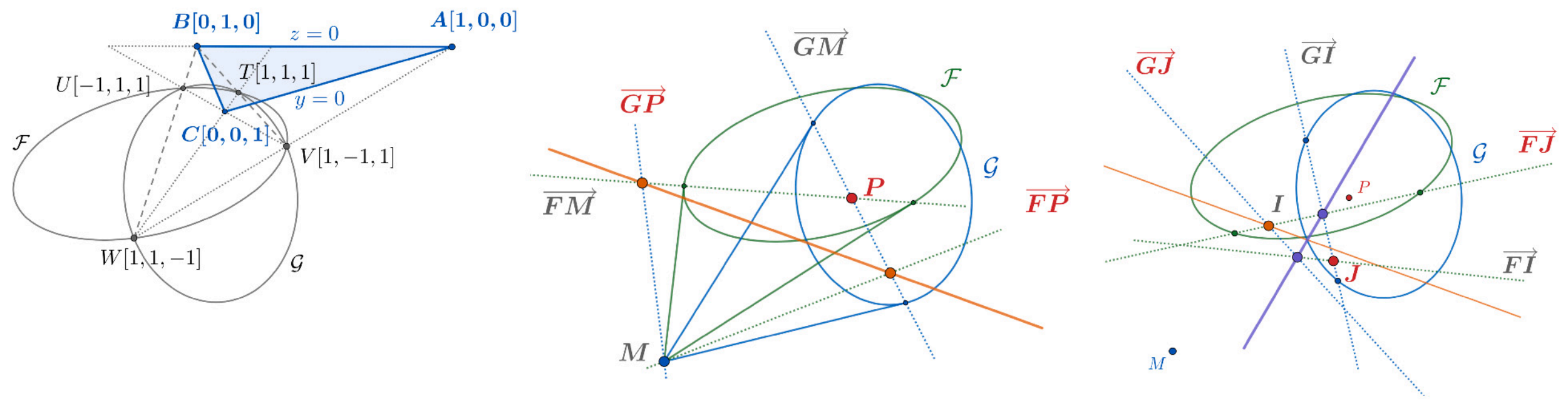

**Fig. 13.** Left: The intersections of two projective conics $\mathcal{F}, \mathcal{G}$ and their coordinates with respect to their diagonal triangle (de Casteljau, 1998) Centre: The connecting line of the polar intersections approximates a common root of $\mathcal{F}$ and $\mathcal{G}$. Right: The polars of $I$ and $J$ intersect at two points; their connection improves the approximation of the common root.

Another one of de Casteljau's challenges in the 1960s was to find intersections between mathematically well-described parts. He has published his ideas in *Le Lissage* as well as in two articles (de Casteljau, 1998; de Faget de Casteljau, 2000), and a first reference dates back to his early notes for the design school (de Casteljau, 1963, p. 38). He proposes one general approach and one approach that works for $n = 2$ only.

In both cases, we want to find a root of the polynomial

$$f(x) = -a_0 + a_1 x + a_2 x^2 + a_3 x^3 + \cdots$$

in the vicinity of $x = 0$, which can also be seen as a MacLaurin power series.

In a common first step, we truncate the series after the linear term, such that $a_0 = a_1 x$ delivers a first value $x_1 = a_0/a_1$. Repetitive reinsertion to the truncated $f$ leads to a better value $x$, which describes the classical Newton method of approximation.

Let us now truncate $f$ after the quadratic term. To make use of a similar method, we replace the square partly by the first solution $x_1$ and get

$$f(x) = -a_0 + a_1 x + a_2 x_1 x$$

Thus $a_0 = (a_1 + a_2 x_1)x$ leads to a solution $x_2 = a_0/(a_1 + a_2 a_0/a_1)$, which is known as the Whittaker method (1918).

In general, a truncation after $a_i$ and step-by-step replacement of the power by its symmetric function leads to a Horner-like procedure:

$$\begin{aligned}
i = 1: &\quad a_0 = a_1 x && \rightarrow \quad x_1 = \frac{a_0}{a_1} && \text{(Newton)}\\
i = 2: &\quad a_0 = (a_1 + a_2 x_1)x && \rightarrow \quad x_2 = \frac{a_0}{a_1 + a_2 \frac{a_0}{a_1}} && \text{(Whittaker)}\\
i = 3: &\quad a_0 = (a_1 + (a_2 + a_3 x_1)x_2)x && \rightarrow \quad x_3 = \frac{a_0}{a_1 + a_2 \frac{a_0}{a_2 + a_3 \frac{a_1}{a_2}}} && \text{(continued fractions)}\\
&\quad etc.,
\end{aligned}$$

which de Faget de Casteljau (1990, p. 32) labels *a polar form in mathematical analysis*. While he observed that a change in the order of the $x_i$ as well as an index repetition (a.k.a. knot insertion) lead to worse results, this still warrants further research.

He could generalise the approach to higher dimensions. For $n = 2$, we have two equations in two variables, $n = 3$ reads

$$\mathbf{f}(\mathbf{x}) = \mathbf{g}(\mathbf{x}) = \mathbf{h}(\mathbf{x}) = 0 \quad \text{with } \mathbf{x} = [x, y, z]$$

The process begins with $\mathbf{f}(\mathbf{x})$ to receive first $x_1$, then $y_1$ and $z_1$, before reusing these approximated values with $\mathbf{g}(\mathbf{x})$ and then using the new results with $\mathbf{h}(\mathbf{x})$.

The case of $n = 2$ allows a different, more geometrical path. The two polynomials represent projective conics, and their coordinates can be simplified using the diagonal triangle (cf. Fig. 13 left), which, however, remains unknown. An analytical approach led de Casteljau towards excessively high degrees. Consequently, he followed his geometric intuition.

Let $M^t F M = ax^2 + 2fxy + by^2 + 2ex + 2dy + c = 0$ be one of the quadratic forms,

$$F = \begin{bmatrix} a & f & e \\ f & b & d \\ e & d & c \end{bmatrix} \qquad M^t = \begin{bmatrix} x & y & z \end{bmatrix} \qquad M = \begin{bmatrix} x \\ y \\ z \end{bmatrix}$$

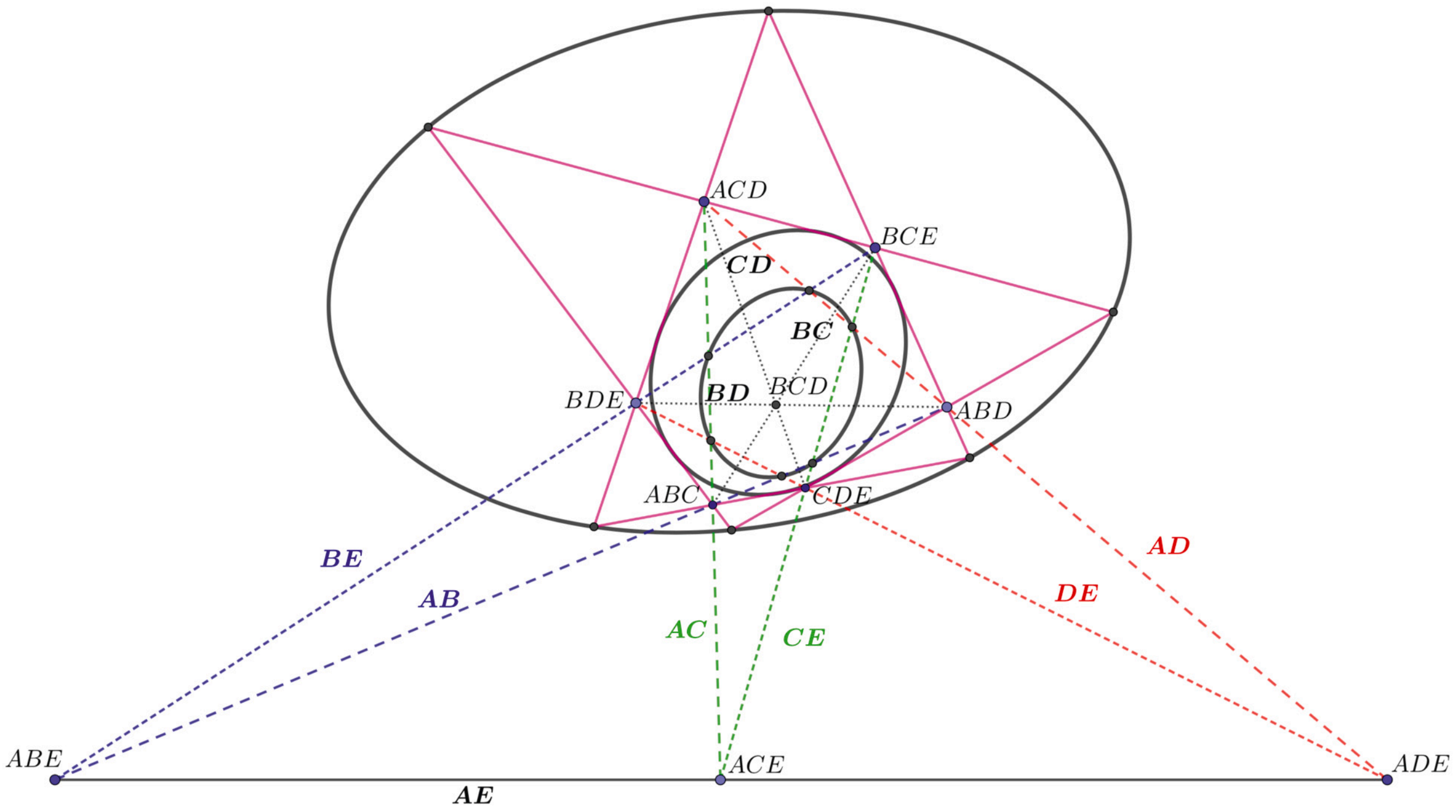

**Fig. 14.** A configuration of 3 conics and 16 lines in Projective Geometry, based on de Faget de Casteljau (1999b, p.100).

We begin with a point $M$, whose polars $\overrightarrow{FM}, \overrightarrow{GM}$ intersect at conjugate point $P$, expressed as a vector product

$$P = \overrightarrow{FM} \wedge \overrightarrow{GM}$$

The polars of $P$ intersect at $M$, which is conjugate by construction. The remaining two intersections of a polar of $P$ with a polar of $M$ form a line, which intersects both conics with quadratic convergence:

$$I = \overrightarrow{FM} \wedge \overrightarrow{GP} + \overrightarrow{FP} \wedge \overrightarrow{GM}$$

The method can be re-iterated starting with point $I$ on this line, whose polars intersect at conjugate point $J$, and the respective polar intersections deliver $\Omega$ with fourth-order convergence:

$$\Omega = \overrightarrow{FI} \wedge \overrightarrow{GJ} + \overrightarrow{FJ} \wedge \overrightarrow{GI}$$

The reader will find further aspects in the Dagstuhl conference article (de Casteljau, 1998). However, de Casteljau did not find a way to generalise the approach to three variables.

### Addendum

Ferdinand Schweins (1825) already observed, that the Taylor representation can also be developed via elementary symmetric functions. De Casteljau probably did not know about this link, although it is indicated by an article by Aitken (1926), which he might have known.

## 9. Projective or Euclidean geometry

*"It is very difficult to make things simple; it is always the complicated things that first come to mind" (Jean Philippe RAMEAU, the musician)* — de Casteljau[26]

De Casteljau's older brother Henri de Faget de Casteljau (1924-1988) was a career military man and historian, known for his genealogy research. Besides various articles on French heraldry, he was fascinated by geometry and had published a booklet on conics, elaborating on Pascal's theorem: *Traité des Coniques – "une pensée de Pascal"* (de Faget de Casteljau, 1973). Despite being

[26] de Casteljau, personal communication (translated from French), 29 September 1998. He used the quotation to introduce his figure of *16D+3C*, short for 16 *droites* (lines) and 3 *coniques* (conics), followed by nine (sic) further dictionary citations on the difficulty of expressing something complicated in simple terms.

quite different thinkers, Paul tried to draw the ideas of his brother Henri, which after months of studies, brought him to a figure that combines some of the elementary theorems of projective geometry, on which his brother Henri had worked, cf. Fig. 14.

The figure, composed of 3 conics and 16 lines (plus two hexagons), can explain the following theorems in an obvious way:

1. Desargues (1636): Five planes $A, B, C, D, E$ intersect in ten lines $AB, AC, \ldots, DE$. While $B, C, D$ intersect in a point, the intersection of the remaining planes $A, E$ forms the axis of perspectivity, each of the hexagon vertices working as centres of perspectivity.
2. Pascal (1639): Six points on the inner conic form three pairs of opposite sides which meet on the line $AE$. The figure includes ten such conics, which pass by six of the 15 intersection points. These 60 lines form the *Hexagrammum Mysticum*.
3. Brianchon (1810): The principal diagonals of the hexagon that is circumscribing the middle conic, intersect in one single point $BCD$. There are ten conics inscribed in such hexagons, each formed by six of the 15 possible lines.
4. Steiner (1827): Three Pascal lines of the inner conic intersect at one single point $BCD$. Assuming a clockwise numbering of the six inner conic points, these hexagons are in the order $163452, 143256, 123654$.
5. Poncelet (1822): The outer hexagon circumscribes the middle conic and is inscribed in the outer conic. Thus, there is an infinite number of such hexagons.
6. Kirkman (1850): Three Pascal lines of the outer conic intersect at one single point $BCD$. Assuming a clockwise numbering of the six outer conic points, these hexagons are in the order $135264, 351426, 513642$.

De Casteljau's article (1999b) also includes a construction scheme to generate the mentioned permutations: to define a hexagon inscribed in a conic we choose six vertices (out of ten) or six edges (out of 15). They generate 60 Pascal lines and 60 Kirkman nodes, 20 Steiner nodes and 20 Cayley lines, which leads to 15 Plücker lines and 15 Salmon points.

Paul de Casteljau often emphasises the limitations of geometries, in particular of projective geometry. He names the notion of degree, which is different to the one of class or genus and is challenged by curves of no degree like Koch snowflakes and other fractals. He also mentions the incompatibility with infinitesimal geometry, as the problem of differentiability at infinity shows, and he argues that geometry should not be reduced to a chapter of algebra.[27] He indicates that the hexagon study above could also be done by metric geometry in three different ways: using angles or trigonometric lines or a parameter (such as the tangent of the half angle).

### Addendum

De Casteljau also relates to Veronese (1877), whose studies on the *multimysticum* have been discussed again recently (Chipalkatti and Ryba, 2020).

## 10. Metric geometry

*We must return to the circle its exceptional role, because of its properties of angles and lengths. It is indeed "Her Majesty the Circle" that we should speak of, and their "Highnesses the Focal Conics" and not curves of second degree.* — de Casteljau[28]

In the mid-1990s, de Casteljau was prepared to write another book, one on metric geometry. The skeleton was to be organised as follows[29]:

1. **Preliminaries** Grids, inverse, polar form, triangular pavings, dualities, trisectors (Morley-Petersen) and extension to the 18 Morley triangles (cf. (Guy, 2007)), orthocentric tetragram
2. **Metric descriptions of alleged projective concepts** Systems of orthogonal circles, confocal conics, harmonic tetragram, eccentric anomaly; possibly also reflections on crystallography, inversion, composition of displacements, yield curves (Mohr envelope), caustics, quaternion
3. **Particular curves** Locus of foci, conical arc, strophoid, Cartesian oval, onduloid, cycloid, hypocycloid, epicycloid
4. **Geometric optics** Conjugate mirrors, Fermat's principle, backdrop canvas, aplanetic mirrors, *ricochet* (rebound), Cartesian oval, duality by refraction
5. **Further problems** Generalisation to $E^3$ or $E^4$ (review of 'confocal' conics), curves defined by three (or four) foci, quaternions as a local orthonormal frame in $E^3$ (or $E^4$), crystallography, caustics, ruled surfaces

We will look at only some phenomena here: grid, polar form, dualities, confocal conics, before elaborating an example from geometric optics, primarily based on an unpublished manuscript on metric geometry by de Faget de Casteljau (1995b).

---

[27] This reminds us of the French mathematician Sophie Germain (1776–1831), who said "Algebra is but written geometry and geometry is but figured algebra.".

[28] de Casteljau, personal communication (translated from French), 12 October 1995.

[29] de Casteljau, personal communication, 26 September 1995, 9 February 1996.

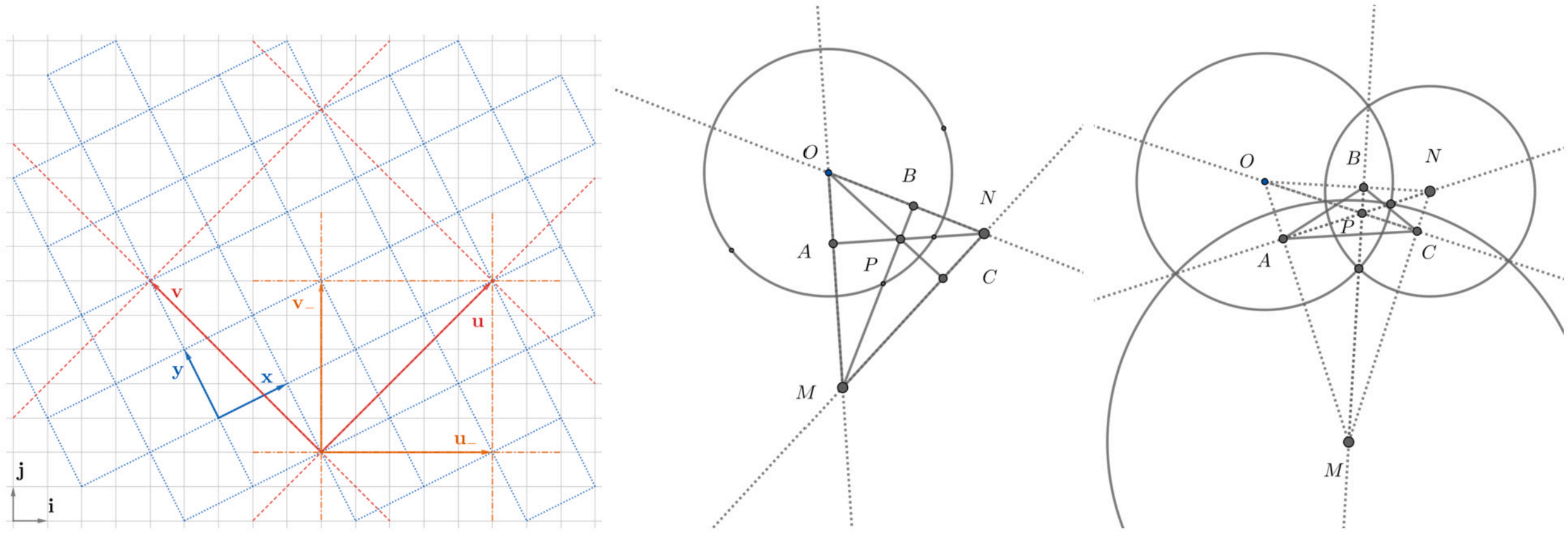

**Fig. 15.** Left: Three orthogonal coordinate systems. Centre: Metric polar form as an 'orthocentric tetragram'. Right: Four orthonormal circles around $M, N, O, P$, one has imaginary radius.

### 10.1. Grid and polar form

We take an orthogonal coordinate system with unit vectors $\mathbf{i}, \mathbf{j}$, on which we build another orthogonal grid:

$$\begin{aligned} \mathbf{x} &= x\,\mathbf{i} + y\,\mathbf{j} \\ \mathbf{y} &= x\,\mathbf{i} - y\,\mathbf{j} \end{aligned}$$

Let us now construct a third orthogonal coordinate system as a linear combination of $\mathbf{x}, \mathbf{y}$:

$$\begin{aligned} \mathbf{u} &= \xi\,\mathbf{x} + \eta\,\mathbf{y} &&= (\xi x - \eta y)\,\mathbf{i} + (\xi y + \eta x)\,\mathbf{j} \\ \mathbf{v} &= \xi\,\mathbf{y} - \eta\,\mathbf{x} &&= (\xi x - \eta y)\,\mathbf{j} - (\xi y + \eta x)\,\mathbf{i} \end{aligned}$$

The specific choice of $\xi = x, \eta = \pm y$ leads to

$$\begin{aligned} \mathbf{u}_+ &= (x^2 - y^2)\,\mathbf{i} + 2xy\,\mathbf{j} & \qquad \mathbf{v}_+ &= (x^2 - y^2)\,\mathbf{j} - 2xy\,\mathbf{i} \\ \mathbf{u}_- &= (x^2 + y^2)\,\mathbf{i} & \mathbf{v}_- &= (x^2 + y^2)\,\mathbf{j} \end{aligned}$$

For symmetry reasons, all four vectors $\mathbf{u}_+, \mathbf{u}_-, \mathbf{v}_+, \mathbf{v}_-$ have the same length $x^2 + y^2$. The homology $|\mathbf{u}_+| : |\mathbf{x}| = \xi|\mathbf{x}| : x|\mathbf{i}| = |\mathbf{x}|$ leads to the length of the 'intermediate' coordinate vector $|\mathbf{x}| = \sqrt{|\mathbf{u}_+|} = \sqrt{x^2 + y^2}$.

By abbreviating $x_0 = \xi x - \eta y$ and $y_0 = \xi y + \eta x$, we get

$$x_0 x + y_0 y = \xi(x^2 + y^2) = \xi|\mathbf{x}|^2,$$

which is the product of the length of $|\mathbf{x}|$ and the projection of $\mathbf{u}$ onto the $\mathbf{x}$-axis, so it is the scalar product $\langle \mathbf{u}, \mathbf{x} \rangle$, and also the polar form of $x^2 + y^2$.

Let us now look at three points $M, N, P$ and a circle $\mathcal{K}$ around $O$ with radius $r^2 = x^2 + y^2$, such that

$$x_M x_N + y_M y_N = r^2 = x_N x_P + y_N y_P = x_P x_M + y_P y_M$$

The projections of $N, M, P$ onto the lines $OM, ON, MN$ will be $A, B, C$ with

$$\begin{aligned} \langle OM, ON \rangle &= \langle OM, OA \rangle &&= \langle OB, ON \rangle \\ \langle ON, OP \rangle &= \langle ON, OB \rangle &&= \langle OC, OP \rangle \\ \langle OP, OM \rangle &= \langle OP, OC \rangle &&= \langle OA, OM \rangle \end{aligned}$$

The four points $M, N, O, P$ are each the orthocentre of the triangle formed by the other three, hence de Casteljau naming this the 'orthocentric tetragram'. There are three ways to group the four points into two pairs, meaning two orthogonal lines, which defines the points $A, B, C$, cf. Fig. 15 (centre). Conversely, there are four ways to intersect the three angle bisectors of the triangle $\triangle ABC$, which results in the in- and excentres $M, N, O, P$ of the triangle; cf. Fig. 15 (right).

### 10.2. Orthocentric tetragram and nine-point circle

Let us consider the circumcircle to the triangle $\triangle IJK$ with angles $2\alpha + 2\beta + 2\gamma = 180°$ such that $\angle I = \beta + \gamma, \angle J = \gamma + \alpha, \angle K = \alpha + \beta$. With diametrically opposite points $ijk$, we can conclude that the central angle to the circular arcs $\widehat{Kj} = \widehat{kJ}$ equals $2\alpha$, and the inscribed angle $\alpha$ from $J$ and $K$ leads to $\widehat{Ik'} = \widehat{Ij'} = \widehat{jK} = \widehat{Jk}$.

We now see that $Ii', Jj', Kk'$ are the altitudes of the triangle and that $IH = jK = kJ, JH = kI = iK, KH = iJ = jI$ with orthocentre $H$. As the diagonals of the three parallelograms $HIkJ, HJiK, HKjI$ form midlines, a similarity with centre $H$ and ratio $1:2$ maps the feet $A, B, C$ of the altitudes onto the points $i', j', k'$. The triangle $\triangle ABC$ is called orthotriangle with angles $2\alpha, 2\beta, 2\gamma$ at $A, B, C$ respectively.

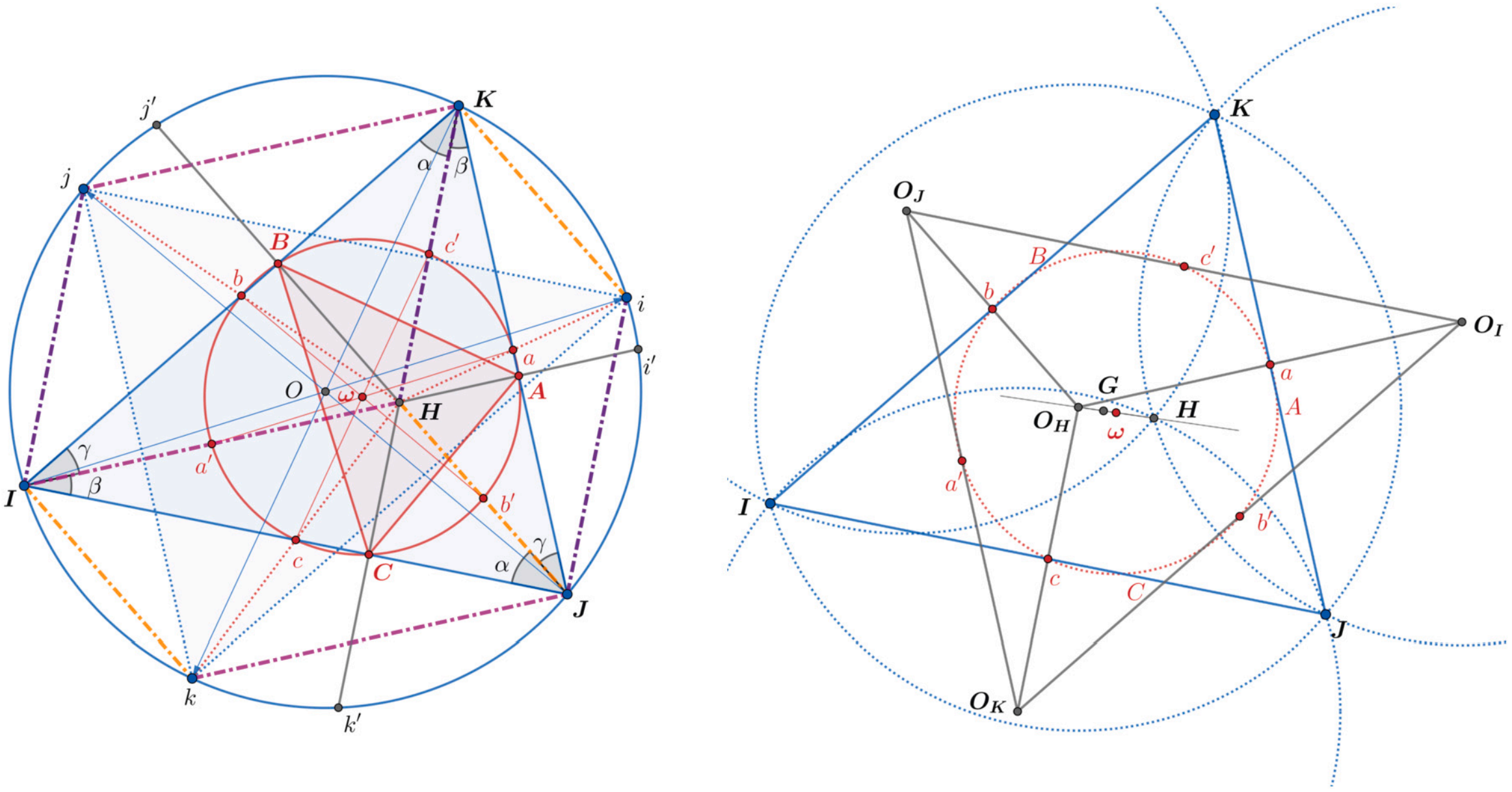

**Fig. 16.** Left: (blue) triangle $IJK$ with (red) nine-point circle $ABC$ Right: Euler line with $O_H, G, \omega, H$. (For interpretation of the colours in the figure(s), the reader is referred to the web version of this article.)

We find the Euler line (Fig. 16) of the triangle $\triangle IJK$ with collinear circumcentre $O_H$, centroid $G$, orthocentre $H$, and centre $\omega$ of the nine-point circle by the equations

$$
\begin{aligned}
4\omega &= \quad H + I + J + K && \\
2\,O_H &= -H + I + J + K = 3G - H && \text{(Feuerbach form)} \\
2\omega + O_H &= \quad 3G && \text{(Euler form)}
\end{aligned}
$$

### *10.3. Dualities*

The notion of duality is related to projective geometry, where we study phenomena from the viewpoint of both points and hyperplanes. In the plane, this duality compares point configurations with line configurations. In planar metric geometry, de Casteljau introduces different aspects on duality, which then lead to observations of growing complexity.

- distance – angle

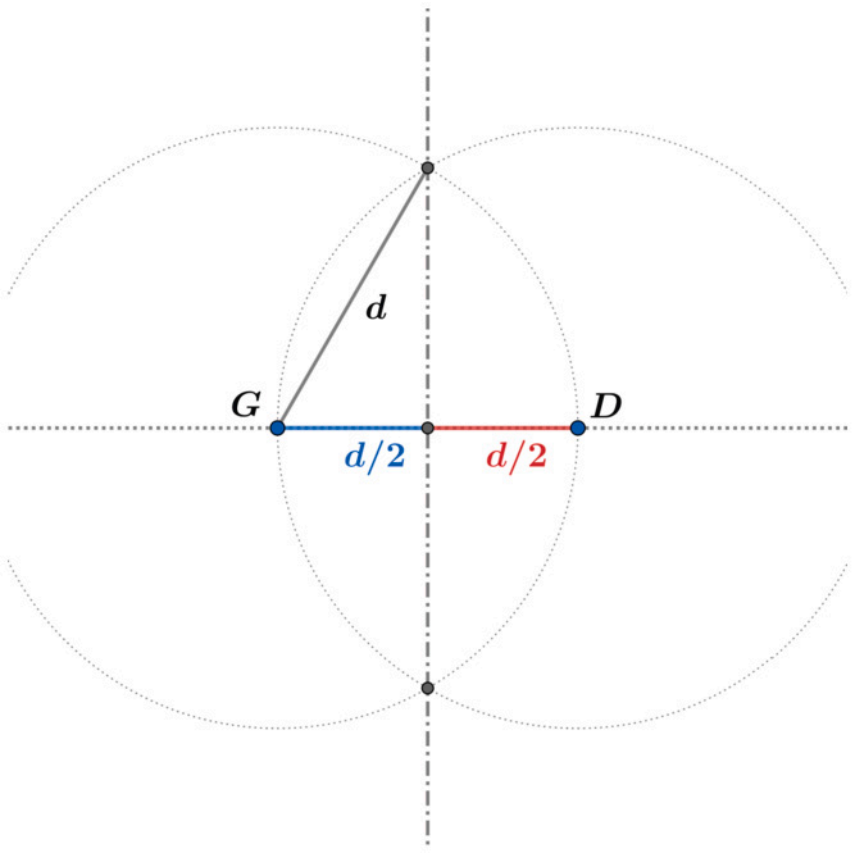

Two points define their *distance*, which is cut into two halves by the perpendicular bisector.

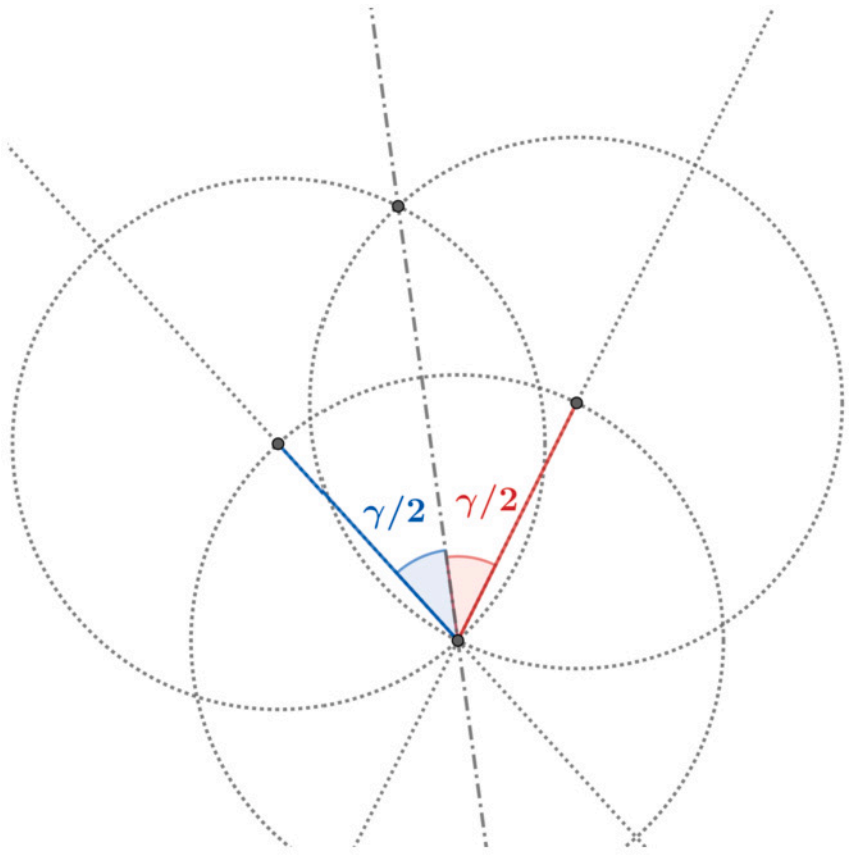

Two lines define their enclosed *angle*, which is cut into two halves by the angle bisector.

- point symmetry – line symmetry
  Assume two points $G, D$ seen from a third one $P$, all three points define a circle $K$.

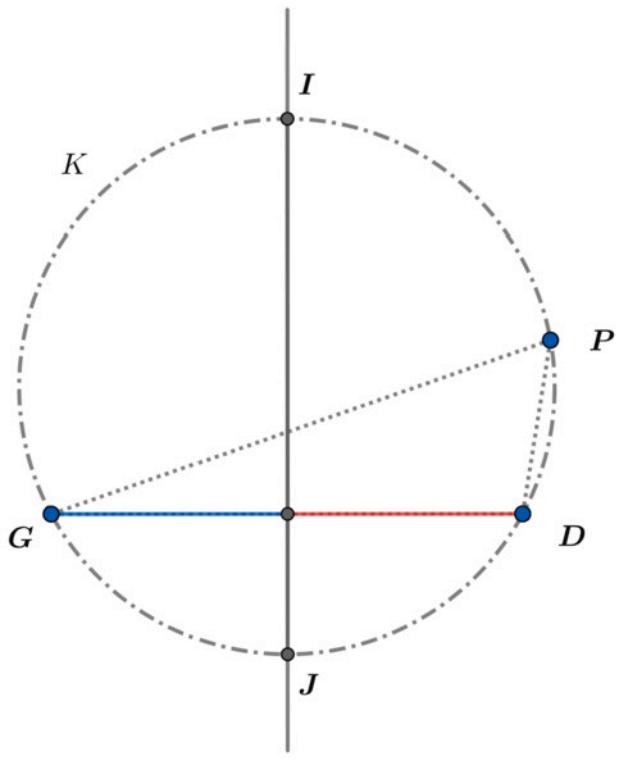

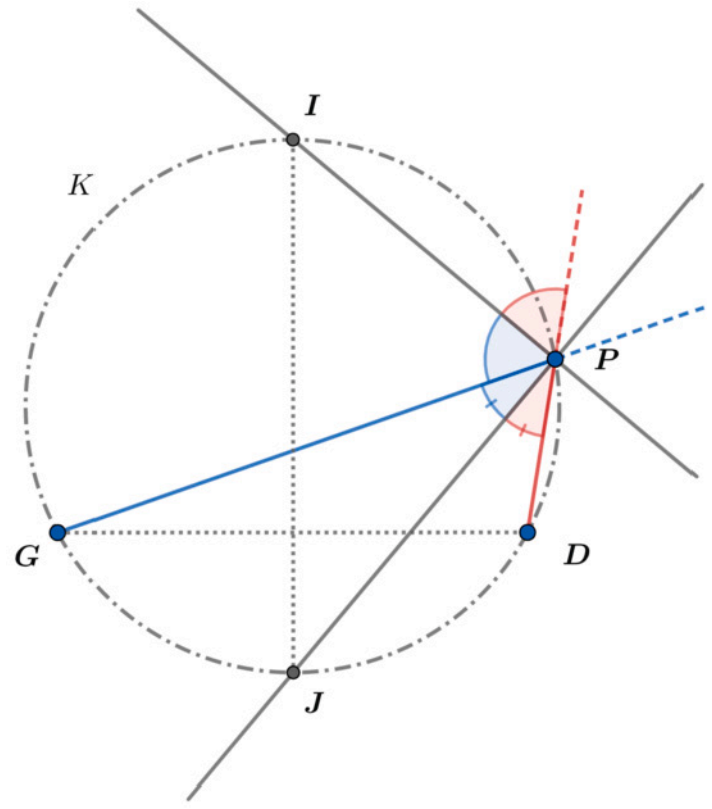

The symmetry axis of the points $G, D$ is their perpendicular bisector, which intersects the circle $K$ in $I, J$.
The line $IJ$ divides the distance $GD$ into two half-distances.

The symmetry axes of the lines $GP, PD$ are their angle bisectors, which intersect the circle $K$ in $I, J$.
The lines $PI, PJ$ divide the angle $GPD$ into two half-angles each.

- tangent – normal
  Conics with same foci $G, D$ ($GD = 2c$) satisfy the equation

$$\frac{x^2}{c^2-\lambda}+\frac{y^2}{\lambda}=1.$$

$\lambda > 0$ results in an ellipse, $-c^2 < \lambda < 0$ is a hyperbola, $\lambda < -c^2$ represents imaginary forms.

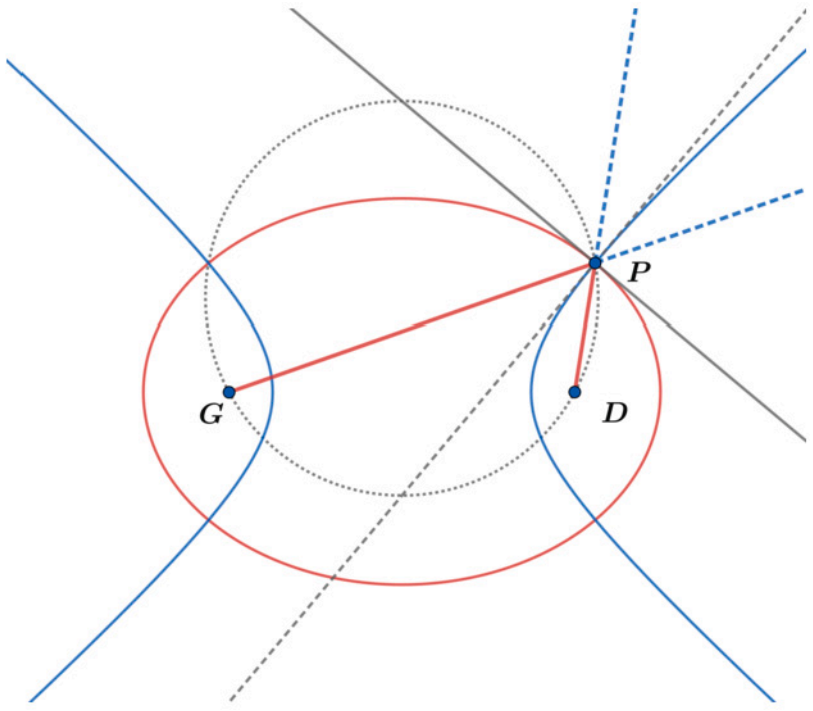

The tangent to the ellipse in $P$ is the normal of the confocal hyperbola in $P$, the ellipse's normal is the hyperbola's tangent.

The tangent to the ellipse in $P$ is the normal of the confocal hyperbola in $P$.

The tangent to the hyperbola in $P$ is the normal of the confocal ellipse in $P$.

- Newton – Huygens

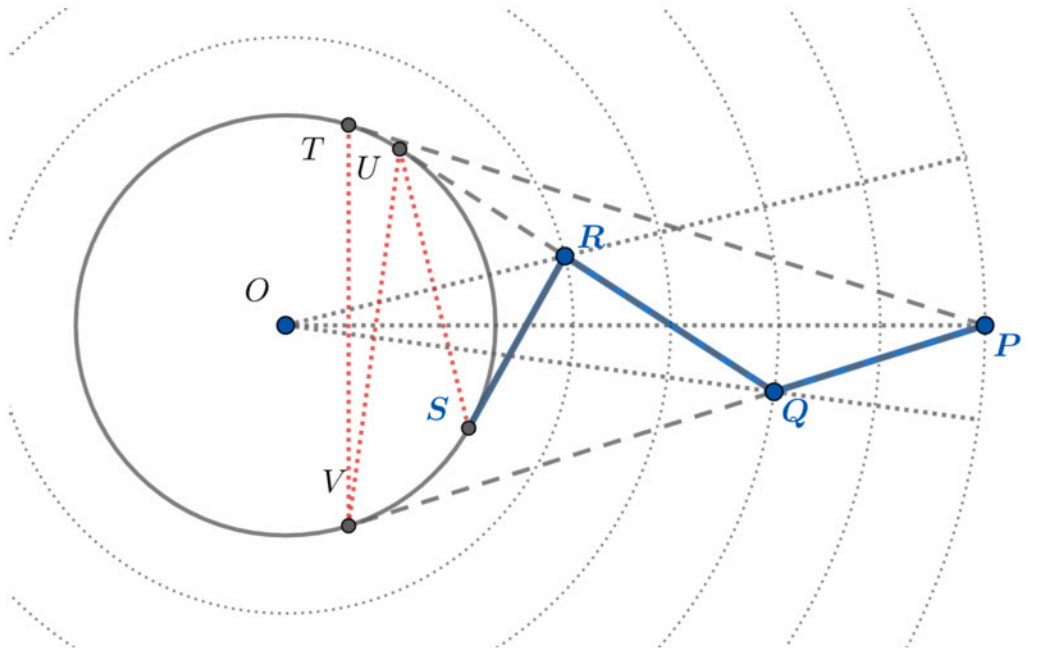

The tangential distance of $P$ to the circle $K$ is equal to the reflection line $PQRS$.

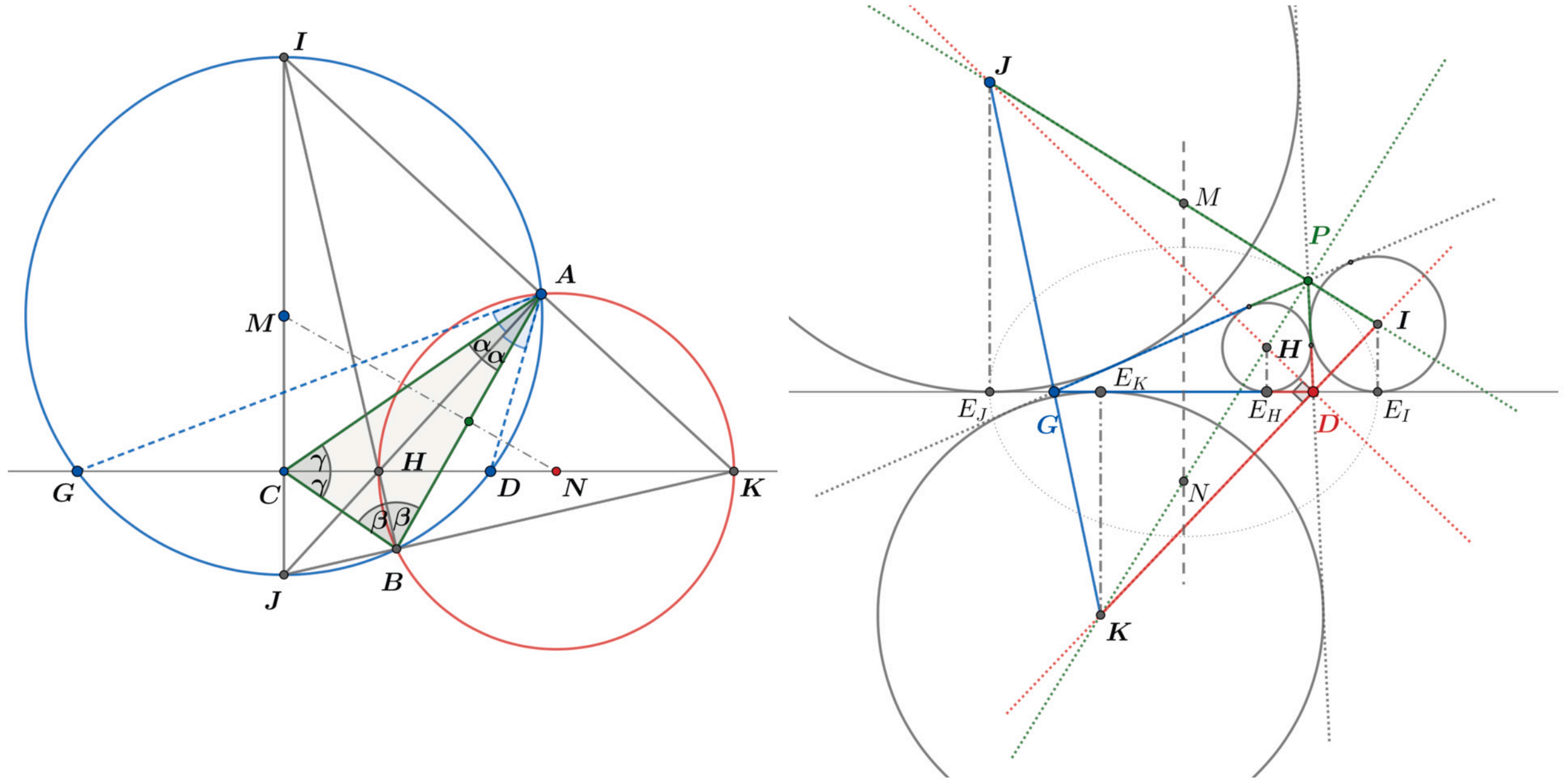

**Fig. 17.** Left: Apollonian circles: (blue) elliptic circle passing through $G, D$, (red) hyperbolic circle with limiting points $G, D$ Right: Ellipse with string length $GP + PD$, hyperbola with $GP - PD$.

| | |
|---|---|
| Light is made out of particles (*corpuscles*), which follow a straight motion, and are thus reflected as rays. | Light is made out of waves, every point of a wavefront serves as a source of secondary wavelets, the new wave front is tangential to all the secondary wavelets. |

### *10.4. Apollonian circles and confocal conics*

After meeting in person for the first time in Oberwolfach in 1992, de Casteljau sent mathematical notes to Boehm, which the latter could not read due to a lack of French language skills. Among them is a metric construction of the duality between Apollonian circles and confocal conics.[30]

**Theorem 1.** *A family of confocal conics with foci $G, D$ and a pencil of (orthogonal) circles with two base points $G, D$ follow a similar description:*

| $\dfrac{x}{\frac{1}{t}+t} = \dfrac{y}{\frac{1}{T}-T} = \cdots$ | $\cdots = \dfrac{c}{\frac{T}{t}+\frac{t}{T}}.$ | $\cdots = \dfrac{2c}{(\frac{1}{t}-t)(\frac{1}{T}+T)}.$ |
|---|---|---|
| $t = const$ : | *elliptic circles passing through $G, D$* | *ellipses with foci $G, D$* |
| $T = const$ : | *hyperbolic circles with limiting points $G, D$* | *hyperbolas with foci $G, D$* |

*Proof with Apollonian circles*

Let $G, D$ form a line, we assume $GD = 2c$ and find $C$ at the perpendicular bisector. The angular bisectors of a point $A = (x, y)$ with $G$ and $D$ meet line $GD$ in $H, K$, and meet its perpendicular bisector in $I, J$. We complete the picture with $B$ on line $IH$, and recognise $HIJK$ as the orthocentric tetragram to the $\triangle ABC$, cf. Fig. 17.

The (elliptic or hyperbolic) Apollonian circles follow the equations

$$x^2 + y^2 - \frac{2cy}{\tan 2\varphi} - c^2 = 0 \quad (elliptic) \qquad x^2 - \frac{2cx}{\cos 2\vartheta} + y^2 + c^2 = 0 \quad (hyperbolic)$$

By setting

$$CH = cT, \quad CK = c/T, \quad CI = c/t, \quad CJ = -ct,$$

[30] de Casteljau, personal communication to W. Boehm, 23 November 1992 and Christmas 1994.

de Casteljau turns the circles into

$$x^2+y^2-c^2=\left(\frac{1}{t}-t\right)cy \quad (elliptic) \qquad x^2+y^2+c^2=\left(\frac{1}{T}+T\right)cx \quad (hyperbolic)$$

which leads to the above stated relation.

*Proof with confocal conics*

If two in $D$ orthogonal lines meet the vertical tangents of the conic in $I, J$, then the tangents from $G, D$ to the circle around $I$ meet at a point $P$ on $IJ$, which lies on the conic with foci $G, D$. It is

$$GP+PD= \ GE_I+DE_I=GE_J+DE_J \ =E_IE_J \quad (ellipse)$$

$$GP-PD=GE_H-DE_H=GE_K-DE_K=E_HE_K \quad (hyperbola)$$

where $H, K$ are the intersections of the angular bisectors at $P$ with the orthogonal lines of $D$.

We again find the orthocentric tetragram $HIJK$ by the centres of the circles inscribed or circumscribed to $\triangle PDG$. De Casteljau lists the tangent lengths from a point $PDG$ to one of the circles in the following table:

| $from$ | $H$ | $I$ | $J$ | $K$ | |
|---|---|---|---|---|---|
| $to:$ | | $G$ | $D$ | $P$ | $=d+g+c$ |
| $to:$ | $G$ | | $P$ | $D$ | $=-d+g+c$ |
| $to:$ | $D$ | $P$ | | $G$ | $=d-g+c$ |
| $to:$ | $P$ | $D$ | $G$ | | $=d+g-c$ |

with $GD=2c, GP=2g, DP=2d$.

Confocal conics follow the equations

$$x^2-\frac{y^2}{\sin^2 2\varphi}=\frac{c^2}{\cos^2 2\varphi} \quad (ellipses) \qquad \frac{x^2}{\cos^2 2\vartheta}-\frac{y^2}{\sin^2 2\vartheta}=c^2 \quad (hyperbolas)$$

A parametrisation of the angles $\varphi, \vartheta$ from $G, D$,

$$\cos 2\varphi=\frac{1-t^2}{1+t^2}, \quad \sin 2\varphi=\frac{2t}{1+t^2}, \quad \cos 2\vartheta=\frac{2T}{1+T^2}, \quad \sin 2\vartheta=\frac{1-T^2}{1+T^2}$$

allows de Casteljau to turn the confocal equations into

$$\left(\frac{x}{1/t+t}\right)^2+\left(\frac{y}{2}\right)^2=\left(\frac{c}{1/t-t}\right)^2 \quad (ellipses) \qquad \left(\frac{x}{2}\right)^2-\left(\frac{y}{1/T-T}\right)^2=\left(\frac{c}{1/T+T}\right)^2 \quad (hyperbolas)$$

which leads to the above stated relation.

**Addendum**

The mapping from the meeting point $A$ of the Apollonian circles onto the intersection of $MN$ with $AB$ is a Joukovsky transformation, which here appears without using complex numbers.[31]

## 11. Bifocal coordinates

De Casteljau already worked in the 1970s on systems with bifocal coordinates, which have supported his works on optics in section 12. We want to examine the building blocks.

Let us look at a circle centered at $T$ from two points $A, B$ at their respective angles $2\alpha, 2\beta$. At its tangents, the circle's radius $r$ forms a right-angled triangle with $TA$ as well as with $TB$, cf. Fig. 18, and we recognise the sine law:

$$r=AT\sin\alpha=BT\sin\beta$$

If we assume two lines rotating around $O$ and $I$, and their angular speed ratio is $n:p$, then their intersection generates the following curves:

[31] de Casteljau, personal communication to W. Boehm, Christmas 1994.

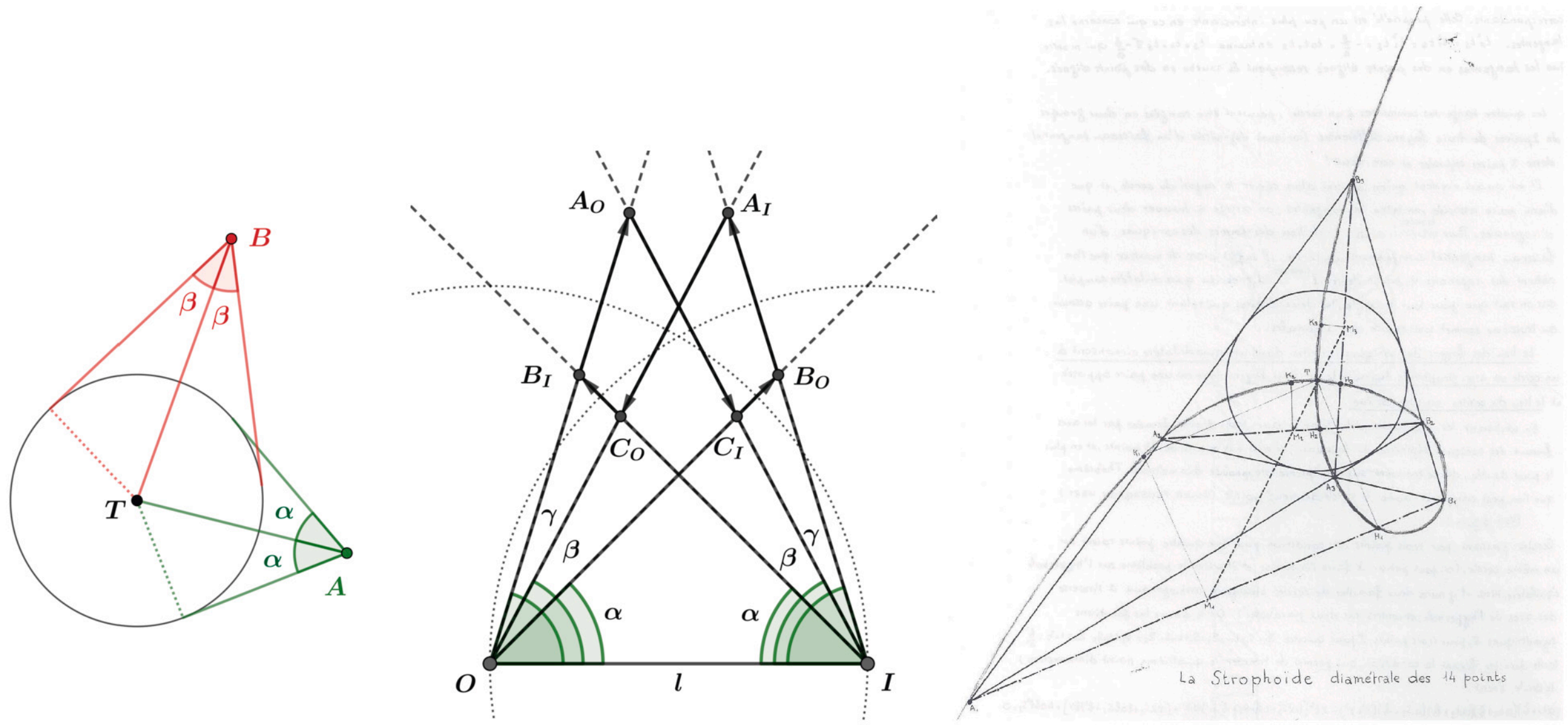

**Fig. 18.** Left: A circle seen from two angles. Centre: Intersections of inverted rays. Right: 14 points of a strophoid; drawing by de Faget de Casteljau (1995b).

| | |
|---|---|
| 0:1 | a line through $O$ |
| 1:1 | a circle with diameter $OI$ |
| -1:1 | a rectangular hyperbola |
| 1:2 | a right strophoid |
| 2:3 | a limaçon of Pascal |
| ... | |

An inversion with centre $O$ and tangential power $OI^2 = l^2$ conserves the angles $\alpha, \beta, \gamma$ with $\alpha + \beta + \gamma = 180°$ and turns the case of $n : p$ into $n : (n-p)$. Conversely, the inversion with centre $I$ leads to the ratio $(p-n) : n$. Fig. 18 shows that six mutual inversions in $(O, l)$ and in $(I, l)$ lead to an identity, e.g., starting at $C_O$:

$$C_O \xrightarrow{I} B_I \xrightarrow{O} A_O \xrightarrow{I} C_I \xrightarrow{O} B_O \xrightarrow{I} A_I \xrightarrow{O} C_O$$

The inverse of a circle is a line ($n = p = 1$, $n - p = 0$), the inverse of a rectangular hyperbola is a right strophoid ($n = 1 = -p$, $n - p = 2$).

Another remarkable finding in de Casteljau's paper (2001b) is a generalisation of the Feuerbach nine-point circle (cf. section 10). Fig. 18 reproduces a drawing by de Casteljau, which shows a circle inscribed in a quadrilateral and a strophoid passing through the circle's centre and the quadrilaterals' intersection points.

**Theorem 2.** ***(De Casteljau's 14-point strophoid)***

*The following 14 points lie on a right strophoid:*

- *the centre $T$ of a circle (counts twice due to the two tangents to the strophoid)*
- *six intersections of a complete quadrilateral that inscribes this circle, at $A_1, A_2, A_3, B_1, B_2, B_3$*
- *three feet $H_1, H_2, H_3$ of the altitudes of the triangles $A_1B_1T$, $B_1B_2T$, $A_3B_3T$*
- *three intersections $K_1, K_2, K_3$ of a parallel of $A_iB_i$ through $T$ with its perpendicular bisector*

De Casteljau worked out a variety of further properties in his last publication (de Faget de Casteljau, 2001b), all based on metric properties. Some of his observations were recently published (independently) by Hellmuth Stachel (2015a,b), who mainly used projective geometry. He also points out further unexplored topics, e.g., the intersection of two strophoids in seven points, or the construction of the curve based on some given elements.

Seen from the metric geometry perspective of Paul de Casteljau, the strophoid with its bifocal coordinates builds directly on the unifocal circle, while conics are seen here as a special case of trifocal Cartesian ovals with a third focus at infinity (cf. the definition by focus and directrix) – quite different to the introduction of conics as curves of second degree in projective geometry.

## 12. Geometric optics

*Then we can admire the extreme elegance of the relations in which $\sin^2\alpha$ is replaced by its polar form $\sin\alpha\sin\beta$ (the same applies for $\cos^2\beta, \cos^2\alpha, \sin^2\beta$).* — de Faget de Casteljau (1995b)

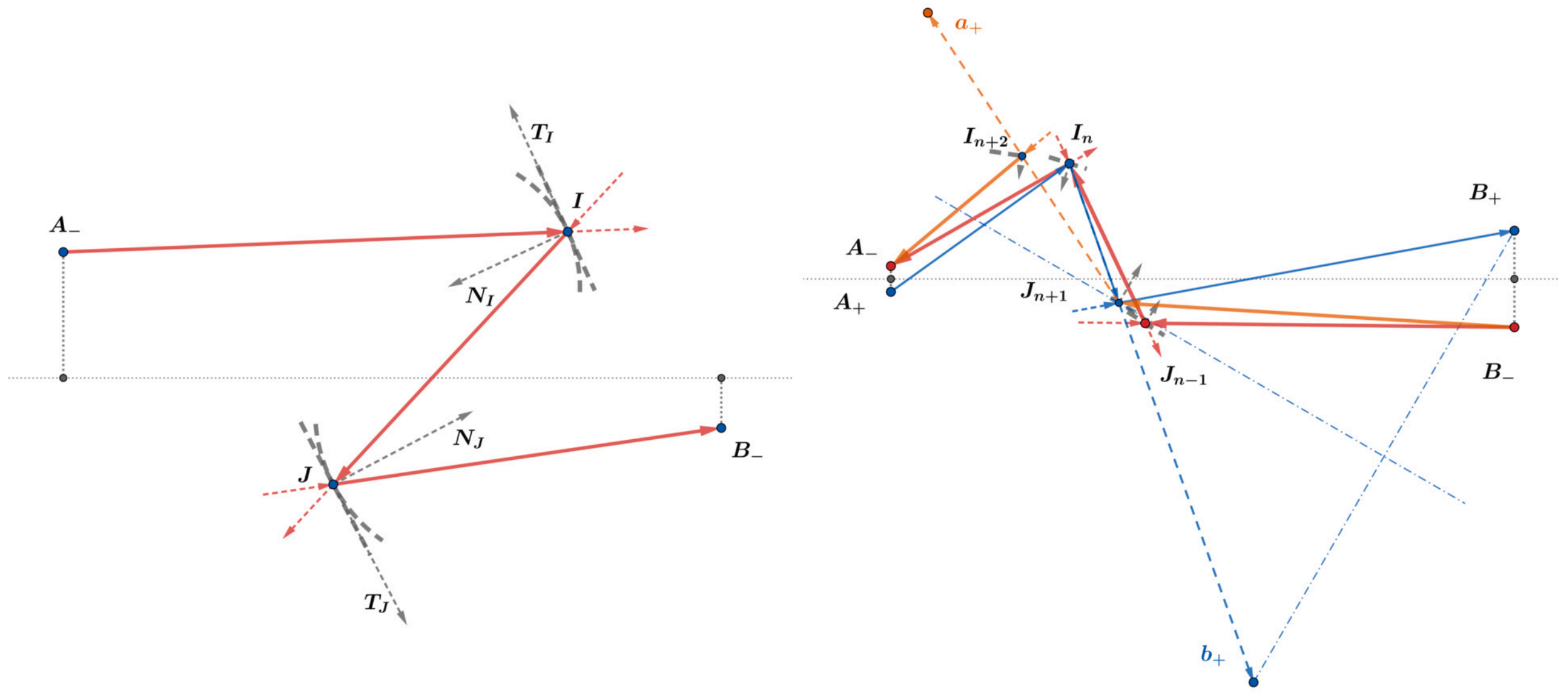

**Fig. 19.** How to define the surfaces, that reflect the light from $A_+$ to $B_+$?

As a well-educated amateur astronomer, de Casteljau possessed a telescope of the highest quality for that time. He was interested in geometric optics and in 1974 developed a subdivision of an ellipse into two ellipses that each share a focus with the outer ellipse, cf. Fig. 20 (left). In *Le Lissage*, he takes as an example the calculation of conjugate mirrors (*aplanats* in geometrical optics), for which he offers two approaches, a precise one via the integration of a differential equation, and an approximate one based on a recurrence relation. Let us first build a common basis:

In the case of a telescope, a ray coming from a star $B_+$ passes the optical system via $J$ and $I$ and arrives at $A_+$, cf. Fig. 19. Let us assume that a second star $B_-$ is emanating light to $A_-$. Then the question occurs as to how the mirrors $I, J$ are defined.[32]

## 12.1. Conjugate mirrors: the recursive approach

We follow a ray that starts in one focus and is then reflected at a conic's tangent such that it arrives at the other focus. We can observe: *If two conics $C_0, C_1$ share one focus $E_0 = E_1$ as well as a tangent at a common point $D$, then the remaining two foci $F_0, F_1$ form a line with $D$.*

De Casteljau applies this to an outer (blue) ellipse and two (red and green) inner ellipses with a side condition: the second focus $J$ of the red ellipse shall lie on the green ellipse, and conversely the second focus $I$ of the green ellipse shall lie on the red one. Then the distances in the (blue) ellipse equal, $PD + DQ = u + (s + t) = u + (w + v) = (u + w) + v = PE + EQ$, where $u + w$ is the chord length of the (red) conic tangent in $E$, and $v + w$ is the chord length of the (green) conic tangent in $D$, cf. Fig. 20.

This conical configuration helps us to now iteratively construct points $I_k, J_k$ of the unknown reflective surface. A first light ray is refracted to the line $B_- J_{n-1} I_n A_-$. We want to travel backwards from $A_+$ to $B_+$, knowing that the distance of the light ray travelling back must be identical to the first ray (Fermat's principle): at $I_n$, the light is reflected along the $I_n$-tangent to some virtual point $b_+$. If reflected at some $J_{n+1}$, the ray returns to $B_+$, so the distances $I_n b_+$ and $I_n J_{n+1} B_+$ are the same. We apply the above duality of bisector and angle bisector: the section bisector $b_+ B_+$ is the angle bisector at a reflection. So, the intersection of the bisector with the ray from $I_n$ to $b_+$ results in the reflection point $J_{n+1}$. The procedure repeats from $B_-$ via $J_{n+1}$ to find $I_{n+2}$ etc.

We call $u_n, w_n, v_n$ the lengths of the reflected line segments $A_+ I_n, I_n J_{n-1}, J_{n-1} B_+$, their angles with the optical axis being $2\alpha_n, 2\gamma_n, 2\beta_n$. Then, the normals have angles $\psi_n = \gamma_n + \alpha_n$ at $I_n$ and $\psi_{n-1} = \gamma_n + \beta_n$ at $J_{n-1}$.

The next trajectory $A_- I_n J_{n+1} B_-$ then constructs $u_{n+1}, w_{n+1}, v_{n+1}, \alpha_{n+1}, \gamma_{n+1}, \beta_{n+1}$ by first identifying $u_{n+1}, \alpha_{n+1}$, then $\gamma_{n+1}$ and point $b_-$ with the perpendicular bisector of $b_- B_-$, thus $\psi_{n+1}$ and $\beta_{n+1}$, and finally the intersection with the line $a_- b_-$ delivers $v_{n+1}$ and ultimately $w_{n+1}$.

## 12.2. Conjugate mirrors: the exact approach

The circle around $O$ with radius $r$ cuts $AB$ in the ratio $\sin 2\beta : \sin 2\alpha$ (cf. Fig. 20); we can thus introduce a value $q$ defined by

$$OA = \frac{1-q}{1+q} r, \qquad OB = \frac{r^2}{OA} = \frac{1+q}{1-q} r$$

[32] A solid understanding of astronomic optics including history and geometric optics can be obtained from Lemaître's book (2009), a deeper understanding of two-mirror systems is developed by Willstrop and Lynden-Bell (2003).

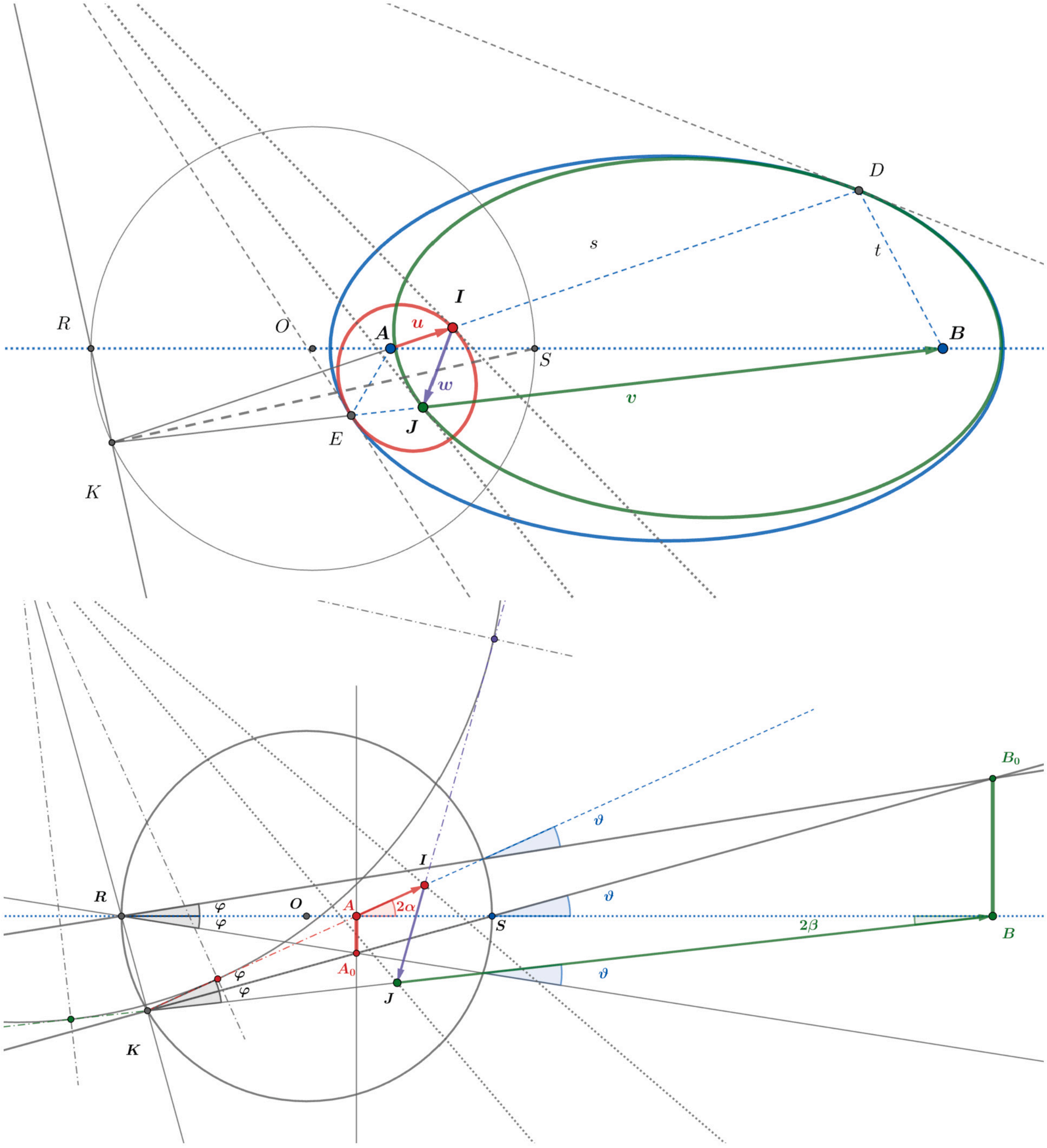

**Fig. 20.** Top: The outer conic is bitangent to the two inner conics. Bottom: Bipolar coordinates $\alpha, \beta$ and $\varphi, \vartheta$.

By similarity considerations, e.g., at $K$ for $\varphi$ or at $O$ for $\vartheta$, we see that

$$2\vartheta = \alpha + \beta, \quad 2\varphi = \alpha - \beta, \quad \alpha = \vartheta + \varphi, \quad \beta = \vartheta - \varphi$$

We set $AB = 2c$ and the length of the optical trajectory $AIJB = 2a$ and remember that the surfaces searched for are defined by polar coordinates

$$u = AI, \qquad v = JB.$$

We then can write

$$(IJ)^2 = w^2 = (2a - u - v)^2 = (u\cos\alpha + v\cos\beta - 2c)^2 + (u\sin\alpha + v\sin\beta)^2$$

which turns into

$$(a^2 - c^2) - u(a - c\cos\alpha) - v(a - c\cos\beta) + uv\sin^2(\frac{\alpha-\beta}{2})$$

Subsequently, we obtain the following relations for the four optical cases (the Gregory and Cassegrain two-mirror configurations only differ by a sign), based on Fermat's condition[33]:

$$\begin{array}{lcccccccc}
RR: & \left[\frac{1}{u} - \frac{\sin^2\alpha}{a+c} - \frac{\cos^2\alpha}{a-c}\right] & \cdot & \left[\frac{1}{v} - \frac{\sin^2\beta}{a+c} - \frac{\cos^2\beta}{a-c}\right] & = & \left[\pm\frac{\sin\alpha\sin\beta}{a+c} + \frac{\cos\alpha\cos\beta}{a-c}\right]^2 & & \\
VR: & \left[\frac{1}{u} + \frac{\cos^2\alpha}{a+c} + \frac{\sin^2\alpha}{a-c}\right] & \cdot & \left[\frac{1}{v} - \frac{\sin^2\beta}{a+c} - \frac{\cos^2\beta}{a-c}\right] & + & \left[\pm\frac{\cos\alpha\sin\beta}{a+c} - \frac{\sin\alpha\cos\beta}{a-c}\right]^2 & = & 0 \\
VV: & \left[\frac{1}{u} + \frac{\cos^2\alpha}{a+c} + \frac{\sin^2\alpha}{a-c}\right] & \cdot & \left[\frac{1}{v} + \frac{\cos^2\beta}{a+c} + \frac{\sin^2\beta}{a-c}\right] & = & \left[\frac{\cos\alpha\cos\beta}{a+c} \pm \frac{\sin\alpha\sin\beta}{a-c}\right]^2 & & \\
RV: & \left[\frac{1}{u} - \frac{\sin^2\alpha}{a+c} - \frac{\cos^2\alpha}{a-c}\right] & \cdot & \left[\frac{1}{v} + \frac{\cos^2\beta}{a+c} + \frac{\sin^2\beta}{a-c}\right] & + & \left[\frac{\sin\alpha\cos\beta}{a+c} \mp \frac{\cos\alpha\sin\beta}{a-c}\right]^2 & = & 0
\end{array} \tag{6}$$

A less intuitive change of the variables

$$\frac{\cos\vartheta}{\cos\varphi} = t, \qquad z = \frac{1}{v} - \frac{a - c\cos\beta}{a^2 - c^2}$$

and the remark of $t\sin\varphi = q\sin\vartheta$ lead de Casteljau to the integrable differential equations of the primary mirror at $J$:

$$\begin{array}{llll}
RR: & \frac{dz}{z} = \frac{(aq-ct)(t^2+q)}{(a-ct)t(t^2-q^2)}dt & VR: & \frac{dz}{z} = \frac{(a-ct)(t^2+q)}{(aq-ct)t(t^2-1)}dt \\
RV: & \frac{dz}{z} = \frac{(a+ct)(t^2+q)}{(aq+ct)t(t^2-1)}dt & VV: & \frac{dz}{z} = \frac{(aq+ct)(t^2+q)}{(a+ct)t(t^2-q^2)}dt
\end{array}$$

These equations can be integrated, once decomposed into the following elements

$$\frac{1}{v} - \frac{a - c\,\cos\beta}{a^2 - c^2} = \frac{K}{t}\,(t - a/c)^{(1-q)\frac{a^2+qc^2}{a^2-q^2c^2}}\,(t+q)^{\frac{1+q}{2}\,\frac{a+c}{a+cq}}\,(t-q)^{\frac{1+q}{2}\,\frac{a-c}{a-qc}} \tag{7}$$

The secondary mirror at $I$ is derived from $z = \frac{1}{u} - \frac{a-c\cos\alpha}{a^2-c^2}$, more details of which can be found in (de Faget de Casteljau, 1990, ch. 11). Fig. 21 shows that these conjugate mirrors are no conics.[34]

De Casteljau argues in *Le Lissage* that both approaches at large – the exact algebraic solution and the recursive approximation – have their advantages. While he develops a recurrence relation for the computation of an ellipsoid's geodesics, he rejects an algebraic approach. He argues inversely for the lines of curvature of a quadric, where he chooses the exact over the approximate. In the above-shown case of the optical systems, both approaches are useful, and the approximation helps to solve an equation in infinity.

### *12.3. The backdrop*

De Casteljau draws our attention to the circle around $O$, which intersects the optical axis in $R$ and $S$. The circle also contains $K$, where $AI$ and $JB$ intersect under the angle of $2\varphi$. The angular bisector from $K$ meets the polar line of $B$ at $A_0$ and the polar of $A_0$ at $B_0$. We now recognise, that $R$ sees the image $AA_0$ under the angle of $\varphi$, and $S$ sees the object $BB_0$ under the angle of $\vartheta$. It becomes obvious that interchanging $R$ with $S$ and $\varphi$ with $\vartheta$ leads from the Gregory case to the Cassegrain case. It is also remarkable, that all reflection or refraction happens within the circle as the 'experiment chamber', which also applies for most of the calculations.

De Casteljau invites us to think beyond the Abbe sine condition $\sin 2\alpha : \sin 2\beta = const$, which he related to $\tan\vartheta : \tan\varphi = konst$ with $K$ describing a circle. In Fig. 22, he relates the position of $K$ to other geometric optics questions. He entitles his observations *la toile de fond* (the backdrop), which describes the shape on which $K$ is located: line, circle, or conic.

$$\begin{array}{lcccl}
\text{Gauss}: & \tan 2\alpha : \tan 2\beta = const & & \sin 2\vartheta : \sin 2\varphi = konst & K \text{ on a line, plane} \\
 & & \swarrow & & \\
\text{Abbe}: & \sin 2\alpha : \sin 2\beta = const & & \tan\vartheta : \tan\varphi = konst & K \text{ on a circle} \\
 & & \swarrow & & \\
 & \tan\alpha : \tan\beta = const & & \sin\vartheta : \sin\varphi = konst & K \text{ on a conic} \\
 & & \swarrow & & \\
\text{Herschel}: & \sin\alpha : \sin\beta = const & & &
\end{array}$$

[33] James Gregory designed a telescope in 1663 where the second mirror $I$ is concave (elliptic) outside the focus of $J$; in 1672 Laurent Cassegrain placed $I$ on a convex (hyperbolic) mirror inside the primary focus of the (parabolic) $J$-mirror. Both constructions were known already to Marin Mersenne (1636, p.61/62), who applied this to the reflection of sound. In 1672, Isaac Newton also presented a telescope with a parabolic primary mirror, only his secondary mirror was planar.

[34] de Casteljau visualised his findings in a letter to Farouki in 1991 (personal communication, 6 September 1991) and at an Oberwolfach talk in 1992, based on his results in *Le Lissage* on pages 103–111.

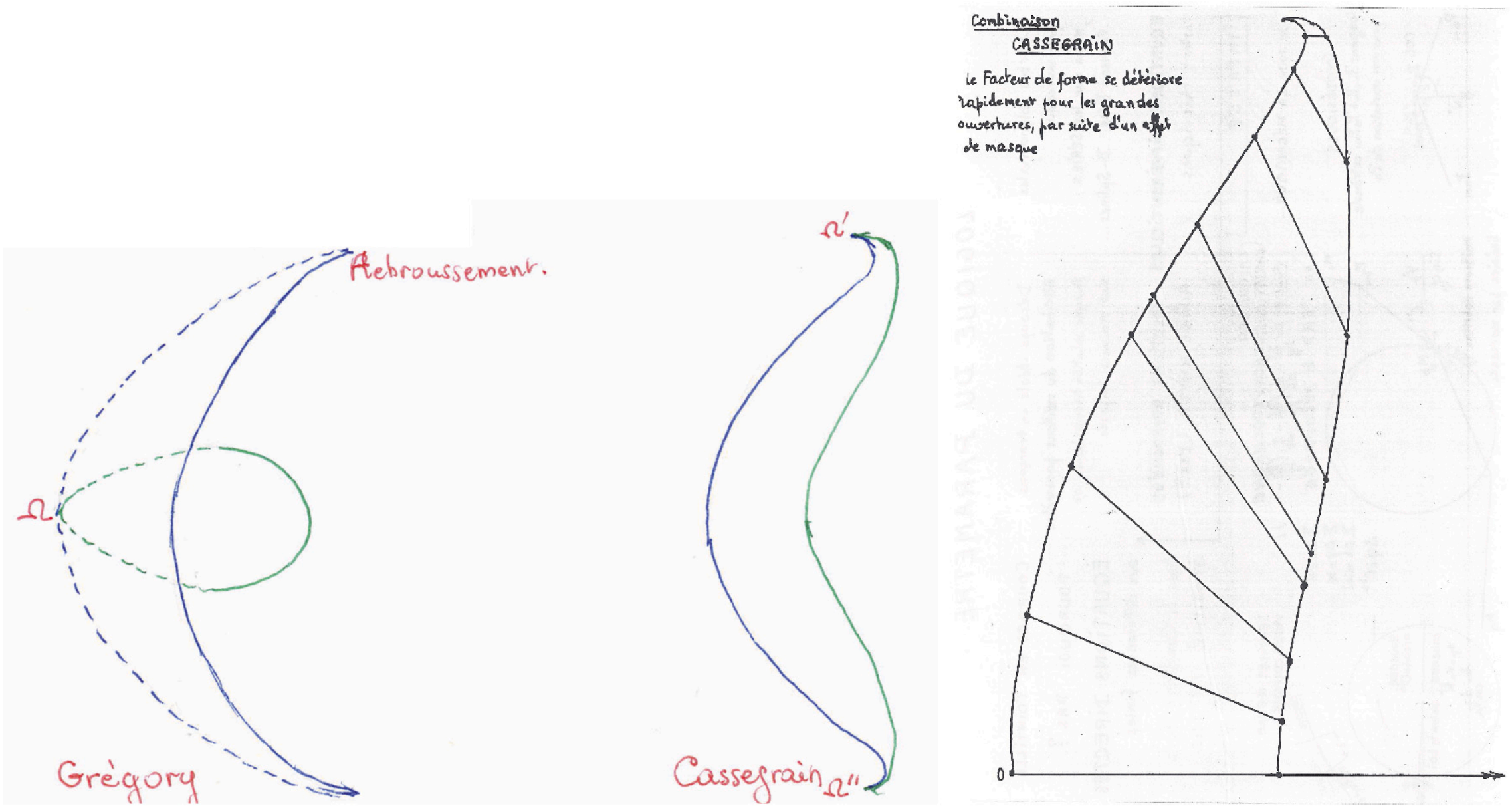

**Fig. 21.** The exact shape of two conjugate mirrors is not conical; its extremities can be calculated with de Casteljau's *ricochet* approach. Figures from a letter to Rida Farouki in 1991 (left) and de Casteljau's talk at Oberwolfach in 1992 (right).

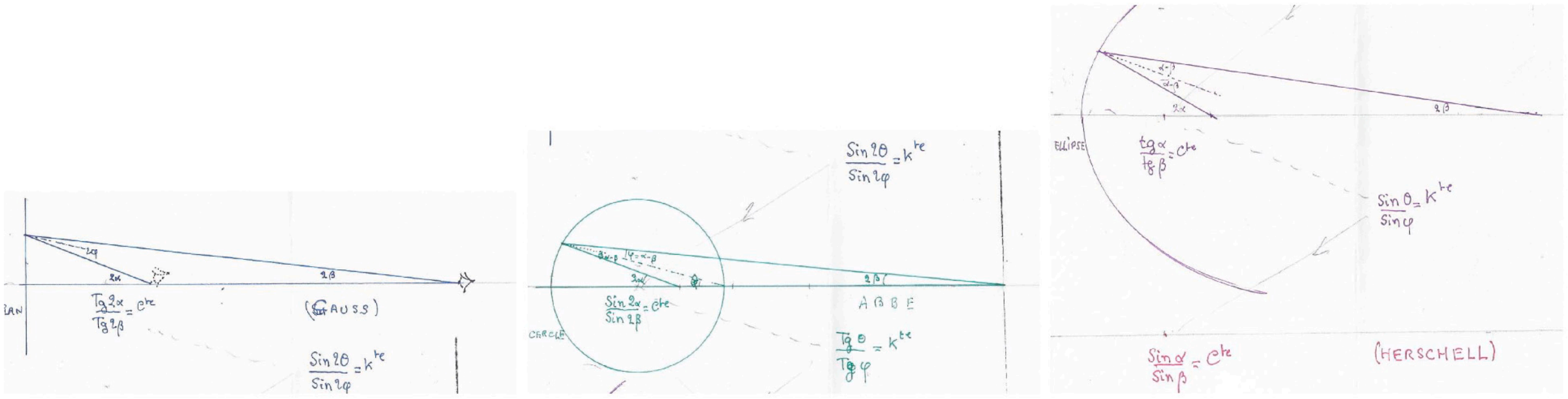

**Fig. 22.** The backdrop; depending on the ratio, point $K$ describes different surfaces; drawings by de Casteljau.

## 13. Regular polygons

*If the effects of age have not diminished my mental faculties too much, I think that this could be the most significant progress since GAUSS* ($z^{2n+1} - 1 = 0$) *concerning the theory of regular polygons with* $2n + 1$ *sides!* — de Casteljau[35]

De Casteljau was fascinated by the geometry of the Byzantine tiles that he once saw in the Alhambra. He worked on meshes in *Les Quaternions*, and there already examined the heptagon (de Faget de Casteljau, 1990, p.155). In 2001, he dedicated a MICAD conference article to the regular 11-gon as an example of a regular polygon that cannot be constructed by ruler, compass and trisector, $n = 11$ being the smallest non-Pierpont prime. He did not cease working on the pentagon, heptagon, nonagon, as further discussed in later letters (de Faget de Casteljau, 2001a).

A remarkable result is his approximation scheme for vertices and diagonals of a regular polygon, which we will now develop for $n = 7$.

### 13.1. Ptolemy's theorem

For the quadrilateral $M, N, O, P$ inscribed in a circle, Ptolemy's theorem states that the product of the diagonals equals the sum of products of opposite sides:

$$ON \cdot MP = OP \cdot MN + OM \cdot PN$$

[35] de Casteljau, personal communication to W. Boehm (translated from French), 5 April 2001.

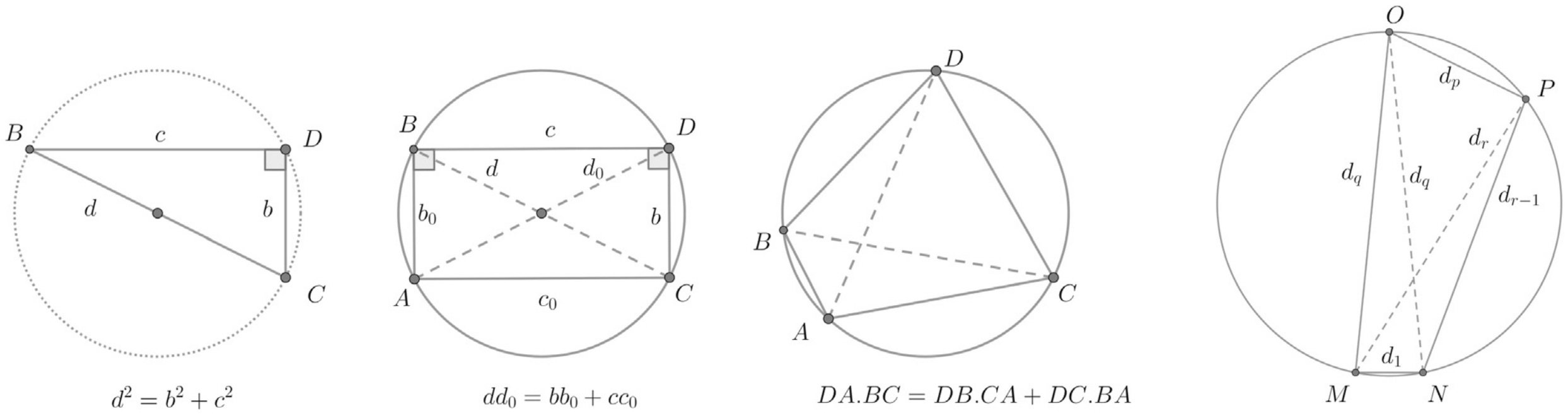

**Fig. 23.** Left: Ptolemy's theorem as a polar form of Pythagoras' theorem. Right: Construction of golden ratios.

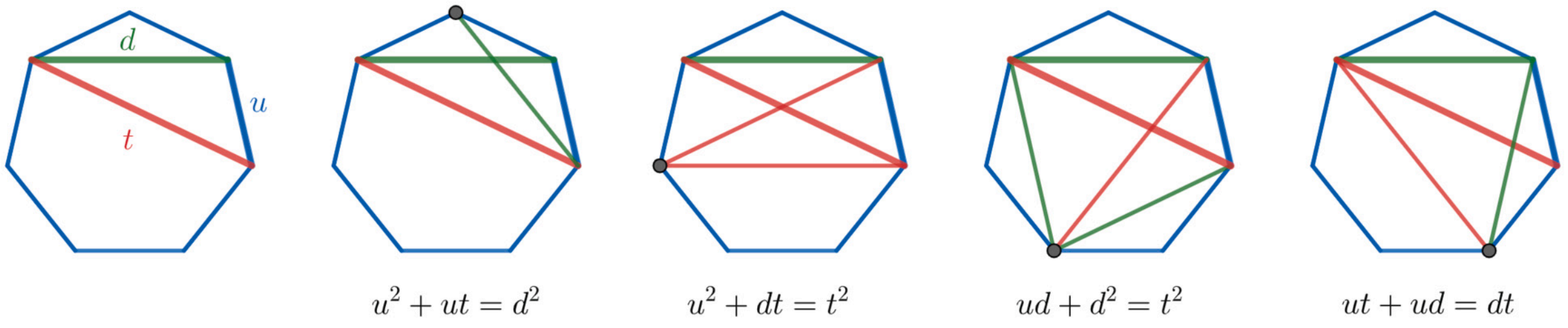

**Fig. 24.** Ptolemy's theorem applied to a heptagon: the base triangle with side lengths $u, d, t$ is extended anti-clockwise by another vertex, which leads to equations I, II, III, 0 from (10).

Let us move $M$ along the polygon. We then find the following Ptolemy relations for chords and diagonals $d_i$, knowing that for symmetry reasons $d_1 = d_{n-1}, \ldots, d_k = d_{n-k}$, cf. Fig. 23:

$$\begin{array}{rcl} d_q d_r &=& d_p d_1 + d_q d_{r-1} \\ d_q d_{r-1} &=& d_{p+1} d_1 + d_q d_{r-2} \\ \vdots && \\ d_q d_2 &=& d_{q-1} d_1 + d_q d_1 \end{array} \quad \rightarrow \quad \begin{array}{rcl} d_q(d_r - d_{r-1}) &=& d_p d_1 \\ d_q(d_{r-1} - d_{r-2}) &=& d_{p+1} d_1 \\ \vdots && \vdots \\ d_q(d_2 - d_1) &=& d_{q-1} d_1 \end{array}$$

The sum of the right column gives $d_q(d_r - d_1) = (\sum_{q-1}^{p} d_i) d_1$ and hence $d_q d_r = (\sum_{q}^{p} d_i) d_1$. This finally leads to the theorem:

**Theorem 3** *(golden ratios). The diagonal lengths $d_i$ in an n-gon, $n = p + q + r$, fulfil the equation*

$$\frac{d_r - d_{r-1}}{d_p} = \frac{d_{r-1} - d_{r-2}}{d_{p+1}} = \cdots = \frac{d_2 - d_1}{d_{q-1}} = \frac{d_1}{d_q} = \frac{d_r}{\sum_{i=q}^{p} d_i} \tag{8}$$

**Addenda**

- For $n = 5$, we get the golden ratio from $\frac{d_1}{d_2} = \frac{d_2}{d_1 + d_2}$, which results in $\Phi = \frac{d_2}{d_1} = \frac{d_1 + d_2}{d_2} = 1 + \frac{1}{\Phi}$.
- Independently, Peter Steinbach (1997, 2000) also worked on golden ratios for $n = 5, 7, 9, 11$, yet without proof, including a relation of the product $\Phi^2 = \Phi + 1$ to the Ammann grid and Penrose tiling; de Faget de Casteljau (2001a, p.26) also mentioned an application to Penrose tilings and quasi-crystals.
- De Casteljau introduced Ptolemy's theorem as a polar form of Pythagoras' theorem, cf. Fig. 23.

### *13.2. The heptagon*

We will now discuss the case of the heptagon and will use de Casteljau's notation: the diagonal lengths will be denoted $u, d, t$ relating to the French integers *un* (one), *deux* (two), *trois* (three). Equation (8) then reads

$$\frac{t - d}{u} = \frac{d - u}{d} = \frac{u}{t} = \frac{t}{u + d + t} \tag{9}$$

Let us look at all Ptolemy relations in a heptagon, as visualised in Fig. 24.

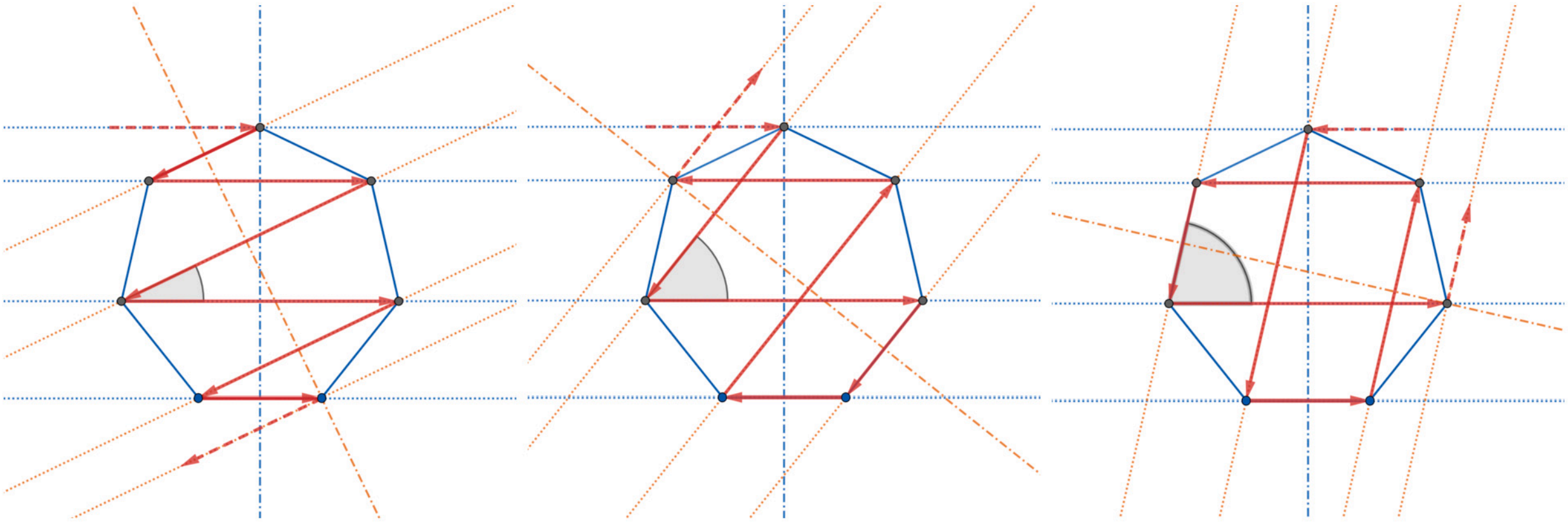

**Fig. 25.** *L'Ange*: Two diagonals of a heptagon form a mesh with inscribed angle $\frac{\pi}{7}, \frac{2\pi}{7}, \frac{3\pi}{7}$ (mirrored from $\frac{4\pi}{7}$), from left to right.

We can find four different quadrilaterals such that $MNOP$ are the vertices of a heptagon. The heptagon shows much more symmetry than polygons with more vertices: the 11-gon already has $20 = 5 \cdot 4$ different quadrilaterals; their geometric structure is discussed by de Faget de Casteljau (2001a).

$$\begin{aligned} d^2 &= u^2 + ut \qquad (I) \\ u^2 &= t^2 - dt \qquad (II) \\ t^2 &= d^2 + ud \qquad (III) \\ dt &= ut + ud \qquad (O) \end{aligned} \tag{10}$$

The four basic relations can be resolved for $u, d, t$:

$$\begin{aligned} \frac{t^2-d^2}{d} &= \frac{dt}{d+t} = u \quad \text{by using } (III) \text{ and } (O) \\ \frac{t^2-u^2}{t} &= \frac{ut}{t-u} = d \quad \text{by using } (II) \text{ and } (O) \\ \frac{d^2-u^2}{u} &= \frac{ud}{d-u} = t \quad \text{by using } (I) \text{ and } (O) \end{aligned}$$

or further

$$(t+d)(t-d)(t+d) - d^2t = t^3 + t^2d - 2d^2t - d^3 = 0$$

$$(t-u)(t+u)(t-u) - t^2u = u^3 - u^2t - 2t^2u + t^3 = 0$$

$$(d-u)(d+u)(d-u) - u^2d = d^3 - d^2u - 2u^2d + u^3 = 0$$

We observe here the equation $x^3 + x^2 - 2x - 1 = 0$ with solutions

$$x_1 = \frac{t}{d}, \qquad x_2 = \frac{-u}{t}, \qquad x_3 = \frac{d}{-u}, \tag{11}$$

which also fulfil

$$x_1x_2x_3 = 1 \quad \text{and} \quad x_1 + x_2 + x_3 = \frac{t}{d} - \frac{u}{t} - d(\frac{1}{d} + \frac{1}{t}) = \frac{t^2 - d^2 - ud - dt}{dt} = -1$$

The geometrical effect of the permutation is depicted in Fig. 25. Obviously, there is a (blue, dashed) symmetry line from the apex to the base of the polygon. Let us now adopt an optical perspective: a light ray enters the heptagon from the apex to one of the other vertices, which defines the (orange) light direction. Depending on the angle of light, we see different reflection patterns – three in the case of the heptagon. The light ray is reflected in an alternating way, mirrored by the (blue) heptagon symmetry axis followed by the (red) normal of the light's direction. The reflection pattern induces a scale on the mesh, which differs in each of the permutations, cf. Fig. 27. Each pair of (blue and red) axes has an inscribed angle, in our case $\frac{\pi}{7}2^k$, with the corresponding mesh intersecting with the polygon's vertices.

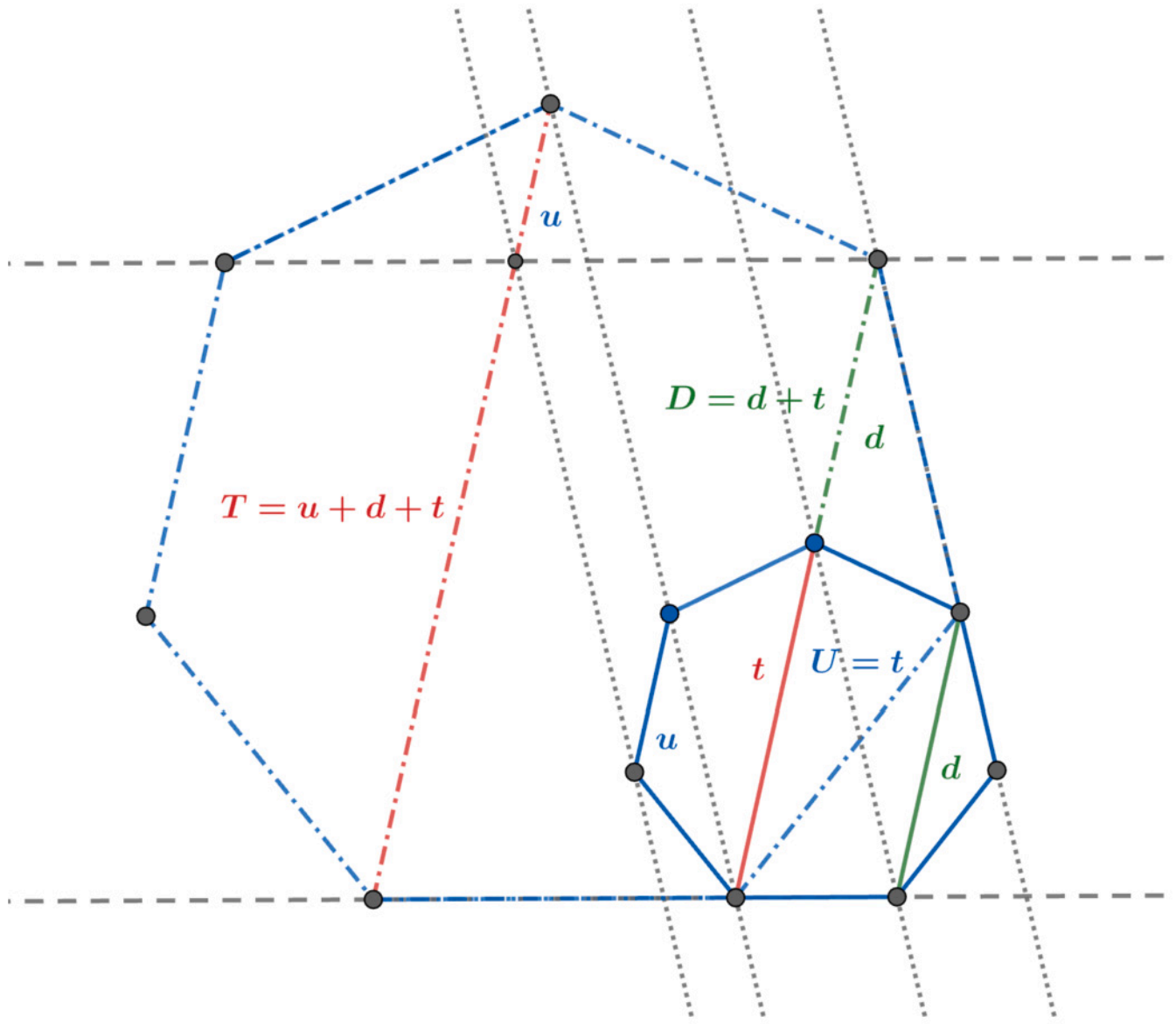

**Fig. 26.** *L'Archange*: The small heptagon extends to the larger one by similarity.

### Addenda

- The equation $(O)$ results from the sum of the previous relations and could also be written as $\frac{1}{u} = \frac{1}{d} + \frac{1}{t}$, which resembles a thin lens with focal length $u$ and distances to object $d$ and image $t$.
- The equations are identical under the permutations $ud \to -dt \to tu \to ud$ and $u^2 \to d^2 \to t^2 \to u^2$.

### 13.3. Towards a golden matrix

The four fundamental Ptolemy relations from equation (10) can well be written by the following matrix

$$\begin{bmatrix} u & d & t \end{bmatrix} \begin{bmatrix} u \\ d \\ t \end{bmatrix} = \begin{bmatrix} u^2 & t^2 - d^2 & d^2 - u^2 \\ ud & d^2 & t^2 - u^2 \\ ut & dt & t^2 \end{bmatrix}$$

We can construct a larger heptagon out of the first, with larger diagonal lengths, as can be seen in Fig. 26 or in equation (9). We identify one edge of the larger heptagon with one $t$-diagonal of the smaller heptagon, $U = t$. The prolongation of the other $t$-diagonal by a $d$-diagonal is identified with the second diagonal, $D = t + d$. Finally, the $T$-diagonal of the larger heptagon can be constructed by a translation, such that $T = u + d + t$.

$$U = t, \qquad D = d + t, \qquad T = u + d + t$$

or in matrix form:

$$\begin{bmatrix} U \\ D \\ T \end{bmatrix} = \begin{bmatrix} 0 & 0 & t \\ 0 & d & t \\ u & d & t \end{bmatrix} = \begin{bmatrix} 0 & 0 & 1 \\ 0 & 1 & 1 \\ 1 & 1 & 1 \end{bmatrix} \begin{bmatrix} u \\ d \\ t \end{bmatrix}$$

The absolute values in the larger heptagon grow, and the approximation of the ratios $\frac{T}{D}, \frac{-U}{T}, \frac{-D}{U}$ improves.

This allows an iterative computation of the polygon diagonals:

$$M = \begin{bmatrix} 0 & 0 & 1 \\ 0 & 1 & 1 \\ 1 & 1 & 1 \end{bmatrix}, \quad M^2 = \begin{bmatrix} 1 & 1 & 1 \\ 1 & 2 & 2 \\ 1 & 2 & 3 \end{bmatrix}, \quad \dots, \quad M^6 = \begin{bmatrix} 14 & 25 & 31 \\ 25 & 45 & 56 \\ 31 & 56 & 70 \end{bmatrix}, \quad \dots, \quad M^\infty = \begin{bmatrix} u^2 & t^2 - d^2 & d^2 - u^2 \\ ud & d^2 & t^2 - u^2 \\ ut & dt & t^2 \end{bmatrix}$$

Hence a parallelogram with the largest two-digit ratios $u : t : d = 31 : 70 : 56$ would give a good approximation for the vertices of a heptagon as its intersections.

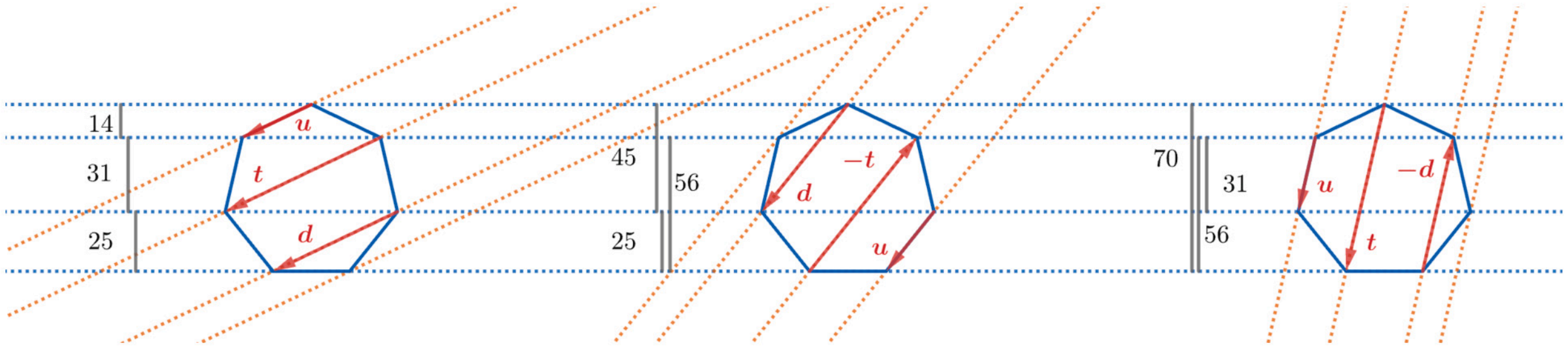

**Fig. 27.** Depending on the angle, the light rays induce a scale on the blue parallels: $(+u+t+d), (+d+u-t), (+t-d+u)$ (from left to right), exemplary ratios from $M^6$.

**Theorem 4.** *The powers $M^k$ of the grad n matrix*

$$M = \begin{bmatrix} 0 & \cdots & 0 & 1 \\ 0 & \cdots & 1 & 1 \\ & \iddots & \iddots & \vdots \\ 1 & 1 & \cdots & 1 \end{bmatrix}$$

*approximate the Golden Matrix $M^\infty$. Its values $(m^k_{nj})$ approximate the diagonals of a regular $(2n+1)$-gon:*

$$m^k_{n1} : m^k_{n2} : \cdots : m^k_{nn} \quad \approx \quad d_1 : d_2 : \cdots : d_n.$$

**Addenda**

- In his article *Au delà du nombre d'Or* (Beyond the Golden Number), de Faget de Casteljau (2001a) describes two flashes of inspiration, which he had at night in bed. The first, which he called *L'Ange* (Angel) made him add to the symmetry lines between opposite vertices the tangent line through its apex, cf. Fig. 25. This allows the vertex construction by a parallelogram. The second flash, only days later, delivered the recurrence relation, which he called *L'Archange* (Archangel), cf. Fig. 26.
- Gauss discovered, that a 17-gon (and some other $n$-gons with $n-1=2^{2^k}$) can be constructed by ruler and compass, which means, that $\cos\frac{2\pi}{17}$ can be calculated by use of addition, multiplication and square roots alone. De Casteljau adds to the exact approach of Gauss his approximative extension, which not only allows the calculation of $\cos\frac{2\pi}{n}$ for any impair $n$, but also delivers a recursive algorithm to approximate the vertices, and a simple construction of this approximation. He titled a one pager of his results *« Je me gausse de GAUSS »* ("I make fun of Gauss"), which is reproduced in Fig. 28.
- In the case of the pentagon, $n=5$, $M^k$ generates the Fibonacci series.
- As another application, de Faget de Casteljau (1990, 1995a) describes an approximation of a cycloid via a rolling regular polygon or the calculation of a brachistochrone.

## On algebra and number theory

This third part will address topics that at first glance might not seem related to geometry. We will look at quaternions, a nearly forgotten method of polynomial root finding, a generalisation of the Euclidean algorithm, the golden matrix with its link to regular polygons, all rounded off by an identity in number theory.

## 14. Quaternions

*When I wrote my book on quaternions, I wanted to scream from the beginning of the book "No, the quaternion is not a number, nor a matrix, nor a quadratic form, it is a system of orthonormal axes."* — de Casteljau[36]

After Leonhard Euler had introduced the imaginary number $i$ with $i^2=-1$ in 1748, Carl Friedrich Gauss, who we owe the term of *complex numbers*, came up in 1813 with the geometric interpretation of a plane with a real and an imaginary coordinate. It was William Rowan Hamilton, who in 1833 studied complex numbers as pairs of two such coordinates and their algebraic properties under addition and multiplication. He was furtheron in search of triplets with one real and two complex coordinates to generalise the Gaussian complex plane into space, but the *trinions* he searched for did not materialise. Instead on 16 October 1843, while walking along the Royal Canal in Dublin with his wife, he discovered that it needs three imaginary terms $i, j, k$, and he could't *"resist the*

[36] de Casteljau, personal communication (translated from French), 14 September 1995.

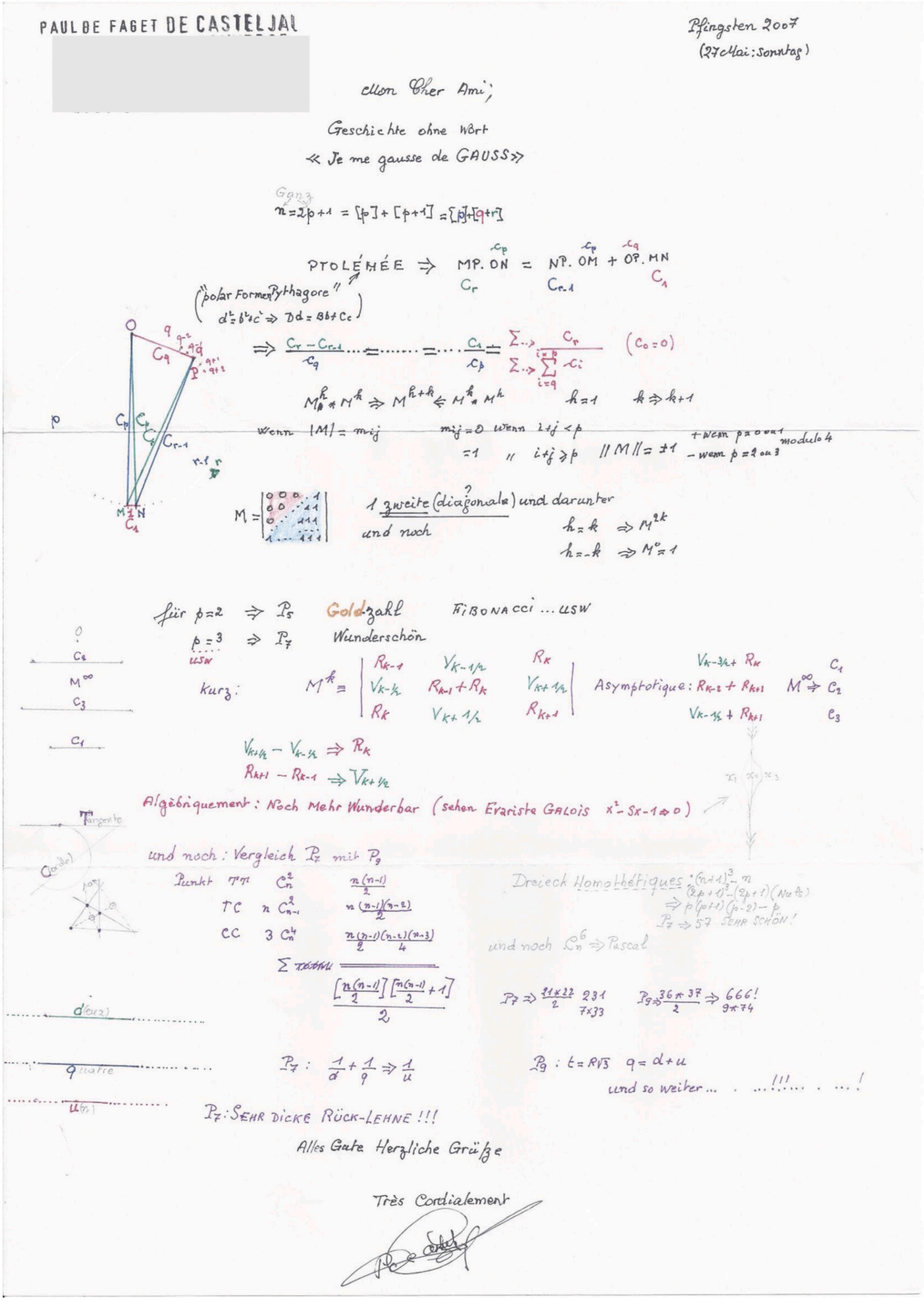
PAUL DE FAGET DE CASTELJAU

Pfingsten 2007
(27 Mai : Sonntag)

Mon Cher Ami;

Geschichte ohne Wort
« Je me gausse de GAUSS »

Ganz
$n = 2p+1 = [p] + [p+1] = [p] + [q+r]$

PTOLÉMÉE ⇒ MP. ON = NP. OM + OP. MN

("polar Formel Pythagore" $d^2 = b^2 + c^2 \Rightarrow Dd = Bb + Cc$)

$\Rightarrow \frac{C_r - C_{r-1}}{c_q} = \dots = \dots = \frac{C_1}{c_p} = \frac{\Sigma C_r}{\Sigma \sum_{i=q}^{i=p} c_i}$ $(C_0 = 0)$

$M^h_p * M^k \Rightarrow M^{h+k} \Leftarrow M^k * M^h$  $k = 1$  $k \Rightarrow k+1$

wenn $|M| = m_{ij}$  $m_{ij} = 0$ wenn $i+j < p$
$= 1$ " $i+j \geq p$  $\|M\| = \pm 1$  + wenn $p = 0$ ou 1, − wenn $p = 2$ ou 3 modulo 4

M = (matrix) 1 zweite (diagonale) und darunter
und noch  $h = k \Rightarrow M^{2k}$
$h = -k \Rightarrow M^0 = 1$

für $p = 2 \Rightarrow P_5$ Goldzahl FIBONACCI … usw
$p = 3 \Rightarrow P_7$ Wunderschön
usw

kurz: $M^k = \begin{vmatrix} R_{k-1} & V_{k-1/2} & R_k \\ V_{k-1/2} & R_{k-1}+R_k & V_{k+1/2} \\ R_k & V_{k+1/2} & R_{k+1} \end{vmatrix}$ Asymptotique: $V_{k-1/2} + R_k$, $R_{k-1} + R_{k+1}$, $V_{k-1/2} + R_{k+1}$ $M^\infty \Rightarrow C_1, C_2, C_3$

$V_{k+1/2} - V_{k-1/2} \Rightarrow R_k$
$R_{k+1} - R_{k-1} \Rightarrow V_{k+1/2}$

Algèbriquement: Noch Mehr Wunderbar (sehen Evariste GALOIS $x^3 - 5x - 1 = 0$)

und noch: Vergleich $P_7$ mit $P_9$

| | | | |
|---|---|---|---|
| Punkt | TT | $C_n^2$ | $\frac{n(n-1)}{2}$ |
| TC | n | $C_{n-1}^2$ | $n\frac{(n-1)(n-2)}{2}$ |
| CC | 3 | $C_n^4$ | $\frac{n(n-1)}{2}\frac{(n-2)(n-3)}{4}$ |

Σ TOTAL $\frac{\left[\frac{n(n-1)}{2}\right]\left[\frac{n(n-1)}{2}+1\right]}{2}$

Dreieck Homothétiques: $(n+1)^3 - n$, $(2p+1)^3 - (2p+1)(N+4)$ ⇒ $p(p+1)(p-2) - p$; $P_7 \Rightarrow 57$ SEHR SCHÖN!
und noch $C_n^6 \Rightarrow$ Pascal

$P_7 \Rightarrow \frac{21 \times 22}{2}$ 231 7×33  $P_9 \Rightarrow \frac{36 \times 37}{2} \Rightarrow$ 666! 9×74

$P_7: \frac{1}{d} + \frac{1}{q} \Rightarrow \frac{1}{u}$  $P_9: t = R\sqrt{3}$ $q = d + u$
und so weiter … !!! … !

$P_7$: SEHR DICKE RÜCK-LEHNE !!!

Alles Gute Herzliche Grüße

Très Cordialement

Fig. 28. «*Je me gausse de GAUSS* », original letter from de Casteljau to Boehm (personal data redacted), entitled *Geschichte ohne Wort[e]* (story with no words [needed]), Pentecost 2007; with kind permission of the de Casteljau family. It is one of his latest purely mathematical notes and demonstrates his preference for distilling an idea to its very essence.

*impulse - unphilosophical as it may have been - to cut with a knife on a stone of Brougham Bridge, as we passed it, the fundamental formula with the symbols, $i, j, k$"*(Hamilton, 1865), and he called those numbers composed of four terms *Quaternions*.

**Definition 1** *(Quaternion)*. A quaternion $q$ is defined as

$$q = t + ix + jy + kz$$

where $i, j, k$ are roots of $-1$, i.e. $i^2 = j^2 = k^2 = -1$.

We consider, that $k$ is result of the non-commutative product of $ij = k$ (from here $ik = iij = -j, kj = ijj = -i, jk = -ikk = i, ki = -kkj = j, ji = jjk = -k, ijk = k^2 = -1$). A norm can be introduced via a conjugate quaternion $\bar{q} = t - ix - jy - kz$, such that $n(q) = q\bar{q} = \bar{q}q = t^2 + x^2 + y^2 + z^2 = |q|^2$.

*14.1. De Casteljau's matrix representation*

Different to Hamilton, de Casteljau introduces a quaternion in its matrix form, with a main advantage in multiplications. Let us associate $q$ with a vector $\mathbf{q} = [t, x, y, z]^t$ to generate the matrix

$$Q = \begin{bmatrix} \mathbf{q} & \mathbf{q}i & \mathbf{q}j & \mathbf{q}k \end{bmatrix}$$

which conveniently includes the coefficients of $q$ and its (hyper)complex right multiples. The underlying imaginary structure associates, e.g., $qi = ti + xi^2 + yji + zki = -x + ti + zj - yk$ with $\mathbf{q}i = [-x, t, z, -y]^t$, which leads to the definition of the matrix quaternion $Q$ as

$$Q = \left[\begin{array}{c:c:c:c} t & -x & -y & -z \\ x & t & -z & y \\ y & z & t & -x \\ z & -y & x & t \end{array}\right]$$

De Casteljau defines a second matrix called anti-quaternion, which is composed of the transposed vectors and their (hyper)complex left multiples,

$$Q^{\star} = \begin{bmatrix} \mathbf{q}^t \\ i\mathbf{q}^t \\ j\mathbf{q}^t \\ k\mathbf{q}^t \end{bmatrix}$$

which associates, e.g., $iq = it + i^2x + ijy + ikz = -x + it - jz + ky$ with $i\mathbf{q}^t = [-x, t, -z, y]$.

**Definition 2** *(Anti-Quaternion)*. An anti-quaternion $Q^{\star}$ (related to the matrix quaternion $Q$) is defined as

$$Q^{\star} = \left[\begin{array}{cccc} t & x & y & z \\ \hdashline -x & t & -z & y \\ \hdashline -y & z & t & -x \\ \hdashline -z & -y & x & t \end{array}\right]$$

De Casteljau then turns a matrix quaternion $Q$ into its anti-quaternion $Q^{\star}$ by employing an auxiliary matrix $v$, which fulfills $Q^{\star} = vQv$ and $Q = vQ^{\star}v$

$$v = \begin{bmatrix} -1 & 0 & 0 & 0 \\ 0 & 1 & 0 & 0 \\ 0 & 0 & 1 & 0 \\ 0 & 0 & 0 & 1 \end{bmatrix}$$

We observe $v^2 = E$, $||v|| = -1$; left multiplication of $v$ inverts the first row, right multiplication of $v$ inverts the first column, which explains the above-mentioned transformation:

$$vQ = \left[\begin{array}{cccc} -t & x & y & z \\ \hdashline x & t & -z & y \\ y & z & t & -x \\ z & -y & x & t \end{array}\right] \qquad Qv = \left[\begin{array}{c:ccc} -t & -x & -y & -z \\ -x & t & -z & y \\ -y & z & t & -x \\ -z & -y & x & t \end{array}\right]$$

**Addenda**

- The transposed matrix quaternion $Q^t$ (identical to the conjugate matrix quaternion $\bar{Q}$) does not preserve the lower $3 \times 3$ matrix of $Q$, where it is different from the anti-quaternion $Q^{\star}$:

$$Q^t = \bar{Q} = \left[\begin{array}{c:ccc} t & x & y & z \\ \hdashline -x & t & z & -y \\ -y & -z & t & x \\ -z & y & -x & t \end{array}\right] \neq \left[\begin{array}{c:ccc} t & x & y & z \\ \hdashline -x & t & -z & y \\ -y & z & t & -x \\ -z & -y & x & t \end{array}\right] = Q^{☆}$$

- We can also translate the imaginary terms $i, j, k$ into matrix quaternions and anti-quaternions:

$$I = \begin{bmatrix} 0 & -1 & 0 & 0 \\ 1 & 0 & 0 & 0 \\ 0 & 0 & 0 & -1 \\ 0 & 0 & 1 & 0 \end{bmatrix} \quad J = \begin{bmatrix} 0 & 0 & -1 & 0 \\ 0 & 0 & 0 & 1 \\ 1 & 0 & 0 & 0 \\ 0 & -1 & 0 & 0 \end{bmatrix} \quad K = \begin{bmatrix} 0 & 0 & 0 & -1 \\ 0 & 0 & -1 & 0 \\ 0 & 1 & 0 & 0 \\ 1 & 0 & 0 & 1 \end{bmatrix}$$

$$I^{☆} = \begin{bmatrix} 0 & 1 & 0 & 0 \\ -1 & 0 & 0 & 0 \\ 0 & 0 & 0 & -1 \\ 0 & 0 & 1 & 0 \end{bmatrix} \quad J^{☆} = \begin{bmatrix} 0 & 0 & 1 & 0 \\ 0 & 0 & 0 & 1 \\ -1 & 0 & 0 & 0 \\ 0 & -1 & 0 & 0 \end{bmatrix} \quad K^{☆} = \begin{bmatrix} 0 & 0 & 0 & 1 \\ 0 & 0 & -1 & 0 \\ 0 & 1 & 0 & 0 \\ -1 & 0 & 0 & 1 \end{bmatrix}$$

  While it still holds that $IJ = K, JI = -K, IJK = -E$, we also have the anti-products $I^{☆}J^{☆} = K^{☆}, J^{☆}I^{☆} = -K^{☆}, I^{☆}I^{☆} = -E$ with the identity matrix $E$.
- De Casteljau offers another approach to quaternions via permutations and magic squares. Let us assume two sets of items $(d, a, b, c)$ and $(t, x, y, z)$ are permuted to the 16 pairs

$$\begin{array}{cccc} dt & ax & by & cz \\ bz & cy & dx & at \\ cx & bt & az & dy \\ ay & dz & ct & bx \end{array}$$

  which can be clustered in six different ways to get four families with all eight items each:

$$lines: \begin{array}{cccc} 0 & 0 & 0 & 0 \\ 1 & 1 & 1 & 1 \\ 2 & 2 & 2 & 2 \\ 3 & 3 & 3 & 3 \end{array} \quad columns: \begin{array}{cccc} 0 & 1 & 2 & 3 \\ 0 & 1 & 2 & 3 \\ 0 & 1 & 2 & 3 \\ 0 & 1 & 2 & 3 \end{array} \quad diagonals: \begin{array}{cccc} 0 & 1 & 2 & 3 \\ 1 & 0 & 3 & 2 \\ 2 & 3 & 0 & 1 \\ 3 & 2 & 1 & 0 \end{array}$$

$$little\ squares: \begin{array}{cccc} 0 & 0 & 2 & 2 \\ 0 & 0 & 2 & 2 \\ 3 & 3 & 1 & 1 \\ 3 & 3 & 1 & 1 \end{array} \quad mid-squares: \begin{array}{cccc} 0 & 2 & 0 & 2 \\ 3 & 1 & 3 & 1 \\ 0 & 2 & 0 & 2 \\ 3 & 1 & 3 & 1 \end{array} \quad others: \begin{array}{cccc} 0 & 2 & 2 & 0 \\ 3 & 1 & 1 & 3 \\ 3 & 1 & 1 & 3 \\ 0 & 2 & 2 & 0 \end{array}$$

  If we attach for example, the values $1, 4, 3, 2$ to $(d, a, b, c)$ and the multiples of four $4, 0, 12, 8$ to $(t, x, y, z)$, then adding two different family members leads to the following magic square with sums of 34:

$$\begin{array}{cccc} 5 & 4 & 14 & 11 \\ 15 & 10 & 8 & 1 \\ 12 & 13 & 3 & 6 \\ 13 & 7 & 9 & 16 \end{array}$$

  We observe that the defined quaternion (and anti-quaternion) follows the diagonal pattern (omitting signs), which to de Casteljau is only one choice out of six. Anyone interested in delving deeper into this subject is encouraged to follow his thoughts in *Les Quaternions* (1987, p.34–37).

### 14.2. Quaternion products

Let us follow de Casteljau and look at the product of two matrix quaternions $Q_1$ and $Q_2$, which can take three different forms: quaternion × quaternion $Q_1Q_2$, anti-quaternion × anti-quaternion $Q_1^{☆}Q_2^{☆}$, and quaternion × anti-quaternion $Q_1^{☆}Q_2$ or $Q_1Q_2^{☆}$, where

$$Q_1^{☆} = \begin{bmatrix} d & a & b & c \\ -a & d & -c & b \\ -b & c & d & -a \\ -c & -b & a & d \end{bmatrix} \qquad Q_2 = \begin{bmatrix} t & -x & -y & -z \\ x & t & -z & y \\ y & z & t & -x \\ z & -y & x & t \end{bmatrix}$$

with norms

$$n_1 = |Q_1|^2 = |Q_1^{☆}|^2 = d^2 + a^2 + b^2 + c^2 \qquad \text{and} \qquad n_2 = |Q_2|^2 = |Q_2^{☆}|^2 = t^2 + x^2 + y^2 + z^2$$

The product of two anti-quaternions is an anti-quaternion:

$$Q_1^{\star}Q_2^{\star} = \begin{bmatrix} K & L & M & N \\ -L & K & -N & M \\ -M & N & K & -L \\ -N & -M & L & K \end{bmatrix}$$

with

$$\begin{array}{lllll} K &=& dt-ax-by-cz &=& 2dt-k-l-m-n \\ L &=& at+dx-cy+bz &=& -2cy+k+l-m-n \\ M &=& bt+cx+dy-az &=& -2az+k-l+m-n \\ N &=& ct-bx+ay+dz &=& -2bx+k-l-m+n \end{array} \quad \text{and } k,l,m,n \text{ as} \quad \begin{array}{lll} 4k &=& (d+a+b+c)(t+x+y+z) \\ 4l &=& (d+a-b-c)(t+x-y-z) \\ 4m &=& (d-a+b-c)(t-x+y-z) \\ 4n &=& (d-a-b+c)(t-x-y+z) \end{array}$$

Making use of the auxiliary matrix $v$ from above, we observe that similarly the product of two quaternions is a quaternion:

$$Q_1Q_2 = vQ_1^{\star}vvAQ_2^{\star}v = vQ_1^{\star}Q_2^{\star}v = \begin{bmatrix} K & -L & -M & -N \\ L & K & -N & M \\ M & N & K & -L \\ N & -M & L & K \end{bmatrix}$$

But what is the product of a matrix quaternion and an anti-quaternion?

**Theorem 5** *(Quaternion products are $E^4$ integers). The commutative product of an integer matrix quaternion and anti-quaternion results in a positive integer matrix in $E^4$.*

Let us look at the product $P = Q_1^{\star}Q_2$ of a matrix quaternion and an arbitrary anti-quaternion:

$$Q_1^{\star}Q_2 = \begin{bmatrix} dt+ax+by+cz & at-dx-cy+bz & bt+cx-dy-az & ct-bx+ay-dz \\ -at+dx-cy+bz & dt+ax-by-cz & -ct+bx+ay-dz & bt+cx+dy+az \\ -bt+cx+dy-az & ct+bx+ay+dz & dt-ax+by-cz & -at-dx+cy+bz \\ -ct-bx+ay+dz & -bt+cx-dy+az & at+dx+cy+bz & dt-ax-by+cz \end{bmatrix} = Q_2Q_1^{\star}$$

It is then

$$PP^t = n_1n_2E = (d^2+a^2+b^2+c^2)(t^2+x^2+y^2+z^2)E = P^tP$$

and $P = Q_1^{\star}Q_2$ is positive.

In particular, the product of a quaternion $Q$ with its proper anti-quaternion $Q^{\star}$ becomes

$$Q^{\star}Q = \begin{bmatrix} d^2+a^2+b^2+c^2 & ad-da-cb+bc & bd+ca-db-ac & cd-ba+ab-dc \\ -ad+da-cb+bc & d^2+a^2-b^2-c^2 & -cd+ba+ab-dc & bd+ca+db+ac \\ -bd+ca+db-ac & cd+ba+ab+dc & d^2-a^2+b^2-c^2 & -ad-da+cb+bc \\ -cd-ba+ab+dc & -bd+ca-db+ac & ad+da+cb+bc & d^2-a^2-b^2+c^2 \end{bmatrix}$$

$$= \begin{bmatrix} d^2+a^2+b^2+c^2 & 0 & 0 & 0 \\ 0 & d^2+a^2-b^2-c^2 & 2(ab-cd) & 2(ac+bd) \\ 0 & 2(ab+cd) & d^2-a^2+b^2-c^2 & 2(bc-ad) \\ 0 & 2(ac-bd) & 2(bc+ad) & d^2-a^2-b^2+c^2 \end{bmatrix} \quad (12)$$

which proves the theorem.

**Theorem 6** *($E^4$ integers are quaternion products). Every integer matrix in $E^4$ can be decomposed into a product of an integer matrix quaternion and an anti-quaternion.*

In *Les Quaternions,* de Casteljau (1987, p.65-68) offers three proofs, which the reader can follow through in further detail. To indicate one of them, let us assume an arbitrary $4 \times 4$ matrix in $E^4$

$$A = \begin{bmatrix} a_{00} & a_{01} & a_{02} & a_{03} \\ a_{10} & a_{11} & a_{12} & a_{13} \\ a_{20} & a_{21} & a_{22} & a_{23} \\ a_{30} & a_{31} & a_{32} & a_{33} \end{bmatrix}$$

and the following matrix $T$, which represents what de Casteljau called a *tetragonal transformation*:

$$T(A) := \frac{1}{2}\begin{bmatrix} a_{00}+a_{11}+a_{22}+a_{33} & a_{01}-a_{10}-a_{23}+a_{32} & a_{02}+a_{13}-a_{20}-a_{31} & a_{03}-a_{12}+a_{21}-a_{30} \\ -a_{01}+a_{10}-a_{23}+a_{32} & a_{00}+a_{11}-a_{22}-a_{33} & -a_{03}+a_{12}+a_{21}-a_{30} & a_{02}+a_{13}+a_{20}+a_{31} \\ -a_{02}+a_{13}+a_{20}-a_{31} & a_{03}+a_{12}+a_{21}+a_{30} & a_{00}-a_{11}+a_{22}-a_{33} & -a_{01}-a_{10}+a_{23}+a_{32} \\ -a_{03}-a_{12}+a_{21}+a_{30} & -a_{02}+a_{13}-a_{20}+a_{31} & a_{01}+a_{10}+a_{23}+a_{32} & a_{00}-a_{11}-a_{22}+a_{33} \end{bmatrix} \quad (13)$$

It is

$$T^2(A) = T(T(A)) = \frac{1}{2}\begin{bmatrix} 4a_{00} & 4a_{01} & 4a_{02} & 4a_{03} \\ 4a_{10} & 4a_{11} & 4a_{12} & 4a_{12} \\ 4a_{20} & 4a_{21} & 4a_{22} & 4a_{22} \\ 4a_{30} & 4a_{31} & 4a_{32} & 4a_{32} \end{bmatrix} = 2A$$

We could also identify the $a_{ij}$ with the products $\mathbf{d}_i^t\mathbf{t}_j$ such that

$$(a_{ij}) = (\mathbf{d}_i^t\mathbf{t}_j) = Q_1^{\star} Q_2 =$$

$$= \begin{bmatrix} dt+ax+by+cz & at-dx-cy+bz & bt+cx-dy-az & ct-bx+ay-dz \\ -at+dx-cy+bz & dt+ax-by-cz & -ct+bx+ay-dz & bt+cx+dy+az \\ -bt+cx+dy-az & ct+bx+ay+dz & dt-ax+by-cz & -at-dx+cy+bz \\ -ct-bx+ay+dz & -bt+cx-dy+az & at+dx+cy+bz & dt-ax-by+cz \end{bmatrix}$$

$$T(Q_1^{\star} Q_2) = 2\begin{bmatrix} dt & at & bt & ct \\ dx & ax & bx & cx \\ dy & ay & by & cy \\ dz & az & bz & cz \end{bmatrix}$$

**Addenda**

- Equation (12) was already discovered by Benjamin Olinde Rodrigues (1840) without the use of (hyper)complex numbers. A more thorough and worthwile exploration of the approaches of both Hamilton and Rodrigues was conducted by Simon Altmann (1989).
- If $b = c = 0$, we find Plato's formula on Pythagorean triples

  $$R = \frac{1}{d^2+a^2}\begin{bmatrix} d^2+a^2 & 0 & 0 & 0 \\ 0 & d^2+a^2 & 0 & 0 \\ 0 & 0 & d^2-a^2 & -2ad \\ 0 & 0 & 2ad & d^2-a^2 \end{bmatrix}$$

  Thus

  $$(d^2-a^2)^2 - (2ad)(-2ad) = (d^2+a^2)^2$$

- De Casteljau also sees Lucas (1876) formula on $3 \times 3$ magic squares, who considers the following table of nine quantities

  $$\begin{matrix} d^2+a^2-b^2-c^2 & 2(ab-cd) & 2(ac+bd) \\ 2(ab+cd) & d^2-a^2+b^2-c^2 & 2(bc-ad) \\ 2(ac-bd) & 2(ad+bc) & d^2-a^2-b^2+c^2 \end{matrix}$$

  and concludes, that the sum of squares of one row or one column equal the constant term $(d^2+a^2+b^2+c^2)^2$.
  If the sum of all nine terms shall sum up to $d^2+a^2+b^2+c^2$, we can deduce

  $$d^2+a^2+b^2+c^2 = \pm(3d^2-a^2-b^2-c^2+4ab+4bc+4ca)$$

  The negative sign leads us to

  $$d^2+ab+bc+ca = 0$$

  which could be understood as a function of $b, c, d$ to determine $a = -(d^2+bc)/(b+c)$

### *14.3. Rotations*

De Casteljau leads us now to an often-used application of quaternions, namely rotations. Normalising the quaternion product results in

$$R = \frac{1}{n} Q^{\star} Q = \frac{1}{n} Q Q^{\star} = \frac{1}{d^2+a^2+b^2+c^2} Q^{\star} Q$$

**Theorem 7** *(rotations). The product of a matrix quaternion $Q$ and its anti-quaternion $Q^{\star}$ is a rotation $R$*

$$R = \frac{1}{n} Q Q^{\star}.$$

The application of $R$ on a vector $\mathbf{a}^t = [0, a, b, c]$ delivers in a first step

$$Q^{\star}\mathbf{a} = \begin{bmatrix} a^2+b^2+c^2 \\ da \\ db \\ dc \end{bmatrix}$$

and finally

$$\frac{Q}{n}Q^{\star}\mathbf{a} = \frac{1}{a^2+b^2+c^2}\begin{bmatrix} 0 \\ a(a^2+b^2+c^2+d^2) \\ b(b^2+a^2+c^2+d^2) \\ c(c^2+a^2+b^2+d^2) \end{bmatrix} = \begin{bmatrix} 0 \\ a \\ b \\ c \end{bmatrix} = \mathbf{a}$$

which shows that $\mathbf{a}$ is unchanged under the operation, thus it is the rotation axis.

The rotation angle can be found by applying the rotation on a unit vector $\mathbf{b} = [0, \alpha, \beta, \gamma]^t$, which is perpendicular to axis $\mathbf{a}$, meaning

$$\mathbf{a}^t\mathbf{b} = a\alpha + b\beta + c\gamma = 0 \qquad \mathbf{b}^t\mathbf{b} = \alpha^2+\beta^2+\gamma^2 = 1$$

We now apply equation (12) and the rotation $R\mathbf{b}$ results in

$$R\mathbf{b} = \frac{1}{n}Q^{\star}Q\mathbf{b} = \frac{1}{n}\begin{bmatrix} 0 \\ (d^2+a^2-b^2-c^2)\alpha + 2(ab-cd)\beta + 2(ac+bd)\gamma \\ 2(ab+cd)\alpha + (d^2-a^2+b^2-c^2)\beta + 2(bc-ad)\gamma \\ 2(ac-bd)\alpha + 2(ad+bc)\beta + (d^2-a^2-b^2+c^2)\gamma \end{bmatrix}$$

$$= \frac{d^2-a^2-b^2-c^2}{d^2+a^2+b^2+c^2}\begin{bmatrix} 0 \\ \alpha \\ \beta \\ \gamma \end{bmatrix} + \frac{2d}{d^2+a^2+b^2+c^2}\begin{bmatrix} 0 \\ b\gamma - c\beta \\ c\alpha - a\gamma \\ a\beta - b\alpha \end{bmatrix} = \mathbf{b}\cos\varphi + \mathbf{c}\sin\varphi$$

with $\cos\varphi = (d^2-a^2-b^2-c^2)/(d^2+a^2+b^2+c^2)$ and $\mathbf{c} = \mathbf{a}\wedge\mathbf{b}$. From here, we get

$$\sin^2\varphi = 1 - \cos^2\varphi = \frac{(d^2+(a^2+b^2+c^2))^2 - (d^2-(a^2+b^2+c^2))^2}{(d^2+a^2+b^2+c^2)^2} = \frac{4d^2(a^2+b^2+c^2)}{(d^2+a^2+b^2+c^2)^2}$$

thus $\sin\varphi = 2d/(d^2+a^2+b^2+c^2)$ with unit vector $\mathbf{a}$ and thus $\mathbf{a}^t\mathbf{a} = a^2+b^2+c^2 = 1$

**Addenda**

- Compared with Euler angles $\varphi, \vartheta, \omega$, a rotation by quaternions requires fewer operations. Euler's rotation matrix has 13 terms:

$$\begin{bmatrix} \cos\omega & -\sin\omega & 0 \\ \sin\omega & \cos\omega & 0 \\ 0 & 0 & 1 \end{bmatrix} \cdot \begin{bmatrix} 1 & 0 & 0 \\ 0 & \cos\vartheta & -\sin\vartheta \\ 0 & \sin\vartheta & \cos\vartheta \end{bmatrix} \cdot \begin{bmatrix} \cos\varphi & -\sin\varphi & 0 \\ \sin\varphi & \cos\varphi & 0 \\ 0 & 0 & 1 \end{bmatrix}$$

$$= \begin{bmatrix} \cos\varphi\cos\omega - \sin\varphi\sin\omega\cos\vartheta & -\sin\varphi\cos\omega - \cos\varphi\sin\omega\cos\vartheta & \sin\omega\sin\vartheta \\ \cos\varphi\sin\omega + \sin\varphi\cos\omega\cos\vartheta & -\sin\varphi\sin\omega + \cos\varphi\cos\omega\cos\vartheta & -\cos\omega\sin\vartheta \\ \sin\varphi\sin\vartheta & \cos\varphi\sin\vartheta & \cos\vartheta \end{bmatrix}$$

  The product of three matrix quaternions instead calculates only four terms:

$$\begin{bmatrix} \cos\frac{\omega}{2} \\ 0 \\ 0 \\ k\sin\frac{\omega}{2} \end{bmatrix} \cdot \begin{bmatrix} \cos\frac{\vartheta}{2} \\ i\sin\frac{\vartheta}{2} \\ 0 \\ 0 \end{bmatrix} \cdot \begin{bmatrix} \cos\frac{\varphi}{2} \\ 0 \\ 0 \\ k\sin\frac{\varphi}{2} \end{bmatrix} = \begin{bmatrix} \cos\frac{\omega+\varphi}{2}\cos\frac{\vartheta}{2} \\ i\cos\frac{\omega+\varphi}{2}\sin\frac{\vartheta}{2} \\ j\sin\frac{\varphi}{2} \\ k\sin\frac{\omega+\varphi}{2}\cos\frac{\vartheta}{2} \end{bmatrix}$$

  with rotation angle $\cos u/2 = \cos\frac{\omega+\varphi}{2}\cos\frac{\vartheta}{2}$ and rotation axis $\left[\cos\frac{\omega-\varphi}{2}\sin\frac{\vartheta}{2},\quad \sin\frac{\omega-\varphi}{2}\sin\frac{\vartheta}{2},\quad \sin\frac{\omega+\varphi}{2}\cos\frac{\vartheta}{2}\right]^t$.
- Here we see the fundamental difference between Hamilton and de Casteljau. While Hamilton wrote a rotation as $X_1 = QX\bar{Q}$, de Casteljau denotes it as $X_1 = (QQ^{\star})X = (Q^{\star}Q)X$, which allows the performance of an abstraction of $X$, giving

$$R = Q^{\star}Q = QQ^{\star}$$

### 14.4. Lorentz transformation

The study of functions of quaternions leads de Casteljau to study derivatives of quaternions. He applies the above-mentioned tetragonal transformation to the Jacobian which results in a form that appears as the electromagnetic Lorentz transformation. (de Casteljau, 1987, p.102) We will reproduce his main argument.

After introducing the matrix operator

$$\nabla = \begin{bmatrix} \frac{\partial}{\partial t} & -\frac{\partial}{\partial x} & -\frac{\partial}{\partial y} & -\frac{\partial}{\partial z} \\ \frac{\partial}{\partial x} & \frac{\partial}{\partial t} & -\frac{\partial}{\partial z} & \frac{\partial}{\partial y} \\ \frac{\partial}{\partial y} & \frac{\partial}{\partial z} & \frac{\partial}{\partial t} & -\frac{\partial}{\partial x} \\ \frac{\partial}{\partial z} & -\frac{\partial}{\partial y} & \frac{\partial}{\partial x} & \frac{\partial}{\partial t} \end{bmatrix}$$

he applies it to four functions $U(\mathbf{t}), V_x(\mathbf{t}), V_y(\mathbf{t}), V_z(\mathbf{t})$ in $t, x, y, z$, using the tetragonal transformation (omitting the factor 1/2):

$$\nabla^{\star}\mathbf{U} = \begin{bmatrix} \frac{\partial U}{\partial t} + \frac{\partial V_x}{\partial x} + \frac{\partial V_y}{\partial y} + \frac{\partial V_z}{\partial z} & -\frac{\partial V_x}{\partial t} + \frac{\partial U}{\partial x} + \frac{\partial V_z}{\partial y} - \frac{\partial V_y}{\partial z} & -\frac{\partial V_y}{\partial t} - \frac{\partial V_z}{\partial x} + \frac{\partial U}{\partial y} + \frac{\partial V_x}{\partial z} & -\frac{\partial V_z}{\partial t} + \frac{\partial V_y}{\partial x} - \frac{\partial V_x}{\partial y} + \frac{\partial U}{\partial z} \\ -\frac{\partial U}{\partial x} + \frac{\partial V_x}{\partial t} - \frac{\partial V_y}{\partial z} + \frac{\partial V_z}{\partial y} & \frac{\partial V_x}{\partial x} + \frac{\partial U}{\partial t} - \frac{\partial V_z}{\partial z} - \frac{\partial V_y}{\partial y} & \frac{\partial V_y}{\partial x} - \frac{\partial V_z}{\partial t} - \frac{\partial U}{\partial z} + \frac{\partial V_x}{\partial y} & \frac{\partial V_z}{\partial x} + \frac{\partial V_y}{\partial t} + \frac{\partial V_x}{\partial z} + \frac{\partial U}{\partial y} \\ -\frac{\partial U}{\partial y} + \frac{\partial V_x}{\partial z} + \frac{\partial V_y}{\partial t} - \frac{\partial V_z}{\partial x} & \frac{\partial V_x}{\partial y} + \frac{\partial U}{\partial z} + \frac{\partial V_z}{\partial t} + \frac{\partial V_y}{\partial x} & \frac{\partial V_y}{\partial y} - \frac{\partial V_z}{\partial z} + \frac{\partial U}{\partial t} - \frac{\partial V_x}{\partial x} & \frac{\partial V_z}{\partial y} + \frac{\partial V_y}{\partial z} - \frac{\partial V_x}{\partial t} - \frac{\partial U}{\partial x} \\ -\frac{\partial U}{\partial z} - \frac{\partial V_x}{\partial y} + \frac{\partial V_y}{\partial x} + \frac{\partial V_z}{\partial t} & \frac{\partial V_x}{\partial z} - \frac{\partial U}{\partial y} + \frac{\partial V_z}{\partial x} - \frac{\partial V_y}{\partial t} & \frac{\partial V_y}{\partial z} + \frac{\partial V_z}{\partial y} + \frac{\partial U}{\partial x} + \frac{\partial V_x}{\partial t} & \frac{\partial V_z}{\partial z} - \frac{\partial V_y}{\partial y} - \frac{\partial V_x}{\partial x} + \frac{\partial U}{\partial t} \end{bmatrix}$$

**Theorem 8.** *The functions $U, V_x, V_y, V_z$ represent a conformal map on $E_4$, iff $\nabla^{\star}\mathbf{U}$ is a rank 1 matrix.*

In this case, he sees nine proportionality equations satisfied, row zero and column zero result in

$$\left[\frac{\partial U}{\partial t} + \mathrm{div}\mathbf{V};\ \varepsilon\left(\frac{\partial \mathbf{V}}{\partial t} - \mathrm{grad}U\right) + \mathrm{rot}\mathbf{V}\right]$$

with $\varepsilon = +1$ for row zero and $\varepsilon = -1$ for column zero. The connection to Lorentz' electromagnetic equations comes by the identification of grad$U$ as electrostatic field $\mathbf{E}$ with induction potential $\frac{\partial \mathbf{V}}{\partial t}$, and rot$\mathbf{V}$ as magnetic field $\mathbf{H}$. The equations then read

$$\mathbf{E} = -\mathrm{grad}U - \frac{\partial \mathbf{V}}{\partial t} \qquad \mathbf{H} = \mathrm{rot}\mathbf{V} \tag{14}$$

based on the assumption that $U$ is seen as a potential and $\mathbf{V}$ as a vector of potentials, such that the delayed potentials are

$$U = \int \rho_{t-r/c}\frac{\mathrm{d}\tau}{\tau} \qquad \mathbf{V} = \int \mathbf{j}_{t-r/c}\frac{\mathrm{d}\tau}{\tau}$$

with distance $r$, local electric charge $\rho$, speed of light $c$, time $t$, current density $\mathbf{j}$, volume element $\mathrm{d}\tau$.

De Casteljau found this result only by using his approach of matrix anti-quaternions, which allows to see all quaternion coefficients at a glance.

An in-depth engagement with *Les Quaternions* is highly recommended for readers interested in algebra, particularly to explore further aspects related to metric geometry.

### Addendum

The previously mentioned auxiliary matrix $\nu$ is directly connected to Minkowski's relativistic $\mathrm{d}s^2$. Assuming $\mathrm{d}\mathbf{u} = \begin{bmatrix} c\mathrm{d}t & \mathrm{d}x & \mathrm{d}y & \mathrm{d}z \end{bmatrix}^t$, we find

$$\mathrm{d}\mathbf{u}^t \nu\, \mathrm{d}\mathbf{u} = -c^2\mathrm{d}t^2 + \mathrm{d}x^2 + \mathrm{d}y^2 + \mathrm{d}z^2$$

However, de Casteljau did not endorse the concept of space reversal with its mirror symmetry, which is represented by the negative $-\mathrm{d}s^2$. He concludes (de Casteljau, 1987, p.103):

*But if we retain the concepts of quaternion and anti-quaternion, there is no need for the aberrant matrix $\nu$ to transition from one to the other. Thus, the notion of a negative square or of imaginaries is evacuated. There remains only the difference between the principal diagonal representing the identity and the three asymmetric diagonals of the matrix quaternion to simulate the fourth roots of unity.*

He instead regarded quaternions as a system of orthonormal axes in $E_4$, a concept akin to the grids discussed in section 10 for $E_2$.

## 15. A (nearly) forgotten way to find roots

> *Vincent's theorem of 1836, which was only recently discovered by the author of this article, is of extreme importance because it constitutes the basis of the fastest method existing for the isolation of the real roots of a polynomial equation (using exact integer arithmetic). —* Akritas (1981)

To find the root of a polynomial, de Casteljau applies a nearly forgotten algorithm by French mathematician Alexandre Vincent (1797-1868), which he might have seen in the readings of his childhood. According to Alkis Akritas, the algorithm was still used by Serret in 1887, but was later forgotten due to the dominance of Sturm's approach. It is today used in an improved form by almost all computer algebra systems, known as Vincent–Akritas–Strzeboński (VAS) method (Akritas et al., 2008).

### *15.1. Vincent's algorithm*

Let us first revisit the original approach by Vincent (1834), which was reprinted 1836 in *Journal de Mathématiques Pures et Appliquées* with its broader audience.

**Theorem 9** *(Vincent 1834). If in a rational polynomial $f(x) = A + Bx + Cx^2 + Dx^3 + \ldots + Nx^n$ without multiple roots we perform successive parameter transformations of the form*

$$x = \alpha_1 + \frac{1}{x_1}, \quad x_1 = \alpha_2 + \frac{1}{x_2}, \quad x_2 = \alpha_3 + \frac{1}{x_3},$$

*with positive integers $\alpha_i$, the resulting polynomial has either no sign variations, or one. In the second case, $P$ has exactly one positive root:*

$$x = \alpha_1 + \cfrac{1}{\alpha_2 + \cfrac{1}{\alpha_3 + \cfrac{1}{\ddots}}}$$

*In the first case, $P$ does not have a root.*

Thus, we can isolate the root of a rational polynomial $f(x) = A + Bx + Cx^2 + Dx^3 + \cdots$ by computing step-by-step $f(x), f(x-1), f(x-2), \ldots,$ as long as the resulting transformations provoke sign variations. If this is the case after $i$ iterations, then the root is between $i$ and $i+1$ and the transform reads

$$f_i(x_i) = A_i x_i^m + B_i x_i^{m-1} + C_i x_i^{m-2} + D_i x_i^{m-3} + \ldots + N_i x_i^0$$

with

$$\begin{array}{rcccccccccc} A_i & = & A_{i-1} & + & \binom{1}{0} B_{i-1} & + & \binom{2}{0} C_{i-1} & + & \cdots & + & N_{i-1} \\ B_i & = & & & \binom{1}{1} B_{i-1} & + & \binom{2}{1} C_{i-1} & + & \cdots & + & \binom{m}{1} N_{i-1} \\ C_i & = & & & & & \binom{2}{2} C_{i-1} & + & \cdots & + & \binom{m}{2} N_{i-1} \\ \vdots & & & & & & & & & & \\ N_i & = & & & & & & & & & \binom{m}{m} N_{i-1} \end{array}$$

**Example**

An example shows the power of the related algorithm. We want to approximate the roots of $1 - 2x - x^2 + x^3 = 0$ and apply Vincent's algorithm until there are no further sign variations:

| $i$ | $A$ | $B$ | $C$ | $D$ |
|---|---|---|---|---|
| 0 | 1 | −2 | −1 | 1 |
| 1 | **−1** | **−1** | **2** | **1** |
| 2 | 1 | 6 | 5 | 1 |

There is a root between 1 and 2, thus $x = 1 + \frac{1}{x_1}$, and the transformed polynomial reads

$$-x_1^3 - x_1^2 + 2x_1 + 1 = 0 \quad \text{with the continued fraction} \quad x = 1 + \frac{1}{x_1} = \frac{x_1 + 1}{x_1}$$

To isolate the roots, we again apply Vincent's algorithm:

| $i$ | $A$ | $B$ | $C$ | $D$ |
|---|---|---|---|---|
| 0 | 1 | 2 | −1 | −1 |
| 1 | **1** | **−3** | **−4** | **−1** |
| 2 | −7 | −14 | −7 | −1 |

One root is between 1 and 2, thus $x_1 = 1 + \frac{1}{x_2}$, and the transformed polynomial reads

$$x_2^3 - 3x_2^2 - 4x_2 - 1 = 0 \quad \text{with the continued fraction} \quad x = 1 + \frac{1}{1+\frac{1}{x_2}} = 1 + \frac{x_2}{x_2+1} = \frac{2x_2+1}{x_2+1}$$

Let us apply the algorithm again to isolate the roots $x_2$:

| $i$ | $A$ | $B$ | $C$ | $D$ |
|---|---|---|---|---|
| 0 | −1 | 4 | −3 | 1 |
| 1 | −7 | −7 | 0 | 1 |
| 2 | −13 | −4 | 3 | 1 |
| 3 | −13 | 5 | 6 | 1 |
| 4 | **−1** | **20** | **9** | **1** |
| 5 | 29 | 41 | 12 | 1 |

There is a root between 4 and 5, thus $x_2 = 4 + \frac{1}{x_3}$, and the transformed polynomial reads

$$-x_3^3 + 20x_3^2 + 9x_3 + 1 = 0 \quad \text{with the continued fraction} \quad x = 1 + \frac{1}{1+\frac{1}{4+\frac{1}{x_3}}} = \frac{9x_3+2}{5x_3+1}$$

**Another example**
In a second example, we apply Vincent's algorithm on $1 - 5x + 6x^2 - x^3 = 0$ to isolate its roots:

| $i$ | $A$ | $B$ | $C$ | $D$ |
|---|---|---|---|---|
| 0 | 1 | −5 | 6 | −1 |
| 1 | 1 | 4 | 3 | −1 |
| 2 | 7 | 7 | 0 | 1 |
| 3 | 13 | 4 | −3 | −1 |
| 4 | 13 | −5 | −6 | −1 |
| 5 | **1** | **−20** | **−9** | **−1** |
| 6 | −29 | −41 | −12 | −1 |

which delivers a root between 5 and 6, thus $x = 5 + \frac{1}{t}$, where $t$ is a root of

$$t^3 - 20t^2 - 9t - 1 = 0 \quad \text{and can be translated into} \quad t = \frac{9t+1}{t(t-20)}$$

### *15.2. Backward substitution*

De Casteljau mentions Vincent's method in *Le Lissage*, and he develops the following substitution table,[37] where we also find convergents of the related continued fractions as $\frac{P+Q}{Q}$

| equation | substitution | $P$ | $Q$ | $P+Q$ |
|---|---|---|---|---|
| $x_{-1}^3 - 41x_{-1}^2 + 558x_{-1} - 2521 = 0$ | $x_{-1} = 13 + 1/x$ | $-x_{-1} + 14$ | $x_{-1} - 13$ | 1 |
| $x^3 - x^2 - 2x + 1 = 0$ | $x = 1 + 1/x_1$ | $x - 1$ | 1 | $x$ |
| $x_1^3 + x_1^2 - 2x_1 - 1 = 0$ | $x_1 = 1 + 1/x_2$ | 1 | $x_1$ | $x_1 + 1$ |
| $x_2^3 - 3x_2^2 - 4x_2 - 1 = 0$ | $x_2 = 4 + 1/x_3$ | $x_2$ | $x_2 + 1$ | $2x_2 + 1$ |
| $x_3^3 - 20x_3^2 - 9x_3 - 1 = 0$ | $x_3 = 20 + 1/x_4$ | $4x_3$ | $5x_3 + 1$ | $9x_3 + 2$ |
| $181x_4^3 - 391x_4^2 - 40x_4 - 1 = 0$ | $\cdots$ | $81x_4 + 4$ | $101x_4 + 5$ | $182x_4 + 9$ |

Even more, he employs a backward substitution, which we can see in the above substitution from $x_3$ upwards. We observe the connection to the diagonals of a heptagon $P_7$ from Section 13, $x = \frac{d}{u}$, and a comparison with a similar substitution table for the nonagon $P_9$ reveals (after suitable reparametrisation) the similarity of backward substitutions

$$\theta^3 - 15\theta^2 + 12\theta - 1 = 0 \quad (P_7) \qquad \theta^3 + 12\theta^2 - 15\theta + 1 = 0 \quad (P_9)$$

## 16. A generalised Euclidean algorithm

*But I will admit to you that I think more often about Euclid's algorithm with n variables, and anti-prime numbers (like $2^\alpha 3^\beta 5^\gamma$ or $2^\alpha 3^\beta 5^\gamma 7^\delta$... etc) of which I gave a small overview at the end of my book "Quaternions". It is fundamental in my mind, or this idea of*

[37] de Casteljau, personal communications to J.-L. Loschutz, 2006, also to H. Prautzsch, 14 December 2010.

*which I had found the germ, at the age of ten, in a book, ended up imposing itself on me, as an example of an unfinished problem.* — de Casteljau[38]

The second part of *Les Quaternions* is dedicated to the Euclidean algorithm, on which de Casteljau already started pondering as a child. He presents a generalisation to multiple variables and deduces a 'machine' to approximate the root of specific cubic polynomials.

### *16.1. Remind Euclid's algorithm and continued fractions*

The Euclidean algorithm is well-known for its computation of the greatest common denominator $GCD$ of two integers. Besides this property, de Casteljau explores the algorithm to make use of continued fractions, which are closely related, and he works out a generative view of the algorithm.

We assume a series of $p+1$ integers $n_0 > n_1 > \ldots > n_p \geq 0$ with

$$n_{i-1} = q_i n_i + n_{i+1}, \qquad q_i = \left\lfloor \frac{n_{i-1}}{n_i} \right\rfloor$$

Euclid's algorithm starts with a first division $n_0/n_1$ and generates the integer series of $q_i$ by repeated divisions, until $i = p$, when the remainder $n_{i+1}$ equals 0 or 1. If $n_{p+1} = 0$, then $n_p$ divides $n_{p-1}$, and as $n_{p-1}$ in the previous step divided $n_{p-2}$, $n_p$ divides all $n_i$, including $n_1$ and $n_0$. The algorithm stops and has found $n_p$ as $GCD(n_0, n_1)$.

If $n_{p+1} = 1$, then both numbers $n_{p-1}$ and $n_p$ (as well as all $n_i$ and thus $n_0$ and $n_1$) do not have a common denominator, they are relatively prime.

The divisions follow the recurrence

$$u_i = \frac{n_{i-1}}{n_i} = q_i + \frac{n_{i+1}}{n_i} = q_i + \frac{1}{n_i/n_{i+1}} = q_i + \frac{1}{u_{i+1}}$$

and define a continued fraction

$$\frac{n_0}{n_1} = q_1 + \cfrac{1}{q_2 + \cfrac{1}{\ddots \quad q_p + \cfrac{n_{p+1}}{n_p}}}$$

The recurrence of the division leads to a recurrence of its components. The observation

$$u_i = q_i + \frac{1}{u_{i+1}} = \frac{q_i \cdot u_{i+1} + 1}{1 \cdot u_{i+1} + 0} \quad =: \quad \frac{S_i u_{i+1} + S_{i-1}}{D_i u_{i+1} + D_{i-1}} = \frac{(S_i q_{i+1} + S_{i-1}) + S_i/u_{i+2}}{(D_i q_{i+1} + D_{i-1}) + D_i/u_{i+2}}$$

leads to the recurrence

$$\begin{cases} S_{i+1} = S_i q_{i+1} + S_{i-1} \\ D_{i+1} = D_i q_{i+1} + D_{i-1} \end{cases}$$

with initial values

$$\begin{cases} S_0 = 1, \quad S_1 = q_1, \quad S_2 = q_1 q_2 + 1, \quad \ldots \\ D_0 = 0, \quad D_1 = 1, \quad D_2 = q_2, \quad \ldots \end{cases}$$

The algorithm can also be stopped after $i$ steps, which truncates $u_i$ to $u_i = q_i$ and defines the $i$-th convergent as

$$\frac{S_i}{D_i} = q_1 + \cfrac{1}{q_2 + \cfrac{1}{\ddots \quad q_i}}$$

It further holds for the determinants, that

$$\begin{vmatrix} S_i & S_{i+1} \\ D_i & D_{i+1} \end{vmatrix} = \begin{vmatrix} S_i & S_i q_{i+1} + S_{i-1} \\ D_i & D_i q_{i+1} + D_{i-1} \end{vmatrix} = \begin{vmatrix} S_i & S_{i-1} \\ D_i & D_{i-1} \end{vmatrix} = - \begin{vmatrix} S_{i-1} & S_i \\ D_{i-1} & D_i \end{vmatrix} = \cdots = (-1)^i \begin{vmatrix} S_0 & S_1 \\ D_0 & D_1 \end{vmatrix} = (-1)^i,$$

which we apply to the identity

[38] de Casteljau, personal communication (translated from French), 10 November 1999.

**Table 3**
Euclidean algorithm for 99/70 and first convergents of $\sqrt{2}$.

| $i$ | -1 | 0 | 1 | 2 | | ⋯ | | $p$ | | ⋯ |
|---|---|---|---|---|---|---|---|---|---|---|
| $n_i$ | | 99 | 70 | 29 | 12 | 5 | 2 | 1 | | |
| $q_i$ | | | 1 | 2 | 2 | 2 | 2 | 2 | 2 | ⋯ |
| $S_i$ | 0 | 1 | $q_1$ | 3 | 7 | 17 | 41 | 99 | 239 | ⋯ |
| $D_i$ | 1 | 0 | 1 | $q_2$ | 5 | 12 | 29 | 70 | 169 | ⋯ |
| $D_i D_{i+1}$ | | | 2 | 10 | 60 | 348 | 2 030 | 11 830 | 68 952 | ⋯ |

$$\frac{S_p}{D_p} = \frac{S_p}{D_p} - \frac{S_{p-1}}{D_{p-1}} + \frac{S_{p-1}}{D_{p-1}} \mp \cdots + \frac{S_2}{D_2} - \frac{S_1}{D_1} + \frac{S_1}{D_1} = \frac{\begin{vmatrix} S_p & S_{p-1} \\ D_p & D_{p-1} \end{vmatrix}}{D_p D_{p-1}} + \cdots = q_1 + \frac{1}{D_1 D_2} \mp \cdots + \frac{(-1)^p}{D_{p-1} D_p} = \frac{n_0}{n_1}.$$

This shows that $\frac{n_0}{n_1} = \frac{S_p}{D_p}$ is approximated from above and from below by the alternating sum – which can easily be visualised by lattice vectors:

$$\mathbf{v}_{i+1} = q_{i+1}\mathbf{v}_i + \mathbf{v}_{i-1}, \qquad \text{with } \mathbf{v}_i = \begin{bmatrix} S_i \\ D_i \end{bmatrix}$$

We have explored the Euclidean algorithm as it divides $n_0$ by $n_1$ in a finite number $p$ of steps, resulting in a rational number. In the infinite case, if the quotients $q_i$ exhibit periodicity, the continued fraction represents an irrational number; otherwise it represents a transcendent number. Both infinite cases ask for a generative principle for the $q_i$, since we cannot deduce them from the division of two integers. For example, in the case of

$$\sqrt{2} = 1 + \cfrac{1}{2 + \cfrac{1}{2 + \cdots}}$$

we have $q_1 = 1$ and $q_i = 2$ for $i < 1$, as shown in Table 3.

### Addenda

- Both perspectives – the deductive Euclidean algorithm and the inductive generative principle – relate to classical Diophantine problems, as de Casteljau explains (1987): For given values $A, B, C, D, \ldots$, find integers $a, b, c, d, \ldots$, such that

$$\frac{a}{A} \approx \frac{b}{B} \approx \frac{c}{C} \approx \frac{d}{D} \approx \cdots \qquad \text{(Hermite)}$$

$$a\,A + b\,B + c\,C + d\,D + \cdots \approx 0 \qquad \text{(Dirichlet)}$$

  The connectedness of both approaches can also be seen using metric geometry, cf. Fig. 29.
- The $i$-th convergent of a continued fraction $n_0/n_1$ reads in matrix form:

$$\begin{bmatrix} S_i \\ D_i \end{bmatrix} = \begin{bmatrix} q_1 & 1 \\ 1 & 0 \end{bmatrix} \begin{bmatrix} q_2 & 1 \\ 1 & 0 \end{bmatrix} \cdots \begin{bmatrix} q_{i-1} & 1 \\ 1 & 0 \end{bmatrix} \begin{bmatrix} q_i \\ 1 \end{bmatrix}$$

### *16.2. Multiple variables*

The generalisation to more than two variables considers both perspectives: the deductive approach from the Euclidean algorithm and the generative approach.

De Casteljau reminds us that the division in Euclid's algorithm is nothing other than a sequence of differences. Therefore, an approach with minimal subtractions always focuses on the two largest numbers instead of, e.g., the largest and the smallest. His algorithm for $p$ terms $a_i$ reads (1987):

1 We sort the $a_i$ in decreasing order, which results in a permutation of the indices $a_i$.
2 From the largest number we subtract the second, as often as necessary.
3 We regroup the remainder to its new position, the other values stay unchanged, which provokes a new permutation.
4 In case of rigorous equality, we receive an exact algebraic Dirichletian relation and we can return to the case of $p-1$ elements, or proceed with a zero element.
5 If the final remainder is 1, the latest algebraic relation generalises Bezout's theorem and the numbers are all relative prime; otherwise we obtain the $GCD$ and the latest relation with an arbitrary second member.
6 All these relations are independent.
7 In the case of linked algebraic numbers, $p$ real roots of an equation of $p$-th degree, or 0- to $p-1$-fold real roots, a periodicity of operations can appear.

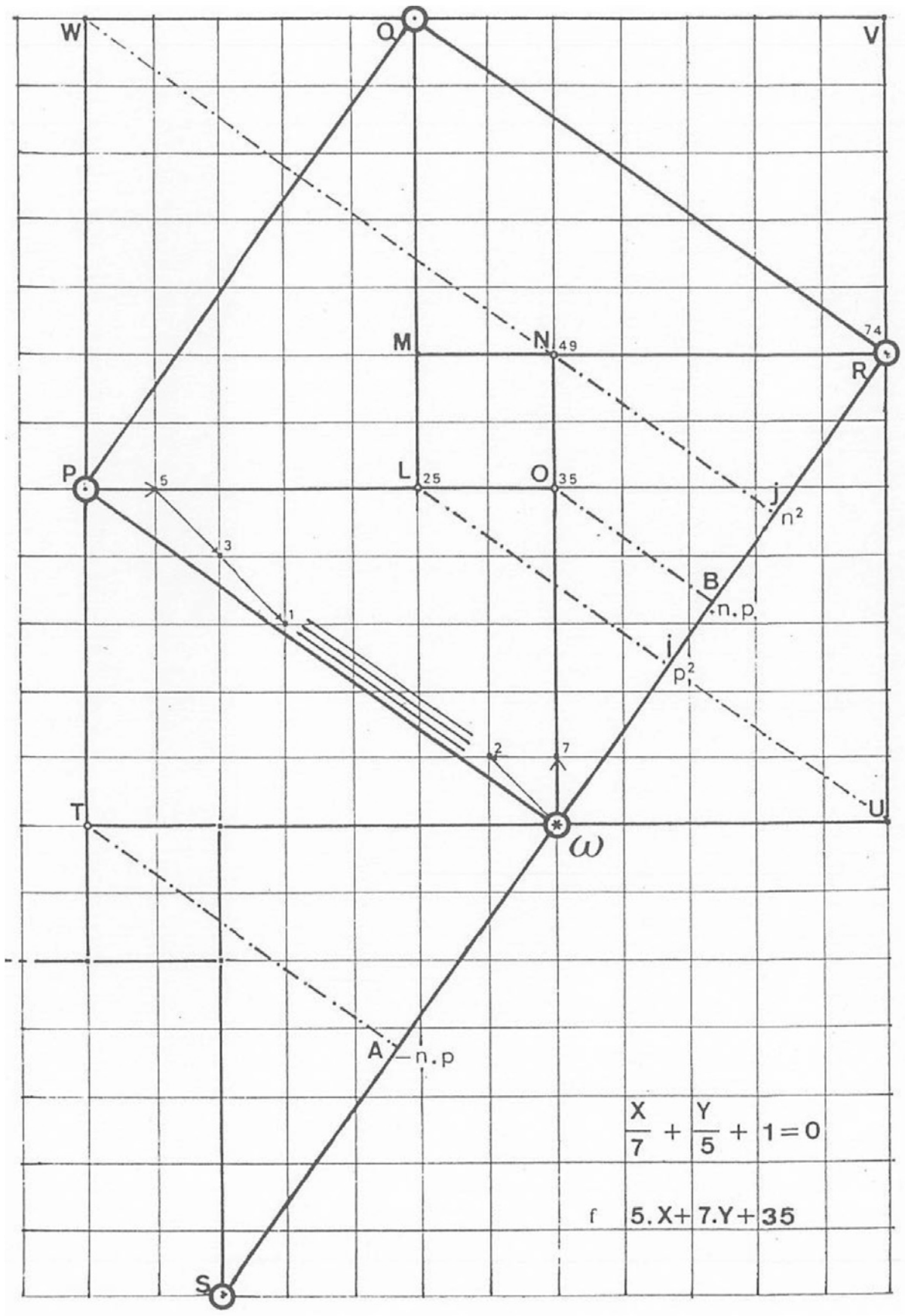

**Fig. 29.** A metric Euclidean algorithm, drawing by de Faget de Casteljau (1995b).

The generative perspective works inversely. The distance of an approximate $a_i$ to a targeted $x_i$ is evaluated against a threshold for all $a_i$, before $i$ is increased. For a deeper exploration and a practical illustration, the interested reader can refer to (de Casteljau, 1987, ch.12), where we also find an application with logarithms of small integers and further insights on the Hermitian approximation result, related to a musical scale with 171 tones (including the dissonant tritone frequency ratio 7:5)[39]:

$$\frac{\ln 2}{171} \approx \frac{\ln 3}{271} \approx \frac{\ln 5}{397} \approx \frac{\ln 7}{480}$$

### 16.3. The heptagon revisited

Let us recall the heptagon from section 13 and the relations (10) between their diagonals $u, d, t$. De Casteljau provides us with a decreasing scale of linear combinations of $u, d, t$ which simultaneously are products (1987, p.157), e.g.,

$$\begin{aligned}
s_0 &= t^3/u^2 = 3t + 2d + u \\
s_{14} &= t^1 = t \\
s_{16} &= d^1 = d \\
s_{21} &= u^1 = u \\
s_{19} &= tu/d = t - u \\
s_{23} &= du/t = d - u \\
s_{28} &= u^2/t = t - d \\
s_{30} &= du^2/t^2 = 2d - u - t = s_{23} - s_{28}
\end{aligned}$$

39 In a letter to Boehm, 22 January 2010, de Casteljau remarks, that a solution spans 80-81 octaves: $7^{28} < 5^{34} < 2^{79} < 3^{50} < 2^{80} < 3^{51} < 2^{81} < 5^{35} < 7^{29}$.

**Table 4**
A double entry table to calculate the values $u, d, t$.

| | $c$ | $b$ | $a$ | $b_{+c}$ | $a_{+b}$ | $b_{+a}$ | $c_{+b}$ | $b_{+c}$ |
|---|---|---|---|---|---|---|---|---|
| $t$ | 1 | 0 | 0 | 1 | **1** | **2** | **3** | ⋮ |
| $d$ | 0 | 1 | 0 | 1 | **1** | **2** | **2** | |
| $u$ | 0 | 0 | 1 | 0 | **1** | **1** | **1** | |
| $t_{-d}$ | **1** | **-1** | **0** | 0 | 0 | 0 | 1 | |
| $d_{-u}$ | 0 | 1 | -1 | 1 | 0 | 1 | 1 | |
| $u_{-d}$ | **0** | **-1** | **2** | -1 | 1 | 0 | 0 | |
| $d_{-t}$ | **-1** | **2** | **-1** | 1 | 0 | 1 | 0 | |
| $t_{-d}$ | … | | | | | | | |

$$s_{35} = u^3/t^2 = 2u - d = s_{21} - s_{23}$$

$$\vdots$$

Using the decreasing order of $t = s_{14}, d = s_{16}, u = s_{21}$, we can now perform the generalised Euclidean algorithm and find:

$$(s_{14}, s_{16}, s_{21}): \quad t_1 = t - d = u^2/t = s_{28}$$

$$(s_{16}, s_{21}, s_{28}): \quad d^* = d - u = ud/t = s_{23}$$

$$(s_{21}, s_{23}, s_{28}): \quad u_1 = u - d^* = u^3/t^2 = s_{35}$$

$$(s_{23}, s_{28}, s_{35}): \quad d_1 = d^* - t_1 = u/t(d-u) = u^2 d/t^2 = s_{30}$$

The 'new' values $t_1, d_1, u_1$ for the next step are proportional to $t, d, u$ by a factor $u^2/t^2$:

$$\frac{t_1}{t} = \frac{d_1}{d} = \frac{u_1}{u} = \frac{u^2}{t^2}$$

Regrouping $s_{28}$, de Casteljau observes

$$d = t - \frac{u^2}{t} = t\left(1 - \frac{u^2}{t^2}\right)$$

which allows us to write the elements related to one value only, $x = u/t$:

$$\begin{array}{cccc} t & & u = xt & d = (1-x^2)t \\ t_1 = x^2 t & d_1^* = x(1-x^2)t & u_1 = x^3 t & d_1 = x^2(1-x^2)t \\ t_2 = x^4 t & d_2^* = x^3(1-x^2)t & u_2 = x^5 t & d_2 = x^4(1-x^2)t \cdots \end{array}$$

This idea can now be taken into a storage table, which takes into account the periodicity of the (generalised) Euclidean algorithm. The columns follow the relations:

$$b_1^* = b + c; \qquad a_1 = a + b_1^*; \qquad b_1 = b_1^* + a_1; \qquad c_1 = c + b_1;$$

or in compact pseudocode

$$b \leftarrow b + c; \qquad a \leftarrow a + b; \qquad b \leftarrow b + a; \qquad c \leftarrow c + b;$$

We see the result in Table 4, which leads to a $4 \times 4$-matrix, after an initialisation of the three variables. The table can be understood as the status of stored variables $(u, d, t)$ while four operations are repeated $(+a, +b, +c, +b)$.

We recognise the inverse matrices

$$M^{-1} = \begin{bmatrix} t_1 \\ d_1 \\ u_1 \end{bmatrix} = \begin{bmatrix} 1 & -1 & 0 \\ -1 & 2 & 1 \\ 0 & 1 & 2 \end{bmatrix} \begin{bmatrix} t \\ d \\ u \end{bmatrix} \quad \text{and} \quad M^1 = \begin{bmatrix} c_1 \\ b_1 \\ a_1 \end{bmatrix} = \begin{bmatrix} 3 & 2 & 1 \\ 2 & 2 & 1 \\ 1 & 1 & 1 \end{bmatrix} \begin{bmatrix} c \\ b \\ a \end{bmatrix}$$

and their similarity with the ones from section 13.

### 16.4. *"Je me gausse..."*

In a one-pager written at Pentecost 2007, which is reproduced in Fig. 28, de Casteljau shows a structure of the matrix $M$ for the case of the heptagon:

$$M^k = \begin{bmatrix} R_{k-1} & V_{k-1/2} & R_k \\ V_{k-1/2} & R_{k-1} + R_k & V_{k+1/2} \\ R_k & V_{k+1/2} & R_{k+1} \end{bmatrix}$$

which could be associated with the understanding of the beam of light, which is reflected inside the circle to generate a real image ($R$), or outside the circle to generate a virtual one ($V$).

The asymptote is described by

$$\begin{bmatrix} V_{k-3/2}+R_k \\ R_{k-2}+R_{k+1} \\ V_{k-1/2}+R_{k+1} \end{bmatrix} \quad \text{and its limit comes by} \quad M^{\infty}=\begin{bmatrix} u \\ d \\ t \end{bmatrix}$$

The recursive calculation is based on

$$V_{k+1/2}-V_{k-1/2}=R_k$$
$$R_{k+1}-R_{k-1}=V_{k+1/2}$$

These integer series are known as A006054 for $R_i$, A006356 for $V_{i+1/2}$, and A187070 for the interwoven sequence of both $R_i, V_{i+1/2}$.

### 16.5. Hermite meets Dirichlet

By making repetitive use of the Euclidean algorithm, de Casteljau constructs the following infinite matrix $DH$ based on the noted submatrices, its initialisation could be seen in Table 4. Its symmetry would also make it possible to store only the $4 \times n$ stripe $\ldots D_4 D_3 D_2 D_1 H_0 H_1 H_2 H_3 H_4 \ldots$, which represents a bi-directional use of the generalised Euclidean algorithm.

$$DH=\begin{bmatrix} H_0 & H_1 & H_2 & H_3 & H_4 & \cdots \\ D_1 & H_0 & H_1 & H_2 & H_3 & \\ D_2 & D_1 & H_0 & H_1 & H_2 & \\ D_3 & D_2 & D_1 & H_0 & H_1 & \\ D_4 & D_3 & D_2 & D_1 & H_0 & \\ \vdots & & & & \ddots & \end{bmatrix}$$

$$H_0=\begin{bmatrix} 0&0&0&1\\ 1&0&1&1\\ -1&1&0&0\\ 1&0&1&0 \end{bmatrix} \quad H_1=\begin{bmatrix} 1&1&2&3\\ 2&2&4&4\\ 0&1&1&1\\ 1&1&2&2 \end{bmatrix} \quad H_2=\begin{bmatrix} 5&6&11&14\\ 9&11&20&25\\ 2&3&5&6\\ 4&5&9&11 \end{bmatrix} \quad H_3=\begin{bmatrix} 25&\mathbf{31}&\mathbf{56}&\mathbf{70}\\ 45&56&101&126\\ 11&\mathbf{14}&\mathbf{25}&\mathbf{31}\\ 20&\mathbf{25}&\mathbf{45}&\mathbf{56} \end{bmatrix}$$

$$D_0=H_0 \quad D_1=\begin{bmatrix} -1&0&-1&1\\ 2&-1&1&0\\ -3&2&-1&0\\ 3&-1&2&-1 \end{bmatrix} \quad D_2=\begin{bmatrix} -4&1&-3&2\\ 6&-3&3&-1\\ -9&5&-4&1\\ 10&-4&6&-3 \end{bmatrix} \quad D_3=\begin{bmatrix} -14&\mathbf{5}&\mathbf{-9}&\mathbf{5}\\ 19&-9&10&-4\\ -28&\mathbf{14}&\mathbf{-14}&\mathbf{5}\\ 33&\mathbf{-14}&\mathbf{19}&\mathbf{-9} \end{bmatrix}$$

In general, the recurrence relations can be described line by line ($n \geq 0$):

$$\begin{aligned} d_{4n+4}&=d_{4n}-d_{4n+3} & \cdots\\ d_{4n+5}&=d_{4n+3}-d_{4n+2} & \cdots\\ d_{4n+6}&=d_{4n+2}-d_{4n+5} & \cdots\\ d_{4n+7}&=d_{4n+5}-d_{4n+4} & \cdots \end{aligned}$$

as well as column by column

$$\begin{matrix} h_{4n+4}=h_{4n+2}+h_{4n+3} & h_{4n+5}=h_{4n+1}+h_{4n+4} & h_{4n+6}=h_{4n+4}+h_{4n+5} & h_{4n+7}=h_{4n+3}+h_{4n+6} \\ \vdots & \vdots & \vdots & \vdots \end{matrix}$$

### 16.6. Characteristic equations and periodicity

The eigenvalues of $M^{-1}$ and $M^{1}$ give the value of the root(s) $x$

$$M\mathbf{h}=x\mathbf{h}$$

which allows to apply $\mathbf{h}=[t,d,u]$ on itself. In particular, the traces $tr(M^{-1})=5$ and $tr(M)=6$ deliver the coefficients of the characteristic equation (here done with Sarrus' rule):

$$det(M-xI)=(3-x)(2-x)(1-x)+2\cdot 2-(2-x)-(3-x)-4(1-x)=-x^3+6x^2-5x+1$$

De Casteljau now makes use of Vincent's algorithm: The substitution $x=1/d^2$ turns the characteristic equation

$$1-5x+6x^2-x^3=0 \qquad \text{into} \qquad d^3+d^2-2d-1=0.$$

From the examples in section 15, we know the roots

$$x_1 = 5 + \frac{1}{t} \qquad \text{and} \qquad d_1 = 1 + \frac{1}{4 + \frac{1}{t}}$$

with

$$t = \frac{9t+1}{t(t-20)} \qquad \text{or} \qquad 61t = 1247 + \frac{1}{61t+13} + \frac{1}{61t+14}$$

De Casteljau proposes also other line substitutions, such as $x = d^3$ with $x^3 + 4x^2 - 8 = 0$, but also column substitutions, such as $y = h^2$, which transforms $h^3 - 3h^2 - 4h - 1$ into $y^3 - 17y^3 + 10y - 1 = 0$. Those substitutions as well as the following cyclic patterns invite further research.

The matrices $M^{\pm p}$ are generated by a pattern of additions or subtractions. Identifying the variables $t, d, u$ and $a, b, c$ both with $ABC$, de Casteljau denotes this pattern $ABCB$, which is short for $1A \pm 1B \pm 1C \pm 1B \pm 1A \pm 1B \ldots$. In chapter 13.3. of *Les Quaternions*, de Casteljau lists various cubics and their characteristic equations, which he sees under the aspect of color coded lengths in decreasing order like $ABCB$ for the periodicity of the diagonal lengths in a heptagon. The beginning of the list is replicated here, where $ApB$ reads: 'subtract $B$ $p$ times from $A$':

$$\begin{array}{rrcl}
ABC & x^3 - 4x^2 + 3x - 1 & = & 0 \\
A2BC & x^3 - 5x^2 + 3x - 1 & = & 0 \\
ABCB & \mathbf{x^3 - 6x^2 + 5x - 1} & = & 0 \text{ (heptagon)} \\
A3BC & x^3 - 6x^2 + 3x - 1 & = & 0 \\
A2B2C & x^3 - 7x^2 + 3x - 1 & = & 0 \\
A2BCB & x^3 - 8x^2 + 6x - 1 & = & 0 \\
ABABC & x^3 - 8x^2 + 5x - 1 & = & 0 \\
A4BC & x^3 - 7x^2 + 3x - 1 & = & 0 \\
\ldots & & &
\end{array}$$

## 17. On the sum of three cubes

*In research you must be versatile! It's all about algebra, analysis, calculus of variations, and all the "geometries"... The Queen of mathematics remains number theory.* — de Casteljau[40]

The story goes that de Casteljau was asked by some colleagues about an approach to proving Fermat's Last Theorem. Unfortunately, his notes about this topic were lost. Despite his mastery of 19th century mathematics – a period marked by significant progress on this question by Sophie Germain and others – he indicated, that we should not expect to find a proof from him.

### *17.1. Ménard's identity*

When de Casteljau worked on the Diophantine equation $X^3 + Y^3 + Z^3 = K$, his particular interest was on $K$ being a cube itself, $K = T^3, |T| > |X| > |Y|$, one of the still unsolved problems in number theory. In letters to Boehm,[41] he links his findings to results reached by Euler and Ramanujan, who (according to de Casteljau) both might have overseen this approach. He starts with the equation

$$(3^6 + 3^2)^3 - (3^5 + 1)^3 = (3^6)^3 - 1^3 = (3^6 - 3^2)^3 + (3^5 - 1)^3. \tag{15}$$

This triple identity is denoted by de Casteljau as the *Ménard identity* named in honour of the Citroën computer scientist Michel Ménard, who directed de Casteljau towards this discovery.

The generalisation

$$(3^{4n-2} + 3^n)^3 - (3^{3n-1} + 1)^3 = (3^{4n-2})^3 - 1 = (3^{4n-2} - 3^n)^3 + (3^{3n-1} - 1)^3$$

quickly leads to large cubes ($n = 1, 2, 3, 4$):

$$\begin{array}{ccc}
12^3 - 10^3 & = 9^3 - 1 = & 6^3 + 8^3 \\
738^3 - 244^3 & = 729^3 - 1 = & 720^3 + 242^3 \\
59076^3 - 6562^3 & = 59049^3 - 1 = & 59022^3 + 6560^3 \\
4783212^3 - 177148^3 & = 4782969^3 - 1 = & 4782726^3 + 177146^3
\end{array}$$

A slightly different notation leads to the above-mentioned sum of three cubes:

[40] de Casteljau, personal communication to H. Prautzsch (translated from French), 14 December 2010.

[41] de Casteljau, personal communication to W. Boehm, December 2008, January 2010 and December 2010.

$$(3^{4n-2}+3^n)^3=(3^{4n-2}-3^n)^3+(3^{3n-1}-1)^3+(3^{3n-1}+1)^3 \tag{16}$$

which for $n=1,2,3,4$ reads:

$$12^3=10^3+8^3+6^3$$
$$738^3=720^3+244^3+242^3$$
$$59076^3=59022^3+6562^3+6560^3$$
$$4783212^3=4782726^3+177148^3+177146^3$$

De Casteljau leaves us the general equation

$$(Z+Y)^3=(Z-Y)^3+(X+1)^3+(X-1)^3$$

which invites further studies.

Obviously this identifies with equation (15) for $X=243, Y=9, Z=729$. A special case is $X=3^m, Y=3^n, Z=3^{m+1}$ which results in the equation $(3^{m+1}+3^n)^3=(3^{m+1}-3^n)^3+(3^m+1)^3+(3^m-1)^3$ with its two solutions $(m,n)=(-1,-1)$ or $(5,2)$.

**Addenda**

- Regarding $(Z+Y)^3-(Z-Y)^3=(X+1)^3+(X-1)^3$ led de Casteljau to $X^3+3X-Y^3-3YZ^2$ with solutions

$$Z^2=\frac{X^3+3X-Y^3}{3Y}$$

- In his late studies on anti-prime (or highly composite) numbers, e.g., 7-smooth numbers $2^\alpha 3^\beta 5^\gamma 7^\delta$, de Casteljau indicated an extension of Catalan's conjecture *(*$X^a=Y^b-1$ *only holds for* $3^2$ *and* $2^3$*)*, which he connects to the sparse matrice products $M^pM^{-q}$ of the above-mentioned generalised Euclidean algorithm.[42]

### *17.2. A link to Euler and Ramanujan*

Let us now look specifically at the Ménard identity with $n=1$, which combines two views on cubes:

$$12^3-10^3=9^3-1=6^3+8^3 \Leftrightarrow A=B=C$$

**A = C:** De Casteljau identifies Euler's finding $3^3+4^3+5^3=6^3$, here constructed by a Ménard identity with halved values (or divided by 8):

$$6^3-5^3=\left(\frac{9+3}{2}\right)^3-\left(\frac{9+1}{2}\right)^3=\left(\frac{9}{2}\right)^3-\left(\frac{1}{2}\right)^3=\left(\frac{9-3}{2}\right)^3-\left(\frac{9-1}{2}\right)^3=3^3+4^3.$$

**A = B:** A quick rearrangement leads him to Ramanujan's famous *taxicab* number $1729=1^3+12^3=9^3+10^3$, which is the smallest number that can be written as two different sums of cubes. Differently to Euler, he notes the more general quadratic forms

$$\begin{array}{ccccccc}(3X^2+5XY-5Y^2)^3 & + & (4X^2-4XY+6Y^2)^3 & + & (5X^2-5XY-3Y^2)^3 & = & (6X^2-4XY+4Y^2)^3\\ (X^2+9XY-Y^2)^3 & + & (12X^2-4XY+2Y^2)^3 & = & (9X^2-7XY-Y^2)^3 & + & (10X^2+2Y^2)^3,\end{array}$$

which obviously lead to the initial identity (15) with $X=1, Y=0$.

De Casteljau displays his understanding of Euler's identity $X^3+Y^3+Z^3=T^3$ by introducing the fraction $Q=\dfrac{X+Y}{T-Z}$, which leads him to

$$2[XZ+TY]=\frac{Q[Y^2-X^2]+T^2-Z^2}{u}=\frac{QXY+TZ}{v}$$

and $Q=u^2+3v^2$, thus

$$\begin{aligned}&\frac{X}{9v^3-3v^2u+3u^2v-u^3+1}=\frac{Y}{9v^3+3v^2u+3u^2v+u^3-1}\\ &=\frac{Z}{u^4+6u^2v^2+9v^4-u-3v}=\frac{-T}{-u^4-6u^2v^2-9v^4+u-3v},\end{aligned}$$

where he highlights the terms $\pm(3v^2u+u^3-1)$ and $\pm(u^4+6u^2v^2+9v^4-u)$.

He further notes that Ramanujan parametrises planar sections that rest on the edges of the tetrahedron $X+Y=0=Z-T$ or similar.

---

42 de Casteljau, personal communication, 10 November 1999.

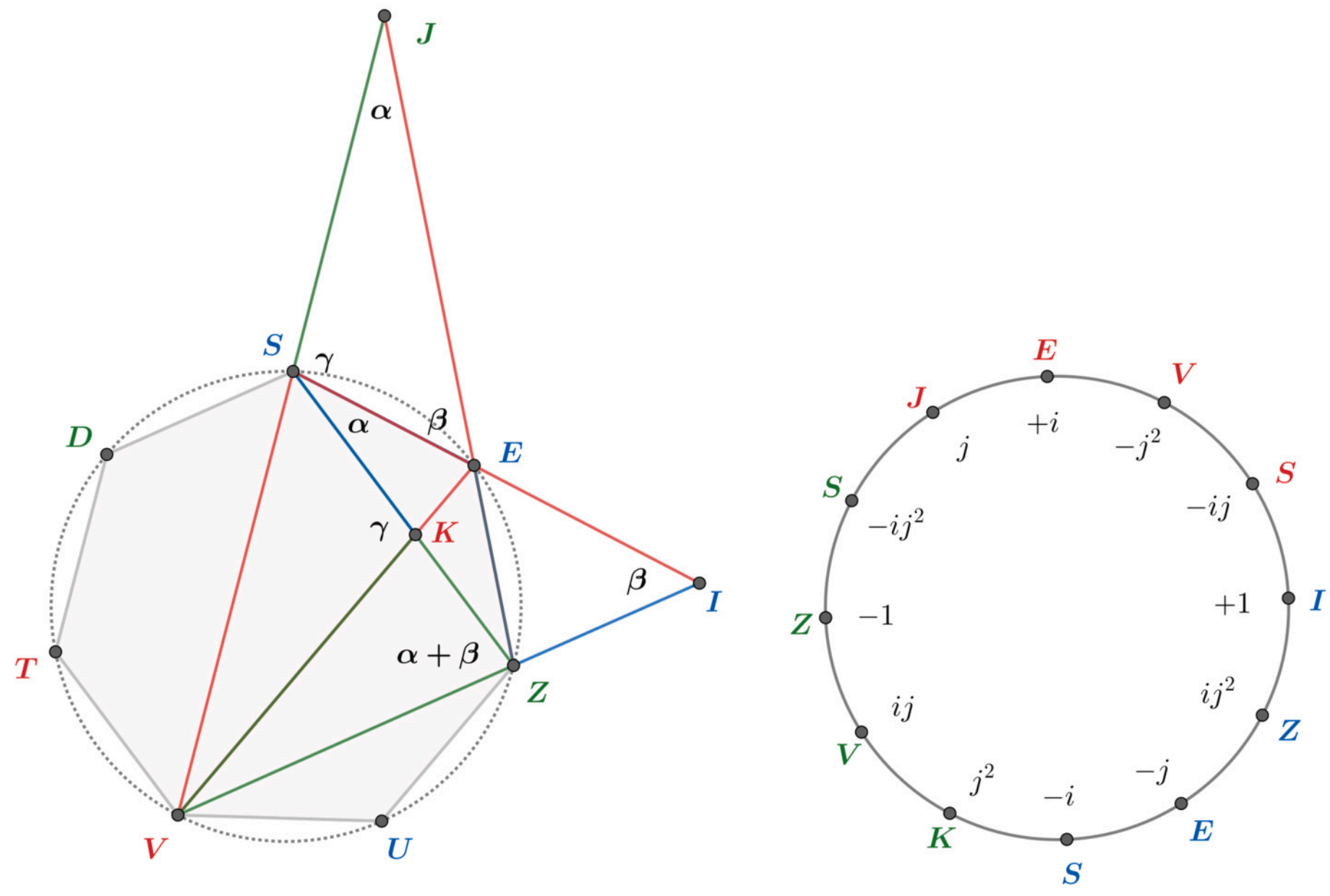

**Fig. 30.** A complete quadrilateral in a heptagon with composing lengths and angles.

De Casteljau hints at a conjecture, that for all prime numbers $6k + 1$, it is possible to define a "cubic residual", similar to the quadratic residual of Legendre. The thoughts in the mentioned letters then walk through the remainders modulo 27, then mod $7, 13, 19, 31, 37, 109$, and leave room for own ideas.

## 18. The heptagon again

$P_7$ *Wunderschön* — de Casteljau[43]

We would like to leave the reader with some ideas on the heptagon, which long fascinated de Casteljau. He denoted the summit by $S$ and took the initial letters of French and German numbers for the remaining six points: *Un (one), Deux (two), Trois (three), Eins (one), Zwei (two), Vier (four)* (remember the symmetry $q \equiv n - q$).

A complete quadrilateral is constructed of four heptagon points (e.g., $S, E, Z, V$) and three further points $I, J, K$; the vertices of the quadrilateral equal the diagonal lengths $u, d, t$ of the heptagon, cf. Fig. 30. We can make a round trip via all vertices in the order of consecutive lengths, $u : IZESK, d : KVZSJ, t : JEVSI$. De Casteljau organised the angles via $i, j$ on this trip around the 'clock', which progress by a multiplication of $ij^2$ clockwise or $-ij$ anti-clockwise, assuming $i^4 = 1$ and $j^3 = 1$. A comparison with the angles $\alpha = \angle SJE, \beta = \angle JES, \gamma = \angle ESJ$ by use of similarities

| | | |
|---|---|---|
| $\alpha$ | $-j^2(EVS), -i(ESK), ij(KVZ), j(SJE)$ | mod $j$ |
| $\beta$ | $+1(SIZ)$ | mod $j^3$ |
| $\gamma$ | $j^2(SKV)$ | mod $j^2$ |
| $\alpha + \beta$ | $-ij(VSI), ij^2(IZE), -1(VZS)$ | mod $i^2$ |
| $\alpha + \gamma$ | $+i(JEV), -j(ZES), -ij^2(ZSJ)$ | mod $i$ |

reveals a discrete version of Fig. 18 (left), where $I, J, K$ are cutting secants out of the circle.

The image also lets de Casteljau indicate the geometrical distinction of the roots $+i$ and $-i$ by their angle, different to the less distinctive algebraic view of $(\pm i)^2 = -1$.

## Conclusion

## 19. Achievements

*Rest assured, I am in no way expecting to be immortalised, in the form of an equestrian statue, with a laurel wreath on my head, neither in front of my birthplace, nor in the courtyard of the Ecole Normale Supérieure, the Academy of Sciences, or even in Oberwolfach... It*

[43] de Casteljau, personal communication to W. Boehm, May 2007. He used the German word *wunderschön* for 'marvellous' when referring to the heptagon in his one-pager *Je me gausse de GAUSS*, cf. Fig. 28.

**Table 5**
De Casteljau's four principles for working with curves and surfaces.

| | |
|---|---|
| **1st** | To avoid the infinite and cyclical points, ALL intermediary construction points must be on the piece of paper, which condemns rational splines (unless weights are composed of squares, therefore positive) |
| **2nd** | The tool works only on one point at a time. That makes the notion of a parameter seem logical, as you do not have to choose between multiple solutions. |
| **3rd** | Only the real domain is accessible. All passage through complex concepts remains dangerous. |
| **4th** | All too often, the physical nature of mathematisation is forgotten, especially the concept of 'error calculation'. It is advisable to avoid indeterminate forms and leave them to Calculus, e.g., the differences of two very large numbers or the quotients of two small quantities. These considerations are well disguised when it comes to matrix calculus or least squares, which have the appearance of confidence-inspiring accuracy while approaching "overflows" or "underflows". I can't warn enough against giving in to the illusion of getting everything right while acrobatically teetering on the edges of the abyss. |

*amuses me greatly to give a presentation in front of an audience of scholars, the least qualified of whom has had his doctorate for three years, while officially they are listening to a modest "licencié ès Sciences" [Science graduate, which was he].* — de Casteljau[44]

We have explored a variety of thoughts that Paul de Faget de Casteljau brought forward during his life, the vast majority after finding his algorithm in late 1958. He demonstrates a broad range of ideas, based on a solid mathematical foundation, always presents calculable approaches and links views from within different geometries (affine, projective, metric) with those from algebra or number theory, often connected to applications as in geometric optics. His preference for integers turned out to be an effective way to maintain highest precision, similar to him refraining from rational curves in order to avoid singularities such as infinity.

### *19.1. A summary of summarised memory aids*

Paul de Faget de Casteljau has left a distinct mark on Computer-Aided Design. His fundamental algorithm is in use in many places, as Hartmut Prautzsch (2023) has recently impressively demonstrated. Based on barycentric calculus (section 2), the algorithm is numerically stable, and he applied it to polynomial curves and surfaces for automotive design as well as to periodic functions as in cam design (section 3). When it was needed, he generalised his original approach to integrate continuity (blossoming, section 4) and restitution (algebraic smoothness, section 5), the latter going unnoticed for a long time.

He was there at the cradle of digitalisation in the automotive sector and had to defy some resistance at Citroën, where the mathematical genius met the pragmatism of the engineers. His talent in working out simple solutions to tricky problems earned him the title of *'Maître'* that his colleagues awarded him. We could see corresponding examples: trigonometric smoothing (section 6), the forward-looking calculations of tolerances (section 7), intersections as an arithmetic polar form (section 8), and even more examples can be found in his book *Le Lissage*.

While his interests increasingly drifted away from affine geometry (which he saw overloaded by projective geometry and its notion of the degree), he turned back once to honour a study by his late brother Henri. The result is a compact visualisation of six classic theorems by 16 lines and three conics $(16D + 3C)$ as well as a combinatorial access to those theorems (section 9).

We had a glimpse into his thoughts on metric geometry (section 10), some of which received some publicity through his article on focal splines. Yet, his findings on Apollonian circles, or his construction of a 14-point strophoid by biangular coordinates (section 11), or on conjugate mirrors in geometric optics (section 12) were barely known or discussed.

Similarly, his finding of an approximation to construct a regular polygon (section 13) did escape the mathematical community.

The most prominent algebraic work is his book *Les Quaternions*, which introduces a matrix view on quaternions, but even more proposes quaternions as a frame of orthonormal axes in $E_4$, as could be seen by the example of the Lorentz equations (section 14). Finally, using fewer steps to calculate a rotation, compared to Euler angles, has its advantages as well.

That same book also uses Vincent's lesser known root-finding method (section 15) and covers a generalisation of Euclid's algorithm (section 16) linking it to the polygon diagonals mentioned previously. He introduces a recursive algorithm to approximate the roots of various polynomials, by simple matrix operations and viewed this quasi-symmetric sequences as a sort of geometric computer with coloured strings.

A deviation on the sum of three cubes (section 17) allowed us to follow the breadth of his ideas, before we returned to one of his latest communications on the heptagon (section 18). We departed from there with a view at a complete quadrilateral in a heptagon and its angles.

### *19.2. In search of an essence*

De Casteljau's œuvre is scarcely accessible: it is mainly not digitised, nor translated into English, and particularly his late works are only accessible in fragments or by (unpublished) personal communication, in French. His work surely deserves to be analysed in more depth.

---

[44] de Casteljau, personal communication (translated from French), 12 October 1995.

Besides his CAD oriented modelling principles[45] (cf. Table 5), de Casteljau's work seems to be motivated by two fundamental quests, which shine through all of his publications. Both are not particularly straightforward to reach or implement:

- Finding a recurrence relation, initialising it properly, and approximating as necessary.
- Exploring a polar form (such as with polynomials, trigonometry, analysis, circle, optics, quaternions).

De Casteljau was often successful in finding an approach that satisfied both. Even more, he could utilise his mastery in the field of continued fractions to calculate approximations with highest precision, which also might explain his interest in Diophantine problems.

This tour d'horizon has been intended to whet the appetite of the interested reader. It cannot be exhaustive, but hopefully it has only omitted minor viewpoints of this fairly complete mathematician and physicist. De Casteljau's publications and ideas invite further, more intensive exploration. It is surely worth the time as you might find more approaches and links than presented here.

While his works were not primarily written to please a less open-minded or hurried reader, they are a blossoming bouquet of ideas that thrive on a variety of mathematical perspectives. Paul de Faget de Casteljau has surely not received all the attention that his work deserves – might his approaches flourish posthumously. Let us conclude with a final quote (de Faget de Casteljau, 1986, p.13):

> *My greatest satisfaction will be to see this presentation inspire new works, complements, or applications, regardless of their field of application.*

## CRediT authorship contribution statement

**Andreas Müller:** Writing – review & editing, Writing – original draft, Visualization, Methodology, Investigation, Conceptualization.

## Declaration of competing interest

The authors declare that they have no known competing financial interests or personal relationships that could have appeared to influence the work reported in this paper.

## Acknowledgements

This article owes its existence to Wolfgang Boehm, who ignited my interest in CAGD, and to Paul de Casteljau, who graciously welcomed me into his home in Versailles and generously shared his insightful ideas with me through personal correspondence.

Sincere thanks also go to Jean-Luc Loschutz and Christophe Rabut for their intensive dialogue on the entire work as well as to Hartmut Prautzsch, Alkis Akritas, Claude Paraponaris, Alvy Ray Smith, Jim Heppelmann, Yvon Gardan, Jannick Rolland, Andrew Rakich, Gérard Lemaître, Christoph Baumgarten, Tom Falter, Johannes Wallner, Rida Farouki, Carl de Boor, Gary Walsh for their feedback on individual sections, to Julia Teschner for her thorough proofreading, to the de Casteljau family for their support, and to the anonymous reviewers for their valuable comments and remarks.

## References

### *The works of Paul de Faget de Casteljau*

[45] de Casteljau, personal communication to J.-L. Loschutz (translated from French), early 2006.

## *Other references*